\documentclass[11pt,leqno]{article}
\normalfont
\renewcommand{\baselinestretch}{1.25}

\usepackage{graphicx}
\usepackage{amsfonts}
\usepackage{mathrsfs}
\usepackage{rotating}
\usepackage{mathdots}
\usepackage{hieroglf}
\usepackage{MnSymbol}
\usepackage[dvipsnames]{xcolor}

\newcommand{\IntroSection}{1}
\newcommand{\Setup}{2}
\newcommand{\ESection}{3}
\newcommand{\ConstructionSection}{4}
\newcommand{\CombinatorialDiscussion}{5}

\newcommand{\EsixEsevenIntroFig}{Figure 1.1}
\newcommand{\ParallelogramFigures}{Figures 1.1, 4.1.A, and 4.2.A}

\newcommand{\CombinatorialSerre}{Proposition 2.1}

\newcommand{\EsixAndEsevenCharacter}{Proposition 3.1}
\newcommand{\PosetFigureEseven}{Figure 3.1}
\newcommand{\ArrayFigureEseven}{Figure 3.2}
\newcommand{\PosetFiguresESix}{Figure 3.3}
\newcommand{\MinusculeDepictions}{Figures 3.1 and 3.3}
\newcommand{\ArrayFiguresESix}{Figure 3.4}
\newcommand{\EsixMinusculeFigures}{Figures 3.3 and 3.4}
\newcommand{\StackedFigureESix}{Figure 3.5}
\newcommand{\StackedFigESix}{Fig.\ 3.5}
\newcommand{\StackedArrayESix}{Figure 3.6}
\newcommand{\StackedFigures}{Figures 3.5 and 3.6}

\newcommand{\EsevenParallelogramFigureAfivePartOne}{Figure 4.1.A}
\newcommand{\EsevenParallelogramFigureAfivePartTwo}{Figure 4.1.B}
\newcommand{\EsevenParallelogramFiguresAfive}{Figures 4.1.A {\footnotesize $\!\!$\&$\!\!$} B}
\newcommand{\JfiveLemma}{Lemma 4.1}
\newcommand{\EsevenParallelogramFigureAsixPartOne}{Figure 4.2.A}
\newcommand{\EsevenParallelogramFigureAsixPartTwo}{Figure 4.2.B}
\newcommand{\EsevenParallelogramFiguresAsix}{Figures 4.2.A {\footnotesize $\!\!$\&$\!\!$} B}
\newcommand{\JsixLemma}{Lemma 4.2}
\newcommand{\JLemmas}{Lemmas 4.1 and 4.2}
\newcommand{\AgreeLemma}{Lemma 4.3}
\newcommand{\ConstructionTheorem}{Theorem 4.4}
\newcommand{\ConstructionCorollaries}{Corollary 4.5}
\newcommand{\CoefficientFiguresA}{Figures 4.1.A and 4.2.A}
\newcommand{\CoefficientFigures}{Figures 4.1.B and 4.2.B}
\newcommand{\CoefficientFiguresPartOne}{Figures 4.1B and}
\newcommand{\CoefficientFiguresPartTwo}{4.2.B}
\newcommand{\AllCoefficientFigures}{Figures 4.1 and 4.2}

\newcommand{\CombinatorialResults}{Proposition 5.1}

\usepackage{amsfonts}
\usepackage{mathrsfs}
\usepackage{rotating}
\usepackage{mathdots}
\usepackage{hieroglf}
\usepackage{MnSymbol}

\newfont{\myscbolditalics}{ecoc0500 at 11pt}

\newfont{\mybolditalics}{ecbi0500 at 11pt}

\newcommand{\myqx}{\mbox{\mybolditalics x}}
\newcommand{\myqy}{\mbox{\mybolditalics y}}
\newcommand{\myqX}{\mbox{\mybolditalics X}}
\newcommand{\myqY}{\mbox{\mybolditalics Y}}
\newcommand{\myqh}{\mbox{\mybolditalics h}}
\newcommand{\myqP}{\mbox{\mybolditalics P}}
\newcommand{\myqQ}{\mbox{\mybolditalics Q}}

\newfont{\mysmallbolditalics}{ecbi0500 at 9pt}

\newcommand{\myqsmallx}{\mbox{\mysmallbolditalics x}}
\newcommand{\myqsmally}{\mbox{\mysmallbolditalics y}}
\newcommand{\myqsmallX}{\mbox{\mysmallbolditalics X}}
\newcommand{\myqsmallY}{\mbox{\mysmallbolditalics Y}}

\newcommand{\myqsmallP}{\mbox{\mysmallbolditalics P}}

\newfont{\eulercursive}{eurm10 at 11pt}

\newcommand{\mym}{\mbox{\eulercursive m}}

\newfont{\smalleulercursive}{eurm10 at 9pt}

\newfont{\smallereulercursive}{eurm10 at 7pt}

\newfont{\myslantcyrillic}{wncyi10 at 11pt}

\newcommand{\QED}{\raisebox{0.5mm}{\fbox{\rule{0mm}{1.5mm}\ }}}

\newcounter{myfn}[page]
\renewcommand{\thefootnote}{\fnsymbol{footnote}}

\newcounter{rone}
\newcounter{rtwo}
\newcounter{rthree}
\newcounter{rfour}
\newcounter{rfive}
\newcounter{rsix}
\newcounter{rseven}
\newcommand{\myA}{\mbox{\sffamily A}}
\newcommand{\mysmallA}{\mbox{\footnotesize \sffamily A}}
\newcommand{\mytinyA}{\mbox{\tiny \sffamily A}}
\newcommand{\myB}{\mbox{\sffamily B}}

\newcommand{\myC}{\mbox{\sffamily C}}

\newcommand{\myD}{\mbox{\sffamily D}}

\newcommand{\myE}{\mbox{\sffamily E}}
\newcommand{\mysmallE}{\mbox{\footnotesize \sffamily E}}
\newcommand{\mytinyE}{\mbox{\tiny \sffamily E}}
\newcommand{\myF}{\mbox{\sffamily F}}

\newcommand{\myG}{\mbox{\sffamily G}}

\newcommand{\mysmallP}{\mbox{\footnotesize \sffamily P}}

\newcommand{\mysmallQ}{\mbox{\footnotesize \sffamily Q}}

\newcommand{\myX}{\mbox{\sffamily X}}

\newcommand{\myvarZ}{\mbox{\scriptsize \sffamily Z}}

\newcommand{\melt}{\mathbf{m}} 
 \newcommand{\pelt}{\mathbf{p}}
\newcommand{\qelt}{\mathbf{q}} \newcommand{\relt}{\mathbf{r}}
\newcommand{\selt}{\mathbf{s}} \newcommand{\telt}{\mathbf{t}}
\newcommand{\uelt}{\mathbf{u}} \newcommand{\velt}{\mathbf{v}}

\newcommand{\wt}{\mbox{\sffamily wt}}
\newcommand{\smallwt}{\mbox{\scriptsize \sffamily wt}}

\newcommand{\mychar}{\mbox{\sffamily char}}

\newcommand{\WGF}{\mbox{\sffamily WGF}}
\newcommand{\RGF}{\mbox{\sffamily RGF}}

\newcommand{\comp}{\mbox{\sffamily comp}}

\newcommand{\myarrow}[1]{\stackrel{#1}{\rightarrow}}

\newcommand{\NEEdgeLabelForLatticeI}[1]{
\setlength{\unitlength}{1.5cm}
\begin{picture}(0,0)
\put(-0.25,0){
\begin{picture}(0,0)
\put(0.4,0.4){\footnotesize #1} 
\end{picture}
}
\end{picture}
}

\newcommand{\NWEdgeLabelForLatticeI}[1]{
\setlength{\unitlength}{1.5cm}
\begin{picture}(0,0)
\put(-0.25,0){
\begin{picture}(0,0)
\put(-0.525,0.4){\footnotesize #1} 
\end{picture}
}
\end{picture}
}

\newcommand{\VerticalEdgeLabelForLatticeI}[1]{
\setlength{\unitlength}{1.5cm}
\begin{picture}(0,0)
\put(-0.25,0){
\begin{picture}(0,0)
\put(-0.05,0.4){\footnotesize #1} 
\end{picture}
}
\end{picture}
}

\begin{document}
\pagenumbering{arabic}
\thispagestyle{empty}%
\vspace*{-0.7in}
\noindent
{\scriptsize To appear in {\em Applicable Algebra in Engineering, Communication, and Computing}}
\hfill {\scriptsize March 12, 2023}

\begin{center}
{\large \bf Explicit constructions of some infinite families of finite-dimensional irreducible\\ representations of the type} $\myE_{6}$ {\bf and} $\myE_{7}$ {\bf simple Lie algebras} 

\vspace*{0.05in}
\renewcommand{\thefootnote}{1}
Robert G.\ Donnelly\footnote{Department of Mathematics and Statistics, Murray State
University, Murray, KY 42071\\ 
\hspace*{0.25in}Email: {\tt rob.donnelly@murraystate.edu}}, 
\renewcommand{\thefootnote}{2} 
\hspace*{-0.07in}Molly W.\ Dunkum\footnote{Department of Mathematics, Western Kentucky University, Bowling Green, KY 42101\\ 
\hspace*{0.25in}Email: {\tt molly.dunkum@wku.edu}},
\renewcommand{\thefootnote}{3} 
\hspace*{-0.07in}and Austin White\footnote{Department of Mathematics, Western Kentucky University, Bowling Green, KY 42101\\ 
\hspace*{0.25in}Email: {\tt akwb@mit.edu}}

\end{center} 

\begin{abstract}
We construct every finite-dimensional irreducible representation of the simple Lie algebra of type $\myE_{7}$ whose highest weight is a nonnegative integer multiple of the dominant minuscule weight associated with the type $\myE_{7}$ root system. 
As a consequence, we obtain constructions of each finite-dimensional irreducible representation of the simple Lie algebra of type $\myE_{6}$ whose highest weight is a nonnegative integer linear combination of the two dominant minuscule $\myE_{6}$-weights. 
Our constructions are explicit in the sense that, if the representing space is $d$-dimensional, then a weight basis is provided such that all entries of the $d \times d$ representing matrices of the Chevalley generators are obtained via explicit, non-recursive formulas. 
To effect this work, we introduce what we call $\myE_{6}$- and $\myE_{7}$-polyminuscule lattices that analogize certain lattices associated with the famous special linear Lie algebra representation constructions obtained by Gelfand and Tsetlin. 

\begin{center}
{\small \bf Mathematics Subject Classification:}\ {\small 17B10 (22E70, 05E10)}\\
{\small \bf Keywords:}\ simple Lie algebra representation, root system, minuscule weight, Weyl group, Weyl symmetric function, ranked poset, diamond-colored distributive lattice, splitting poset, weight basis supporting graph$/$representation diagram

\end{center} 
\end{abstract}

\noindent {\bf \S \IntroSection\ Introduction}. 
A fundamental problem in the study of Lie groups and Lie algebras is the classification and construction of their linear representations. 
For any finite-dimensional semisimple Lie algebra over the complex number field, its finite-dimensional representations are completely reducible, and the finite-dimensional irreducible representations are nicely classified by dominant integral weights. 
The problem of constructing these irreducible representations can be addressed in many ways. 
Verma modules are a standard approach; 
Borel--Weil theory is another. 
However, such standard historical approaches are typically not computationally explicit. 

An exemplar of the kind of explicitness we seek is the famous constructions by Gelfand and Tsetlin of the irreducible representations of the special linear Lie algebras \cite{GT}. 
Their result was stated in 1950, but thorough proofs were not available until a decade or so later. 
(For details on some of this history, with ample references, see \cite{MoSurvey} and \cite{HL}.) 
Even so, the Gelfand--Tsetlin (`GT') constructions have several distinguishing and desirable features. 
In order to frame our overall perspective on explicitly constructing representations, we describe these features in terms of a generic rank $n$ simple Lie algebra $\mathfrak{g}(\myX_{n})$ (where $\myX \in \{\myA,\myB,\myC,\myD,\myE,\myF,\myG\}$ refers to the classification of simple Lie algebras by Dynkin diagrams as in \cite{Hum}) with Chevalley generators $\{\myqx_{i},\myqy_{i},\myqh_{i}\}_{i \in \{1,\ldots,n\}}$ satisfying the so-called Serre relations. 

\setcounter{footnote}{1}

Here are the explicitness features we desire. 
First, for any given dominant integral weight, a set $R$ of combinatorial objects should be prescribed as an indexing set for a `weight basis' $\{v_{\relt}\}_{\relt \in R}$ of the corresponding representing space.\footnote{That  $\{v_{\relt}\}_{\relt \in R}$ is a weight basis means that each $v_{\relt}$ is an eigenvector for each $\myqh_{i}$; the associated eigenvalues are integers.} 
Moreover, as an initial step in prescribing generator actions, the elements of $R$ should be naturally related by some directed edges that are `colored' by the set $\{1,\ldots,n\}$, with directed edges visualized as pointing upward.\footnote{The algebraic and combinatorial requisites for our explicit constructions ultimately necessitate that $R$ be a ranked poset whose covering relations are precisely our prescribed set of colored and directed edges.}
Second, each basis vector $v_{\relt}$ should be an eigenvector for each $\myqh_{i}$ in such a way that the associated integer eigenvalue can be computed by a (hopefully simple and combinatorial) rule applied to $\relt$. 
That is, for any weight basis vector $v_{\relt}$ and any color $i$, we should have $\myqh_{i}.v_{\relt} = \mym_{i}(\relt)v_{\relt}$ for some explicitly defined integer-valued function $\mym_{i}$ on $R$. 
Third, each $\myqx_{i}$ should act as a `raising' operator with respect to our directed edges so that $\myqx_{i}.v_{\relt} = \sum \myqX_{\selt,\relt}^{(i)}v_{\selt}$, where this sum is over all $\selt \in R$ such that $\relt \myarrow{i} \selt$ and the scalar coefficient\footnote{The constructions of this paper require square roots of positive rational numbers and thus are, in general, valid only over the smallest extension of $\mathbb{Q}$ that contains square roots of all prime numbers; despite this drawback, we obtain a representing matrix for each $\myqsmally_{i}$ that is, beneficently, the transpose of the representing matrix for its companion $\myqsmallx_{i}$.} $\myqX_{\selt,\relt}^{(i)}$ is specified by an explicit formula in terms of $\relt$ and $\selt$; similarly, each $\myqy_{i}$ should act as a `lowering' operator via $\myqy_{i}.v_{\selt} = \sum \myqY_{\relt,\selt}^{(i)}v_{\relt}$, where this sum is over all $\relt \in R$ such that $\relt \myarrow{i} \selt$ and $\myqY_{\relt,\selt}^{(i)}$ is given by an explicit formula in terms of $\relt$ and $\selt$. 

These salutary features of the GT constructions of the irreducible $\mathfrak{g}(\myA_{n})$-representations, whose characters align naturally with Schur functions, were highlighted and utilized in \cite{PrGZ} and \cite{HL}, and in \cite{DD1} they were extended to the much larger family of $\mathfrak{g}(\myA_{n})$-representations associated with all skew Schur functions. 
For a concrete example, see \EsixEsevenIntroFig.  
Moreover, these features are characteristic of the remarkable explicit constructions of weight bases for all irreducible representations of the other classical simple Lie algebras $\big(\mathfrak{g}(\myC_{n})$, $\mathfrak{g}(\myD_{n})$, $\mathfrak{g}(\myB_{n})\big)$ obtained by Molev in the series of papers \cite{Mo1}, \cite{Mo2}, \cite{Mo3}. 

For the exceptional simple Lie algebras -- i.e.\ of type $\myE_{6}$, $\myE_{7}$, $\myE_{8}$, $\myF_{4}$, and $\myG_{2}$ --  there remains much work to be done. 
The only explicit constructions we know for an infinite family of irreducible representations of an exceptional simple Lie algebra are the constructions in \cite{DLP1} of the `one-rowed' representations of $\mathfrak{g}(\myG_{2})$, i.e.\ those irreducible representations whose highest weights are nonnegative integer multiples of the dominant weight for the $7$-dimensional fundamental representation. 
However, those constructions were obtained by serendipity: It happens that, in the case $n=3$, the one-rowed representations of $\mathfrak{g}(\myB_{n})$ remain irreducible when viewed as $\mathfrak{g}(\myG_{2})$-modules via the natural inclusion $\mathfrak{g}(\myG_{2}) \hookrightarrow \mathfrak{g}(\myB_{3})$. 

Our purpose here is to present explicit constructions of some infinite families of irreducible representations of $\mathfrak{g}(\myE_{6})$ and $\mathfrak{g}(\myE_{7})$. 
In particular, we explicitly construct all irreducible representa-

\newcommand{\VertexParallelogramEsixEsevenIntro}[6]{\parbox[c]{1.2cm}{\setlength{\unitlength}{0.3cm}
\begin{picture}(3,5)
\put(0,2.2){\circle*{0.6}} 
\put(#5,#6){
\begin{picture}(0,0)
\put(0,-2){\begin{picture}(0,0)
\put(0.25,2.5){\tiny 0} \put(1.25,2){\tiny #1} \put(2.25,1.5){\tiny #3} \put(3.25,1){\tiny 3} 
\put(0.25,1.5){\tiny 2} \put(1.25,1){\tiny #2} \put(2.25,0.5){\tiny #4} \put(3.25,0){\tiny 3} 
\end{picture}}
\end{picture}} \end{picture}}
}
\newcommand{\VertexBlankEsixEsevenIntro}{\parbox[c]{1.2cm}{\setlength{\unitlength}{0.3cm}
\begin{picture}(3,5)
\put(0,2.2){\circle*{0.6}} 
\end{picture}}
}

\begin{figure}[ht]
\begin{center}
{\small {\bf \EsixEsevenIntroFig}\ \  The skew-tabular lattice $L :=$`$L_{\mytinyA_{2}}^{\mbox{\tiny skew}}({\setlength{\unitlength}{0.125cm}\begin{picture}(3,0)\put(0,0){\line(0,1){1}} \put(1,0){\line(0,1){1}} \put(2,0){\line(0,1){2}} \put(3,0){\line(0,1){2}} \put(0,0){\line(1,0){3}} \put(0,1){\line(1,0){3}} \put(2,2){\line(1,0){1}}\end{picture}})$' from \cite{DD1} for skew shape $\mysmallP/\mysmallQ = (3,3)/(2,0) = {\setlength{\unitlength}{0.125cm}\begin{picture}(3,0)\put(0,0){\line(0,1){1}} \put(1,0){\line(0,1){1}} \put(2,0){\line(0,1){2}} \put(3,0){\line(0,1){2}} \put(0,0){\line(1,0){3}} \put(0,1){\line(1,0){3}} \put(2,2){\line(1,0){1}}\end{picture}}$\\ 
built using the `Gelfand--Tsetlin $3$-parallelograms' {\parbox[t]{2cm}{\setlength{\unitlength}{0.3cm}
\begin{picture}(5,2.25)
\put(0,0.75){
\begin{picture}(0,0)
\put(0,-2){\begin{picture}(0,0)
\put(0.25,2.5){\tiny 0} \put(1.25,2){\tiny $g_{1,0}$} \put(3.25,1.5){\tiny $g_{2,1}$} \put(5.25,0.75){\tiny 3} 
\put(0.25,1.5){\tiny 2} \put(1.25,1){\tiny $g_{1,1}$} \put(3.25,0.5){\tiny $g_{2,2}$} \put(5.25,-0.25){\tiny 3} 
\end{picture}}
\end{picture}} 
\end{picture}}} framed by $\mysmallP/\mysmallQ$.} 
\end{center}
\hspace*{0.5in}
\setlength{\unitlength}{1.5cm}
\begin{picture}(4,6.5)
\thicklines
\put(1,0){\textcolor{CornflowerBlue}{\vector(-1,1){1}}}
\put(1,0){\textcolor{WildStrawberry}{\vector(1,1){1}}}
\put(2,1){\textcolor{CornflowerBlue}{\vector(1,1){1}}}
\put(0,1){\textcolor{WildStrawberry}{\vector(1,1){1}}}
\put(1,2){\textcolor{CornflowerBlue}{\vector(1,1){1}}}
\put(2,1){\textcolor{CornflowerBlue}{\vector(-1,1){1}}}
\put(3,2){\textcolor{CornflowerBlue}{\vector(-1,1){1}}}
\put(2,1){\textcolor{WildStrawberry}{\vector(0,1){1}}}
\put(1,2){\textcolor{WildStrawberry}{\vector(0,1){1}}}
\put(3,2){\textcolor{WildStrawberry}{\vector(0,1){1}}}
\put(2,3){\textcolor{WildStrawberry}{\vector(0,1){1}}}
\put(2,2){\textcolor{CornflowerBlue}{\vector(-1,1){1}}}
\put(3,3){\textcolor{CornflowerBlue}{\vector(-1,1){1}}}
\put(4,4){\textcolor{CornflowerBlue}{\vector(-1,1){1}}}
\put(2,2){\textcolor{CornflowerBlue}{\vector(1,1){1}}}
\put(1,3){\textcolor{CornflowerBlue}{\vector(1,1){1}}}
\put(0,4){\textcolor{CornflowerBlue}{\vector(1,1){1}}}
\put(3,3){\textcolor{CornflowerBlue}{\vector(1,1){1}}}
\put(2,4){\textcolor{CornflowerBlue}{\vector(1,1){1}}}
\put(1,5){\textcolor{CornflowerBlue}{\vector(1,1){1}}}
\put(1,3){\textcolor{WildStrawberry}{\vector(-1,1){1}}}
\put(2,4){\textcolor{WildStrawberry}{\vector(-1,1){1}}}
\put(3,5){\textcolor{WildStrawberry}{\vector(-1,1){1}}}
\thinlines
\put(2,6){\VertexParallelogramEsixEsevenIntro{2}{3}{3}{3}{0.0}{2.5}}
\put(1,5){\VertexParallelogramEsixEsevenIntro{1}{3}{3}{3}{-4.75}{4}}
\put(3,5){\VertexParallelogramEsixEsevenIntro{2}{3}{2}{3}{0.0}{2.5}}
\put(0,4){\VertexParallelogramEsixEsevenIntro{0}{3}{3}{3}{-4.75}{4}}
\put(2,4){\VertexParallelogramEsixEsevenIntro{1}{3}{2}{3}{-6.75}{3}}
\put(1.55,4){\vector(1,0){0.35}}
\put(4,4){\VertexParallelogramEsixEsevenIntro{2}{2}{2}{3}{0.0}{2.5}}
\put(1,3){\VertexParallelogramEsixEsevenIntro{0}{3}{2}{3}{-6.75}{3}}
\put(0.55,3){\vector(1,0){0.35}}
\put(2,3){\VertexParallelogramEsixEsevenIntro{1}{3}{2}{3}{2}{6.5}}
\put(2.5,3.75){\vector(-2,-3){0.45}}
\put(3,3){\VertexParallelogramEsixEsevenIntro{1}{2}{2}{3}{0.25}{0.75}}
\put(1,2){\VertexParallelogramEsixEsevenIntro{0}{3}{1}{3}{-6.75}{3}}
\put(0.55,2){\vector(1,0){0.35}}
\put(2,2){\VertexParallelogramEsixEsevenIntro{0}{2}{2}{3}{-6.5}{-2}}
\put(1.55,1.1){\vector(1,2){0.375}}
\put(3,2){\VertexParallelogramEsixEsevenIntro{1}{2}{1}{3}{0.25}{0.75}}
\put(0,1){\VertexParallelogramEsixEsevenIntro{0}{3}{0}{3}{-4.75}{4}}
\put(2,1){\VertexParallelogramEsixEsevenIntro{0}{2}{1}{3}{0.25}{0.75}}
\put(1,0){\VertexParallelogramEsixEsevenIntro{0}{2}{0}{3}{-6.75}{3}}
\put(0.55,0){\vector(1,0){0.35}}
\put(1,5){\NEEdgeLabelForLatticeI{\textcolor{CornflowerBlue}{\large \bf \em 1}}}
\put(2.95,5){\NWEdgeLabelForLatticeI{\textcolor{WildStrawberry}{\large \bf \em 2}}}
\put(0,4){\NEEdgeLabelForLatticeI{\textcolor{CornflowerBlue}{\large \bf \em 1}}}
\put(1.95,4){\NWEdgeLabelForLatticeI{\textcolor{WildStrawberry}{\large \bf \em 2}}}
\put(2,4){\NEEdgeLabelForLatticeI{\textcolor{CornflowerBlue}{\large \bf \em 1}}}
\put(4,4){\NWEdgeLabelForLatticeI{\textcolor{CornflowerBlue}{\large \bf \em 1}}}
\put(1,3){\NEEdgeLabelForLatticeI{\textcolor{CornflowerBlue}{\large \bf \em 1}}}
\put(0.95,3){\NWEdgeLabelForLatticeI{\textcolor{WildStrawberry}{\large \bf \em 2}}}
\put(1.95,3){\VerticalEdgeLabelForLatticeI{\textcolor{WildStrawberry}{\large \bf \em 2}}}
\put(3,3){\NWEdgeLabelForLatticeI{\textcolor{CornflowerBlue}{\large \bf \em 1}}}
\put(3,3){\NEEdgeLabelForLatticeI{\textcolor{CornflowerBlue}{\large \bf \em 1}}}
\put(0.95,2){\VerticalEdgeLabelForLatticeI{\textcolor{WildStrawberry}{\large \bf \em 2}}}
\put(1.25,2.25){\NEEdgeLabelForLatticeI{\textcolor{CornflowerBlue}{\large \bf \em 1}}}
\put(2.2,1.8){\NWEdgeLabelForLatticeI{\textcolor{CornflowerBlue}{\large \bf \em 1}}}
\put(1.8,1.8){\NEEdgeLabelForLatticeI{\textcolor{CornflowerBlue}{\large \bf \em 1}}}
\put(2.95,2){\VerticalEdgeLabelForLatticeI{\textcolor{WildStrawberry}{\large \bf \em 2}}}
\put(2.75,2.25){\NWEdgeLabelForLatticeI{\textcolor{CornflowerBlue}{\large \bf \em 1}}}
\put(0,1){\NEEdgeLabelForLatticeI{\textcolor{WildStrawberry}{\large \bf \em 2}}}
\put(1.95,1){\VerticalEdgeLabelForLatticeI{\textcolor{WildStrawberry}{\large \bf \em 2}}}
\put(2,1){\NWEdgeLabelForLatticeI{\textcolor{CornflowerBlue}{\large \bf \em 1}}}
\put(2,1){\NEEdgeLabelForLatticeI{\textcolor{CornflowerBlue}{\large \bf \em 1}}}
\put(1,0){\NWEdgeLabelForLatticeI{\textcolor{CornflowerBlue}{\large \bf \em 1}}}
\put(1,0){\NEEdgeLabelForLatticeI{\textcolor{WildStrawberry}{\large \bf \em 2}}}
\put(5.5,5.25){\parbox[t]{2.5in}{\footnotesize The integers in these parallelogram-shaped arrays are required to weakly increase along diagonals from NW to SE and from NE to SW. 
For notational consistency with the constructions of \cite{DD1}, we let $(g_{0,0},g_{0,-1}) := (2,0) = \mysmallQ$ and $(g_{3,3},g_{3,2}) := (3,3) = \mysmallP$; these columns are fixed for all arrays. 
One can see that the implied partial ordering in the figure to the left is componentwise comparison: $\relt \leq \selt$ for arrays $\relt$ and $\selt$ in $L$ if and only if $g_{pq}(\relt) \leq g_{pq}(\selt)$ at all positions $(p,q)$ in the array. 
Also observe that there is a directed edge of color $i$ from an array $\relt$ to an array $\selt$, i.e.\ $\relt \myarrow{i} \selt$, if and only if there is a position $(i,j)$ in the array with $g_{ij}(\relt)+1=g_{ij}(\selt)$ while $g_{pq}(\relt) = g_{pq}(\selt)$ at all other positions $(p,q)$.}}
\end{picture}

\noindent 
\hspace*{0.5in}
\setlength{\unitlength}{1.5cm}
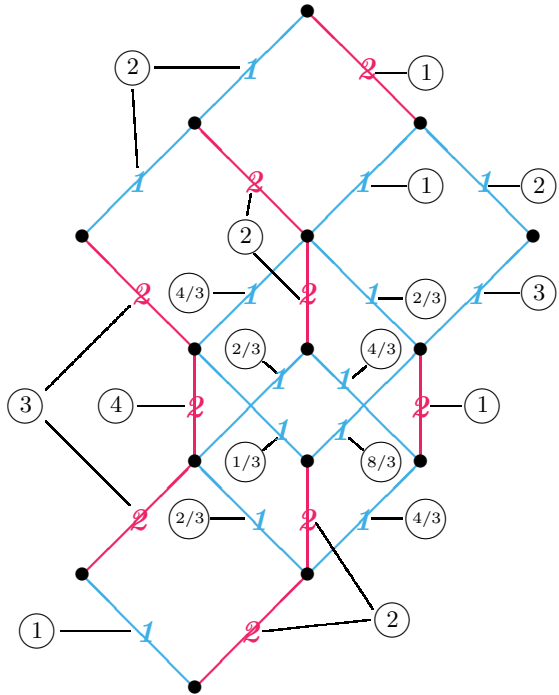
\begin{picture}(4,6.5)
\thicklines
\put(1,0){\textcolor{CornflowerBlue}{\line(-1,1){1}}}
\put(1,0){\textcolor{WildStrawberry}{\line(1,1){1}}}
\put(2,1){\textcolor{CornflowerBlue}{\line(1,1){1}}}
\put(0,1){\textcolor{WildStrawberry}{\line(1,1){1}}}
\put(1,2){\textcolor{CornflowerBlue}{\line(1,1){1}}}
\put(2,1){\textcolor{CornflowerBlue}{\line(-1,1){1}}}
\put(3,2){\textcolor{CornflowerBlue}{\line(-1,1){1}}}
\put(2,1){\textcolor{WildStrawberry}{\line(0,1){1}}}
\put(1,2){\textcolor{WildStrawberry}{\line(0,1){1}}}
\put(3,2){\textcolor{WildStrawberry}{\line(0,1){1}}}
\put(2,3){\textcolor{WildStrawberry}{\line(0,1){1}}}
\put(2,2){\textcolor{CornflowerBlue}{\line(-1,1){1}}}
\put(3,3){\textcolor{CornflowerBlue}{\line(-1,1){1}}}
\put(4,4){\textcolor{CornflowerBlue}{\line(-1,1){1}}}
\put(2,2){\textcolor{CornflowerBlue}{\line(1,1){2}}}
\put(1,3){\textcolor{CornflowerBlue}{\line(1,1){2}}}
\put(0,4){\textcolor{CornflowerBlue}{\line(1,1){2}}}
\put(1,3){\textcolor{WildStrawberry}{\line(-1,1){1}}}
\put(2,4){\textcolor{WildStrawberry}{\line(-1,1){1}}}
\put(3,5){\textcolor{WildStrawberry}{\line(-1,1){1}}}
\thinlines
\put(2,6){\VertexBlankEsixEsevenIntro}
\put(1,5){\VertexBlankEsixEsevenIntro}
\put(3,5){\VertexBlankEsixEsevenIntro}
\put(0,4){\VertexBlankEsixEsevenIntro}
\put(2,4){\VertexBlankEsixEsevenIntro}
\put(4,4){\VertexBlankEsixEsevenIntro}
\put(1,3){\VertexBlankEsixEsevenIntro}
\put(2,3){\VertexBlankEsixEsevenIntro}
\put(3,3){\VertexBlankEsixEsevenIntro}
\put(1,2){\VertexBlankEsixEsevenIntro}
\put(2,2){\VertexBlankEsixEsevenIntro}
\put(3,2){\VertexBlankEsixEsevenIntro}
\put(0,1){\VertexBlankEsixEsevenIntro}
\put(2,1){\VertexBlankEsixEsevenIntro}
\put(1,0){\VertexBlankEsixEsevenIntro}
%
\put(-0.5,2.5){\circle{0.3}}
\put(-0.55,2.44){\footnotesize 3}
\put(-0.32,2.64){\qbezier(0,0)(0.37,0.37)(0.74,0.74)}
\put(-0.32,2.36){\qbezier(0,0)(0.37,-0.37)(0.74,-0.74)}
\put(2.75,0.6){\circle{0.3}}
\put(2.7,0.54){\footnotesize 2}
\put(2.55,0.6){\qbezier(0,0)(-0.47,-0.05)(-0.94,-0.1)}
\put(2.575,0.7){\qbezier(0,0)(-0.25,0.38)(-0.5,0.76)}
\put(-0.4,0.5){\circle{0.3}}
\put(-0.45,0.44){\footnotesize 1}
\put(-0.18,0.5){\qbezier(0,0)(0.3,0)(0.6,0)}
\put(0.3,2.5){\circle{0.3}}
\put(0.25,2.44){\footnotesize 4}
\put(0.5,2.5){\qbezier(0,0)(0.2,0)(0.4,0)}
\put(0.95,1.5){\circle{0.34}}
\put(0.83,1.46){\tiny 2/3}
\put(1.15,1.5){\qbezier(0,0)(0.14,0)(0.28,0)}
\put(1.45,2){\circle{0.34}}
\put(1.33,1.96){\tiny 1/3}
\put(1.61,2.11){\qbezier(0,0)(0.06,0.05)(0.12,0.10)}
\put(1.45,3){\circle{0.34}}
\put(1.33,2.96){\tiny 2/3}
\put(1.61,2.89){\qbezier(0,0)(0.06,-0.05)(0.12,-0.10)}
\put(2.65,2){\circle{0.34}}
\put(2.53,1.96){\tiny 8/3}
\put(2.49,2.11){\qbezier(0,0)(-0.06,0.05)(-0.12,0.10)}
\put(2.65,3){\circle{0.34}}
\put(2.53,2.96){\tiny 4/3}
\put(2.52,2.85){\qbezier(0,0)(-0.06,-0.05)(-0.12,-0.10)}
\put(3.05,1.5){\circle{0.34}}
\put(2.93,1.46){\tiny 4/3}
\put(2.83,1.5){\qbezier(0,0)(-0.13,0)(-0.26,0)}
\put(3.55,2.5){\circle{0.3}}
\put(3.5,2.44){\footnotesize 1}
\put(3.36,2.5){\qbezier(0,0)(-0.13,0)(-0.26,0)}
\put(4.05,3.5){\circle{0.3}}
\put(4,3.44){\footnotesize 3}
\put(3.87,3.5){\qbezier(0,0)(-0.15,0)(-0.30,0)}
\put(0.95,3.5){\circle{0.34}}
\put(0.83,3.46){\tiny 4/3}
\put(1.17,3.5){\qbezier(0,0)(0.13,0)(0.26,0)}
\put(3.05,3.45){\circle{0.34}}
\put(2.93,3.41){\tiny 2/3}
\put(2.83,3.45){\qbezier(0,0)(-0.10,0)(-0.20,0)}
\put(1.45,4){\circle{0.3}}
\put(1.4,3.94){\footnotesize 2}
\put(1.465,4.195){\qbezier(0,0)(0.015,0.09)(0.03,0.18)}
\put(1.525,3.85){\qbezier(0,0)(0.2,-0.2)(0.4,-0.4)}
\put(4.05,4.45){\circle{0.3}}
\put(4,4.39){\footnotesize 2}
\put(3.87,4.45){\qbezier(0,0)(-0.13,0)(-0.26,0)}
\put(3.05,4.45){\circle{0.3}}
\put(3,4.39){\footnotesize 1}
\put(2.87,4.45){\qbezier(0,0)(-0.15,0)(-0.3,0)}
\put(3.05,5.45){\circle{0.3}}
\put(3,5.39){\footnotesize 1}
\put(2.87,5.45){\qbezier(0,0)(-0.13,0)(-0.26,0)}
\put(0.45,5.5){\circle{0.3}}
\put(0.4,5.44){\footnotesize 2}
\put(0.65,5.5){\qbezier(0,0)(0.37,0)(0.74,0)}
\put(0.45,5.3){\qbezier(0,0)(0.02,-0.34)(0.04,-0.68)}
\put(1,5){\NEEdgeLabelForLatticeI{\textcolor{CornflowerBlue}{\large \bf \em 1}}}
\put(2.95,5){\NWEdgeLabelForLatticeI{\textcolor{WildStrawberry}{\large \bf \em 2}}}
\put(0,4){\NEEdgeLabelForLatticeI{\textcolor{CornflowerBlue}{\large \bf \em 1}}}
\put(1.95,4){\NWEdgeLabelForLatticeI{\textcolor{WildStrawberry}{\large \bf \em 2}}}
\put(2,4){\NEEdgeLabelForLatticeI{\textcolor{CornflowerBlue}{\large \bf \em 1}}}
\put(4,4){\NWEdgeLabelForLatticeI{\textcolor{CornflowerBlue}{\large \bf \em 1}}}
\put(1,3){\NEEdgeLabelForLatticeI{\textcolor{CornflowerBlue}{\large \bf \em 1}}}
\put(0.95,3){\NWEdgeLabelForLatticeI{\textcolor{WildStrawberry}{\large \bf \em 2}}}
\put(1.95,3){\VerticalEdgeLabelForLatticeI{\textcolor{WildStrawberry}{\large \bf \em 2}}}
\put(3,3){\NWEdgeLabelForLatticeI{\textcolor{CornflowerBlue}{\large \bf \em 1}}}
\put(3,3){\NEEdgeLabelForLatticeI{\textcolor{CornflowerBlue}{\large \bf \em 1}}}
\put(0.95,2){\VerticalEdgeLabelForLatticeI{\textcolor{WildStrawberry}{\large \bf \em 2}}}
\put(1.25,2.25){\NEEdgeLabelForLatticeI{\textcolor{CornflowerBlue}{\large \bf \em 1}}}
\put(2.2,1.8){\NWEdgeLabelForLatticeI{\textcolor{CornflowerBlue}{\large \bf \em 1}}}
\put(1.8,1.8){\NEEdgeLabelForLatticeI{\textcolor{CornflowerBlue}{\large \bf \em 1}}}
\put(2.95,2){\VerticalEdgeLabelForLatticeI{\textcolor{WildStrawberry}{\large \bf \em 2}}}
\put(2.75,2.25){\NWEdgeLabelForLatticeI{\textcolor{CornflowerBlue}{\large \bf \em 1}}}
\put(0,1){\NEEdgeLabelForLatticeI{\textcolor{WildStrawberry}{\large \bf \em 2}}}
\put(1.95,1){\VerticalEdgeLabelForLatticeI{\textcolor{WildStrawberry}{\large \bf \em 2}}}
\put(2,1){\NWEdgeLabelForLatticeI{\textcolor{CornflowerBlue}{\large \bf \em 1}}}
\put(2,1){\NEEdgeLabelForLatticeI{\textcolor{CornflowerBlue}{\large \bf \em 1}}}
\put(1,0){\NWEdgeLabelForLatticeI{\textcolor{CornflowerBlue}{\large \bf \em 1}}}
\put(1,0){\NEEdgeLabelForLatticeI{\textcolor{WildStrawberry}{\large \bf \em 2}}}
\put(5.5,5.5){\parbox[t]{2.5in}{\footnotesize In the picture to the left, the directed edges between vertices are represented as segments rather than as arrows, with the understanding that all edges are directed `up'. 
Assigned to each edge $\relt \myarrow{i} \selt$ is a circled number we denote $\myqsmallP^{(i)}_{\relt,\selt}$.  
These numbers were obtained using formulas analogous to those of \CoefficientFiguresA. 
Set $\myqsmallX^{(i)}_{\selt,\relt} := \sqrt{\myqsmallP^{(i)}_{\relt,\selt}\, } =: \myqsmallY^{(i)}_{\relt,\selt}$. 
These `edge-coefficients' satisfy the diamond and crossing relations of \S \Setup. 
In fact, \CombinatorialSerre\ applies here, so that we get a well-defined action of the simple Lie algebra $\mathfrak{g}(\mysmallA_{2}) \cong \mathfrak{sl}(3,\mathbb{C})$ on the 15-dimensional vector space spanned by the weight basis vectors $\{v_{\relt}\}_{\relt \in L}$ if we define generator actions by the rules $\myqsmallx_{i}.v_{\relt} := \sum \myqsmallX^{(i)}_{\selt,\relt}v_{\selt}$ and $\myqsmallx_{i}.v_{\selt} := \sum \myqsmallY^{(i)}_{\relt,\selt}v_{\relt}$,.}}
\end{picture}

\end{figure}
\clearpage

\noindent 
tions of $\mathfrak{g}(\myE_{7})$ whose highest weights are nonnegative integer multiples of $\myE_{7}$'s dominant minuscule weight, cf.\ \ConstructionTheorem. 
As a consequence, we obtain explicit constructions of all irreducible representations of $\mathfrak{g}(\myE_{6})$ whose highest weights are nonnegative integer linear combinations of $\myE_{6}$'s two dominant minuscule weights, cf.\ \ConstructionCorollaries. 
Our constructions are from scratch but rely crucially on some reducible $\mathfrak{g}(\myA_{n})$-representation constructions, with $n \in \{5,6\}$, obtained in \cite{DD1}. 
The $\mathfrak{g}(\myA_{n})$-representation constructions of that paper took place within the combinatorial setting of what we call skew-tabular lattices.   
The key aspects of those constructions are illustrated in \EsixEsevenIntroFig. 
For our work here, we introduce what we call $\myE_{6}$- and $\myE_{7}$-polyminuscule lattices to serve as type $\myE$ analogs of the skew-tabular lattices for type $\myA$. 
Our choice of this combinatorial setting, and indeed our overall perspective, was strongly influenced by Proctor's work in \cite{PrEur} and \cite{PrGZ}.

\vspace*{0.2in}
\noindent {\bf \S \Setup\ General set-up}. 
In this section we set our language $\!\!$/$\!\!$ notation and declare the general algebraic and combinatorial environment we will be working in. 
The first-time reader is encouraged to lightly browse this section and use it as a reference when reading later sections. 
See, for example, \cite{Jac}, \cite{Hum}, or \cite{FH} for further details on standard notions from Lie algebra representation theory; for order-theoretic combinatorics, see, for example, \cite{Aigner} or \cite{StanText}. 
Our particular perspective on how these two areas can work together is more fully developed in \cite{DonSupp}, \cite{DonDistributive}, and \cite{DonPosetModels}, but \cite{DD1} is briefer and most immediately relevant to our work here. 

Unless otherwise stated, all vector spaces in this paper are complex and finite-dimensional and all partially ordered sets (`posets') are finite. 
Only idiosyncratic notions$/$language will be {\em italicized}. 

Fix a semisimple Lie algebra $\mathfrak{g}$ and Cartan subalgebra $\mathfrak{h}$. 
Say $\dim \mathfrak{h} = n$, so $\mathfrak{g}$ has rank $n$. 
Note that $\mathfrak{g}$ is uniquely identified (up to isomorphism) by its Dynkin diagram, which is a simple graph on $n$ nodes whose edges carry additional information from which we can easily re-construct $\mathfrak{g}$ by applying the Serre relations to our set $\{\myqx_{i},\myqy_{i},\myqh_{i}\}_{i \in I}$ of Chevalley generators. 
(Here, $I$ is some indexing set of size $n$.) 
So, the set $\{\myqh_{i}\}_{i \in I}$ is a basis for $\mathfrak{h}$. 

Associated to $\mathfrak{g}$ is a root system $\Phi$ residing in an $n$-dimensional Euclidean space $\mathfrak{E}$ with inner product $\langle \cdot,\cdot \rangle$. 
For convenience, we sometimes use `$\Phi$' as an argument or super/subscript to contextualize a given object, e.g.\ `$\mathfrak{g}(\Phi)$'. 
There is a basis $\{\alpha_{i}\}_{i \in I} \subseteq \Phi$ for $\mathfrak{E}$ wherein any $\alpha \in \Phi$ can be expressed as $\alpha = \sum_{i \in I}k_{i}\alpha_{i}$ for integers $k_{i}$ that are either all nonnegative (in which case $\alpha \in \Phi^{+}$ is `positive') or all nonpositive (in which case $\alpha \in \Phi^{-}$ is `negative'). 
The $\alpha_{i}$'s are simple roots. 
For any $\alpha \in \Phi$, its coroot $\alpha^{\vee}$ is defined as the vector $\alpha/\langle \alpha,\alpha \rangle$. 
The matrix of integers $(\langle \alpha_{i},\alpha_{j}^{\vee} \rangle)_{i,j \in I}$ is the Cartan matrix. 
Let $\{\omega_{i}\}_{i \in I}$ be the basis dual to $\{\alpha_{i}^{\vee}\}_{i \in I}$, so $\langle \omega_{i},\alpha_{j}^{\vee} \rangle = \delta_{i,j}$. 
The $\omega_{i}$'s are fundamental weights. 
The lattice of weights $\Lambda$ is the $\mathbb{Z}$-span of $\{\omega_{i}\}_{i \in I}$. 
A weight $\lambda = \sum_{i \in I}\lambda_{i}\omega_{i} \in \Lambda$ is dominant if each $\lambda_{i}$ is nonnegative. 
A nonzero dominant weight $\lambda$ is minuscule if $\langle \lambda,\alpha^{\vee} \rangle \in \{0,\pm 1\}$ for all $\alpha \in \Phi^{+}$; we sometimes say that $\lambda$ is {\em $\Phi$-minuscule} to emphasize the provenance of $\lambda$ as a $\Phi$-related weight.

The Weyl group $\mathscr{W} = \mathscr{W}(\Phi)$ is generated by the reflections $\{S_{i}: \mathfrak{E} \longrightarrow \mathfrak{E}\}_{i \in I}$ wherein $S_{i}(v) = v-\langle v,\alpha_{i}^{\vee} \rangle\alpha_{i}$. 
Note that $S_{i}(\alpha_{j}) = \alpha_{j} - \langle \alpha_{j},\alpha_{i}^{\vee} \rangle\alpha_{i}$, so $\mathscr{W}$ preserves the root lattice $\mathbb{Z}\Phi$. 
Also, $S_{i}(\omega_{j}) = \omega_{j} - \langle \omega_{j},\alpha_{i}^{\vee} \rangle\alpha_{i}$ and $\alpha_{i} = \sum_{j \in I}\langle \alpha_{i},\alpha_{j}^{\vee} \rangle\omega_{j}$, so $\mathscr{W}$ preserves the lattice of weights $\Lambda$. 
One can easily check that $S_{i}^{2}$ is the identity and that $\det(S_{i}) = -1$. 
Moreover, $\det:\mathscr{W} \longrightarrow \{\pm 1\}$ is a (well-defined) homomorphism. 

An irreducible representation $V$ of $\mathfrak{g}$ is identified by the dominant integral weight of its unique (up to scalar multiples) maximal vector. 
That is, if a nonzero vector $v \in V$ has the property that $\myqx_{i}.v = 0$ for each $i$, then we get $\myqh_{i}.v = \lambda_{i} v$ for some nonnegative integers $\lambda_{i}$, in which case the dominant weight $\lambda := \sum \lambda_{i}\omega_{i}$ is the highest weight associated with $V$. 
Let $\mu = \sum \mu_{i}\omega_{i} \in \Lambda$ and let $W$ be any $\mathfrak{g}$-module. 
The $\mu$ weight space $W_{\mu}$ is the subspace of $W$ consisting of all vectors $v$ such that $\myqh_{i}.v = \mu_{i}v$ for all $i \in I$. 
Then $W = \bigoplus_{\nu \in \Lambda}W_{\nu}$, and any basis respecting this decomposition is a weight basis. 

Let $\{z_{i}\}_{i \in I}$ be a set of indeterminates, and for any $\mu \in \Lambda$ write $\myvarZ^{\mu} := \prod_{i \in I}z_{i}^{\mu_{i}}$, a Laurent monomial. 
The character $\mychar(\lambda) = \mychar_{\Phi}(\lambda)$ of our highest-weight-$\lambda$ irreducible representation $V$ is the formal sum $\sum_{\mu \in \Lambda}\dim(V_{\mu})\myvarZ^{\mu}$. 
Any given $\mathfrak{g}$-module $W$ decomposes uniquely as a direct sum $\bigoplus_{i=1}^{m}V_{i}$ of irreducible sub-modules, where each irreducible $V_{i}$ corresponds to some dominant integral weight $\lambda^{(i)}$. 
Then $\mychar(W) :=  \sum_{\mu \in \Lambda}\dim(W_{\mu})\myvarZ^{\mu} = \sum_{i=1}^{m}\mychar(\lambda^{(i)})$. 
We can naturally extend the action of $\mathscr{W}$ on weights $\mu$ to an action of $\mathscr{W}$ on monomials $\myvarZ^{\mu}$ by the rule $\sigma.\myvarZ^{\mu} := \myvarZ^{\sigma.\mu}$ for all $\sigma \in \mathscr{W}$. 
In turn, this extends naturally to an action of $\mathscr{W}$ on the set `$\mathbb{Z}[\Lambda]$' all Laurent polynomials $\sum_{\mu \in \Lambda}c_{\mu}\myvarZ^{\mu}$ having integer coefficients $c_{\mu}$ only finitely many of which are nonzero.  
A basic fact about any representation $W$ of $\mathfrak{g}$ is that for any $\sigma \in \mathscr{W}$, we have $\dim(W_{\mu}) = \dim(W_{\sigma.\mu})$.  
That is, $\mychar(W)$ is $\mathscr{W}$-invariant. 
Within $\mathbb{Z}[\Lambda]$, the $\mathscr{W}$-invariant elements are {\em Weyl symmetric functions}. 
Let $\varrho := \sum_{i \in I}\omega_{i}$, and for any dominant weight $\lambda$ set $\mathcal{A}(\myvarZ^{\lambda+\varrho}) := \sum_{\sigma \in \mathscr{W}}\det(\sigma)\myvarZ^{\sigma.(\lambda+\varrho)}$. 
Note that the quantity $\mathcal{A}(\myvarZ^{\lambda+\varrho})$ is an {\em alternant} in the sense that $S_{i}.\mathcal{A}(\myvarZ^{\lambda+\varrho}) = -\mathcal{A}(\myvarZ^{\lambda+\varrho})$ for all $i \in I$. 
The famous Weyl character formula can be viewed as the assertion that the unique Laurent polynomial $\chi$ satisfying the equation $\mathcal{A}(\myvarZ^{\varrho}) \chi = \mathcal{A}(\myvarZ^{\lambda+\varrho})$ is the $\mathscr{W}(\Phi)$-invariant character $\chi = \mychar_{\Phi}(\lambda)$, sometimes also denoted `$\chi_{_{\lambda}}^{\Phi}$' and called a {\em Weyl bialternant}. 
It is well-known that Weyl bialternants comprise a $\mathbb{Z}$-basis for the ring of Weyl symmetric functions. 

Our representation constructions make use of some rudimentary concepts related to posets. 
That said, no advanced poset theory is needed, and the notions we require are intuitive. 
A poset is a set $R$ together with a relation `$\leq$' that is reflexive, transitive, and anti-symmetric (i.e.\ $\relt \leq \selt$ and $\selt \leq \relt$ $\Longrightarrow$ $\relt=\selt$). 
Although the poset is specified by the pair $(R,\leq)$, we refer to the poset simply as $R$ when the partial order is understood. 
In $R$, say $\selt$ covers $\relt$ and write $\relt \rightarrow \selt$ if there is no $\telt$ in $R$ such that $\relt < \telt < \selt$. 
The {\em covering digraph} (aka `Hasse diagram') of $R$ is the directed graph whose vertex set $\mathcal{V} = \mathcal{V}(R)$ is the elements of $R$ and whose directed-edge set $\mathcal{E} = \mathcal{E}(R)$ is the set of covering relations in $R$. 
As a convenient abuse of notation, we often identify a poset with its covering digraph. 
In this way, concepts that apply to digraphs extend to posets. 
So, a poset is connected if its covering digraph is (weakly) connected, we may speak of adjacency of poset elements, etc. 
An `undirected path' in $R$ is a path in the covering digraph of $R$ where we regard edges to be undirected. 
Typically, we depict the directed edges of the covering digraph as pointing upward; if so, we omit the arrowhead on the directed edge. 
Given a set $J$, to be thought of as colors, a vertex-coloring of $R$ is a function $\mathscr{V}_{\mbox{\tiny color}}: \mathcal{V}(R) \longrightarrow J$, and an edge-coloring of $R$ is a function $\mathscr{E}_{\mbox{\tiny color}}: \mathcal{E}(R) \longrightarrow J$. 
In our work, a poset usually has a vertex-coloring or an edge-coloring, but not both; colors, for us, most often refer to the nodes of a Dynkin diagram or a set of simple roots. 
If $\psi: J \longrightarrow J'$ is a set mapping, then $R^{\psi}$ denotes the vertex-colored (respectively, edge-colored) poset colored by $\psi \circ \mathscr{V}_{\mbox{\tiny color}}$ (resp., $\psi \circ \mathscr{E}_{\mbox{\tiny color}}$). 

A poset $R$ is ranked if there is a surjective function $\rho: R \longrightarrow \{0,1,\ldots,\ell\}$ such that $\rho(\relt)+1=\rho(\selt)$ whenever $\relt \rightarrow \selt$, in which case $\rho$ is its rank function and $\ell$ is the length of $R$ with respect to $\rho$. 
The associated depth function is the mapping $\delta: R \longrightarrow \{0,1,\ldots,\ell\}$ given by $\delta(\relt) := \ell - \rho(\relt)$. 
Observe that a connected ranked poset has a unique rank function and a unique depth function. 
The rank generating function for $R$ is the $q$-polynomial $\RGF(R;q) := \sum_{\relt \in R}q^{\rho(\relt)}$. 
For example, for the ranked poset $L$ of \EsixEsevenIntroFig, $\RGF(L;q) = 1+2q+3q^{2}+3q^{3}+3q^{4}+2q^{5}+q^{6}$. 
Say our ranked poset $R$ is edge-colored by a set $J$. 
For any $K \subseteq J$ and any $\relt \in R$, $\comp_{K}(\relt)$ is the edge-colored directed graph of all elements of $R$ that can be reached from $\relt$ via undirected paths whose edges only have colors from $K$ together with the directed edges of all such possible paths. 
Notice that $\comp_{K}(\relt)$ is the connected edge-colored covering digraph of a ranked poset. 
When $K$ is the singleton set $\{k\}$ consisting of only the color $k$, we let $\rho_{k}$ denote the unique rank function of $\comp_{k}(\relt)$ and $\delta_{k}$ its unique depth function. 
For any $\selt \in \comp_{k}(\relt)$, we let $\mym_{k}(\selt) := \rho_{k}(\selt) - \delta_{k}(\selt)$. 

A poset $L$ is a lattice if each pair of poset elements has a unique least upper bound and a unique greatest lower bound. 
That is, for all $\relt, \selt \in L$ there are elements, denoted $\relt \vee \selt$ (the join of $\relt$ and $\selt$) and $\relt \wedge \selt$ (the meet of $\relt$ and $\selt$), such that $\relt \vee \selt \leq \telt$ whenever $\relt \leq \telt$ and $\selt \leq \telt$ and $\qelt \leq \relt \wedge \selt$ when $\qelt \leq \relt$ and $\qelt \leq \selt$. 
So, $L$ is necessarily connected. 
The lattice $L$ is {\em topographically balanced} if, for all $\relt$ and $\selt$ in $L$, (i) $\relt \leftarrow \qelt \rightarrow \selt$ for some $\qelt \in L$ $\Longrightarrow$ there exists a unique $\telt$ in $L$ such that $\relt \rightarrow \telt \leftarrow \selt$ and (ii) $\relt \rightarrow \telt \leftarrow \selt$ for some $\telt \in L$ $\Longrightarrow$ there exists a unique $\qelt$ in $L$ such that $\relt \leftarrow \qelt \rightarrow \selt$. 
A topographically balanced lattice is typically called a modular lattice. 
Such a lattice is necessarily ranked and has the property that $\rho(\relt \vee \selt) + \rho(\relt \wedge \selt) = \rho(\relt)+\rho(\selt)$. 
If $\relt \wedge (\selt \vee \telt) = (\relt \wedge \selt) \vee (\relt \wedge \telt)$ and $\relt \vee (\selt \wedge \telt) = (\relt \vee \selt) \wedge (\relt \vee \telt)$ for all $\relt, \selt, \telt$ in some lattice $L$, then we say $L$ is a distributive lattice. 
Any distributive lattice is modular. 
We say an edge-colored modular or distributive lattice $L$ is {\em diamond-colored} if, whenever \parbox{1.4cm}{\begin{center}
\setlength{\unitlength}{0.2cm}
\begin{picture}(5,3)
\put(2.5,0){\circle*{0.5}} \put(0.5,2){\circle*{0.5}}
\put(2.5,4){\circle*{0.5}} \put(4.5,2){\circle*{0.5}}
\put(0.5,2){\line(1,1){2}} \put(2.5,0){\line(-1,1){2}}
\put(4.5,2){\line(-1,1){2}} \put(2.5,0){\line(1,1){2}}
\put(1.25,0.55){\em \small k} \put(3.2,0.7){\em \small l}
\put(1.2,2.7){\em \small i} \put(3.25,2.55){\em \small j}
\end{picture} \end{center}} is a diamond of edges in $L$, then $i=l$ and $j=k$. 
A diamond-colored modular (respectively, distributive) lattice is a `DCML' (resp.\ `DCDL') for short. 

The {\em type} $\myA_{n}$ {\em skew-tabular lattices} $L_{\mytinyA_{n}}^{\mbox{\tiny skew}}(\mysmallP/\mysmallQ)$ of \cite{DD1} are a useful example of a family of diamond-colored distributive lattices.  
The non-increasing integer $m$-tuples $\mysmallP = (\mysmallP_{1},\ldots,\mysmallP_{m})$ and $\mysmallQ = (\mysmallQ_{1},\ldots,\mysmallQ_{m})$ are to be viewed as partitions with $\mysmallP_{i} \geq \mysmallQ_{i} \geq 0$ for each $1 \leq i \leq m$. Without losing generality, we may take $m > n$.  
As a set, $L_{\mytinyA_{n}}^{\mbox{\tiny skew}}(\mysmallP/\mysmallQ)$ is comprised of nonnegative integer arrays $(g_{i,j})$ where $i \in \{0,1,\ldots,n+1\}$ and $j \in C_{i} := \{i,i-1,\ldots,i-(m-1)\}$ satisfying $g_{i-1,j} \geq g_{i,j} \geq g_{i+1,j}$ and $g_{i-1,j-1} \leq g_{i,j} \leq g_{i+1,j+1}$ whenever any of the foregoing entries $g_{p,q}$ is part of our array.
It is advantageous to depict these arrays as parallelograms, as in \ParallelogramFigures, and we refer to these arrays as {\em GT $(n+1)$-parallelograms framed by} $\mysmallP/\mysmallQ$. 
The partial ordering on $L_{\mytinyA_{n}}^{\mbox{\tiny skew}}(\mysmallP/\mysmallQ)$ and the resulting colored covering relations are detailed in \EsixEsevenIntroFig. 
It is routine to verify that $L_{\mytinyA_{n}}^{\mbox{\tiny skew}}(\mysmallP/\mysmallQ)$ is a DCDL. 
 
Diamond-colored distributive lattices have some structural properties that make them especially nice to work with.  
Say an element $\velt$ of a diamond-colored distributive lattice $L$ is join irreducible if it covers exactly one other element $\uelt$ in $L$, and declare that $\mathscr{V}_{\mbox{\tiny color}}(\velt) := \mathscr{E}_{\mbox{\tiny color}}(\uelt \rightarrow \velt)$. 
Let the vertex-colored set $P$ of join irreducible elements of $L$ have the induced partial order. 
We use the notation $P = \mathbf{j}_{\mbox{\tiny color}}(L)$ and call $P$ the {\em compression poset} of $L$. 
On the other hand, given a vertex-colored poset $P$, a down-set (also called an ideal) $\mathcal{I}$ has the property that $v\in\mathcal{I}$ and $u \leq v$ in $P$ means that $u \in \mathcal{I}$. 
Let $L$ be the collection of down-sets of $P$, partially ordered by subset containment. 
It is easy to see that joins are unions and meets are intersections, so $L$ is distributive. 
Now, $\mathcal{I} \rightarrow \mathcal{J}$ in $L$ if and only if there exists a unique $v \in P$ such that $\mathcal{J} \setminus \mathcal{I} = \{v\}$; in this case, we set $\mathscr{E}_{\mbox{\tiny color}}(\mathcal{I} \rightarrow \mathcal{J}) := \mathscr{V}_{\mbox{\tiny color}}(v)$. 
It is easy now to see that $L$ is a diamond-colored distributive lattice, which we notate as $\mathbf{J}_{\mbox{\tiny color}}(P)$ and call the lattice of down-sets of $P$. 
Analogous to Birkhoff's famous Representation Theorem for distributive lattices, we have, for all vertex-colored posets $P$ and diamond-colored distributive lattices $L$:  $\mathbf{J}_{\mbox{\tiny color}}(\mathbf{j}_{\mbox{\tiny color}}(L)) \cong L$ and $\mathbf{j}_{\mbox{\tiny color}}(\mathbf{J}_{\mbox{\tiny color}}(P)) \cong P$. 

The above properties of DCDL's and their companion compression posets help to characterize what are known as minuscule lattices and minuscule posets, which we henceforth call, respectively, {\em minuscule splitting DCDL's} and {\em minuscule compression posets}. 
In the $\myE_{6}/\myE_{7}$ cases, minuscule compression posets and minuscule DCDL's provide the crucial framework for the `$k=1$' versions of our constructions. 
For $k>1$, the diamond-colored distributive lattices we utilize in our constructions are part of a more general family of lattices whose investigation in \cite{DD2} is effected by these advantageous structural properties.  

Assume now that $R$ is a ranked poset with edges colored by our set $I$ which indexes simple roots within $\Phi$ and Chevalley generators within $\mathfrak{g}$. 
For any $\relt \in R$, $\wt(\relt)$ is the weight $\sum_{i \in I} \mym_{i}(\relt)\omega_{i}$ and $\WGF(R;\myvarZ) := \sum_{\relt \in R}\myvarZ^{\smallwt(\relt)}$ is the weight generating function for $R$. 
Say $R$ is $\Phi${\em -structured} if $\wt(\relt)+\alpha_{i}=\wt(\selt)$ whenever $\relt \myarrow{i} \selt$ in $R$, i.e.\ for all $j \ne i$, we have $\mym_{j}(\relt)+\langle \alpha_{j},\alpha_{i}^{\vee} \rangle = \mym_{j}(\selt)$. 
If $R$ is $\Phi$-structured and $\WGF(R;\myvarZ)$ is $\mathscr{W}$-invariant, then we call $R$ a {\em splitting poset} for the $\mathscr{W}$-symmetric function $\WGF(R;\myvarZ)$. 

Suppose now that we are given a $\mathfrak{g}$-module $V$ and a weight basis $\{v_{\relt}\}_{\relt \in R}$ for $V$ indexed by some set of objects $R$. 
We can depict this basis using an edge-colored directed graph in the following way. 
For $\relt$ and $\selt$ in $R$, we place a directed edge $\relt \myarrow{i} \selt$ of color $i$ from $\relt$ to $\selt$ if at least one of $\myqX^{(i)}_{\selt,\relt}$ or $\myqY^{(i)}_{\relt,\selt}$ is nonzero when we write $\myqx_{i}.v_{\relt} = \sum_{\telt \in R}\myqX^{(i)}_{\telt,\relt}v_{\telt}$ and $\myqy_{i}.v_{\selt} = \sum_{\qelt \in R}\myqY^{(i)}_{\qelt,\selt}v_{\qelt}$. 
This edge-colored directed graph is the {\em supporting graph} for the given weight basis, and it is the {\em representation diagram} for the weight basis if, in addition, we attach the scalar pair $(\myqX^{(i)}_{\selt,\relt},\myqY^{(i)}_{\relt,\selt})$ to each edge $\relt \myarrow{i} \selt$. 
In \cite{DonSupp} it is observed that $R$ is a splitting poset for the $\mathscr{W}$-symmetric function $\mychar(V)$, so $R$ is necessarily $\Phi$-structured with $\WGF(R;\myvarZ) = \mychar(V)$. 

Using our combinatorial perspective, we can synthetically produce representation diagrams in the following way. 
For simplicity, assume $L$ is a DCML whose edges are colored by $I$. 
To each edge $\relt \myarrow{i} \selt$, attach a pair of (complex) scalars $(\myqX^{(i)}_{\selt,\relt},\myqY^{(i)}_{\relt,\selt})$, at least one of which is nonzero. 
For all $\pelt,\qelt \in L$, regard each of the scalars $\myqX^{(i)}_{\qelt,\pelt}$ and $\myqY^{(i)}_{\pelt,\qelt}$ to be zero if there is no color $i$ directed edge from $\pelt$ to $\qelt$. 
Let $V[L]$ be the vector space freely generated by the set of symbols $\{v_{\relt}\}_{\relt \in L}$. 
Define actions of $\myqx_{i}$ and $\myqy_{i}$ on $V[L]$ by the rules 
\begin{equation}
\myqx_{i}.v_{\relt} := \sum_{\telt : \relt \myarrow{i} \telt}\myqX^{(i)}_{\telt,\relt}v_{\telt}\hspace*{0.25in} \mbox{ and }\hspace*{0.25in} \myqy_{i}.v_{\selt} := \sum_{\qelt : \qelt \myarrow{i} \selt}\myqY^{(i)}_{\qelt,\selt}v_{\qelt}.
\end{equation}
We say that the scalar pairs assigned to $L$ satisfy the {\em diamond relations} if, for all $i,j \in I$ and all diamonds \parbox{1.4cm}{\begin{center}
\setlength{\unitlength}{0.2cm}
\begin{picture}(5,3)
\put(2.5,0){\circle*{0.5}} \put(0.5,2){\circle*{0.5}}
\put(2.5,4){\circle*{0.5}} \put(4.5,2){\circle*{0.5}}
\put(0.5,2){\line(1,1){2}} \put(2.5,0){\line(-1,1){2}}
\put(4.5,2){\line(-1,1){2}} \put(2.5,0){\line(1,1){2}}
\put(1.25,0.55){\em \small j} \put(3.2,0.7){\em \small i}
\put(1.2,2.7){\em \small i} \put(3.25,2.55){\em \small j}
\put(3,-0.75){\footnotesize $\qelt$} \put(5.25,1.75){\footnotesize $\selt$}
\put(3,4){\footnotesize $\telt$} \put(-1,1.75){\footnotesize $\relt$}
\end{picture} \end{center}}, we have 
\[\myqY^{(j)}_{\selt,\telt}\myqX^{(i)}_{\telt,\relt} = \myqX^{(i)}_{\selt,\qelt}\myqY^{(j)}_{\qelt,\relt}\] 
and satisfy the {\em crossing relations} if, for every $i \in I$ and $\relt \in R$, we have 
\[\mym_{i}(\relt) = \sum_{\qelt : \qelt \myarrow{i} \relt}\bigg(\myqX^{(i)}_{\relt,\qelt}\myqY^{(i)}_{\qelt,\relt}\bigg) - \sum_{\selt : \relt \myarrow{i} \selt}\bigg(\myqX^{(i)}_{\selt,\relt}\myqY^{(i)}_{\relt,\selt}\bigg).\]
Within this setting, the following result `combinatorializes' key aspects of the structure of $\mathfrak{g}$-modules. 
This result is a direct translation of Propositions 3.1 and 3.2 of \cite{DD1}.

\noindent
{\bf \CombinatorialSerre}\ \ {\sl Keep the set-up of the preceding paragraph. 
Let $d = |L|$, and let $\mathscr{B}$ be the basis $\{v_{\relt}\}_{\relt \in L}$ of $V[L]$. 
With generator actions as defined on $V[L]$ by equations (1) above, we have the following equivalent conditions: 
The diamond-colored modular lattice $L$ is $\Phi$-structured and the scalar pairs assigned to the edges of $L$ satisfy the diamond and crossing relations if and only if $V[L]$ is a $\mathfrak{g}$-module, $\mathscr{B}$ is a weight basis for $V[L]$, $L$ together with the specified scalar pairs is the representation diagram for this weight basis, $\myqh_{i}.v_{\relt} = \mym_{i}(\relt)v_{\relt}$ for all $i \in I$ and $\relt \in L$, and for each $i \in I$  the $d \times d$ representing matrices for $\myqx_{i}$ and $\myqy_{i}$ with respect to $\mathscr{B}$ are, respectively, $\big(\myqX^{(i)}_{\selt,\relt}\big)_{\relt,\selt \in L}$ and $\big(\myqY^{(i)}_{\relt,\selt}\big)_{\relt,\selt \in L}$. 
In this case, let $\lambda$ be the (dominant) weight of the unique maximal element of $L$. 
Then $V[L]$ is irreducible with highest weight $\lambda$ if and only if} $\WGF(L;\myvarZ) = \chi_{_{\lambda}}^{\Phi}$. 

As an example, the preceding result was applied in \cite{DD1} to the type $\myA_{n}$ skew-tabular lattices. 
In particular, we demonstrated that all type $\myA_{n}$ skew-tabular lattices are $\myA_{n}$-structured. 
Edge coefficients within a given skew-tabular lattice were defined in terms of the GT $(n+1)$-parallelograms comprising the lattices, and all diamond and crossing relations were verified. 
(However, the representations realized by skew-tabular lattices are not, in general, irreducible.) 

We close this section with some brief observations about so-called minuscule representations. 
(These are the topic of \cite{Green}.) 
When a dominant weight $\lambda$ is minuscule, many special simplifying properties are conferred upon the corresponding minuscule representation, including: The weights of all nonzero weight spaces comprise exactly the $\mathscr{W}$-orbit of $\lambda$; all nonzero weight spaces have dimension one; and there is, up to scaling, only one weight basis. 
It is easy to see, then, that a minuscule representation has exactly one supporting graph which coincides with the unique splitting poset for the associated Weyl symmetric function $\chi_{_{\lambda}}^{\Phi}$. 
One can readily discern from \cite{PrEur} that this unique supporting graph $\!\!$/$\!\!$ splitting poset is a DCDL. 
For a proof of this latter fact that uses only general principles, see Theorem 12.6 of \cite{DonDistributive}. 
Henceforth, we denote by $L_{\Phi}(\lambda)$ the minuscule splitting distributive lattice associated with a $\Phi$-minuscule dominant weight $\lambda$. 
In \cite{PrEur}, Proctor observed that if we attach the scalar pair $(1,1)$ to each edge of a minuscule splitting DCDL, we get the representation diagram for a weight basis of the associated minuscule representation. 
The vertex-colored poset $P_{\Phi}(\lambda) := \mathbf{j}_{\mbox{\tiny color}}(L_{\Phi}(\lambda))$ denotes the associated minuscule compression poset. 
In what follows, each minuscule splitting DCDL $L_{\mytinyE_{6}}(\omega_{1})$, $L_{\mytinyE_{6}}(\omega_{6})$, and $L_{\mytinyE_{7}}(\omega_{1})$ is the $k=1$ instance for families of diamond-colored distributive lattices that will be shown, via \CombinatorialSerre, to be representation diagrams for families of irreducible representations of $\mathfrak{g}(\myE_{6})$ and $\mathfrak{g}(\myE_{7})$. 
The minuscule compression posets here denoted $P_{\mytinyE_{6}}(\omega_{1'})$, $P_{\mytinyE_{6}}(\omega_{6'})$, and $P_{\mytinyE_{7}}(\omega_{1})$, which are depicted in \MinusculeDepictions, will play a large role in our representation constructions.

\vspace*{0.2in}
\noindent {\bf \S \ESection\ Our} $\myE_{6}$- {\bf and} $\myE_{7}$- {\bf polyminuscule lattices.}  
In this section we exactly specify the $\myE_{6}$- and $\myE_{7}$-colored directed graphs which are the setting for our representation constructions. 
These edge-colored directed graphs are actually covering digraphs for some diamond-colored distributive lattices that have (for the most part) appeared elsewhere in the literature. 
Our versions of these lattices -- which we call `$\myE_{6}$- and $\myE_{7}$-polyminuscule lattices' -- are presented so as to effect the representation constructions of the next section. 

Our numbering of the Dynkin diagram nodes for $\myE_{6}$ and $\myE_{7}$ is unconventional but helps make crucial connections with results from \cite{DD1}:
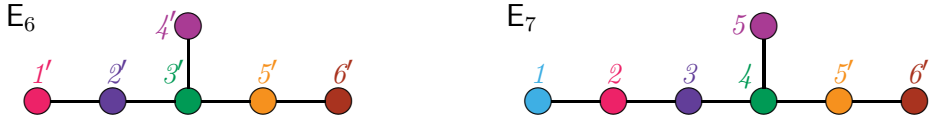
\begin{center}
\setlength{\unitlength}{0.5cm}
\begin{picture}(10,4)
\put(0.15,3){\large $\myE_{6}$}
{\color{WildStrawberry}\put(1,1){\circle*{0.7}}}
\put(1,1){\circle{0.7}}
{\color{RoyalPurple}\put(3,1){\circle*{0.7}}}
\put(3,1){\circle{0.7}}
{\color{ForestGreen}\put(5,1){\circle*{0.7}}}
\put(5,1){\circle{0.7}}
{\color{Mulberry}\put(5,3){\circle*{0.7}}}
\put(5,3){\circle{0.7}}
{\color{BurntOrange}\put(7,1){\circle*{0.7}}}
\put(7,1){\circle{0.7}}
{\color{Mahogany}\put(9,1){\circle*{0.7}}}
\put(9,1){\circle{0.7}}
\thicklines
\put(1.35,1){\line(1,0){1.3}}
\put(3.35,1){\line(1,0){1.3}}
\put(5.35,1){\line(1,0){1.3}}
\put(7.35,1){\line(1,0){1.3}}
\put(5,1.35){\line(0,1){1.3}}
\put(0.8,1.5){\textcolor{WildStrawberry}{\em 1$'$}}
\put(2.8,1.5){\textcolor{RoyalPurple}{\em 2$'$}}
\put(4.3,1.5){\textcolor{ForestGreen}{\em 3$'$}}
\put(4.1,2.8){\textcolor{Mulberry}{\em 4$'$}}
\put(6.8,1.5){\textcolor{BurntOrange}{\em 5$'$}}
\put(8.8,1.5){\textcolor{Mahogany}{\em 6$'$}}
\end{picture}
\hspace*{0.5in}
\setlength{\unitlength}{0.5cm}
\begin{picture}(12,4)
\put(0.15,3){\large $\myE_{7}$}
{\color{CornflowerBlue}\put(1,1){\circle*{0.7}}}
\put(1,1){\circle{0.7}}
{\color{WildStrawberry}\put(3,1){\circle*{0.7}}}
\put(3,1){\circle{0.7}}
{\color{RoyalPurple}\put(5,1){\circle*{0.7}}}
\put(5,1){\circle{0.7}}
{\color{ForestGreen}\put(7,1){\circle*{0.7}}}
\put(7,1){\circle{0.7}}
{\color{Mulberry}\put(7,3){\circle*{0.7}}}
\put(7,3){\circle{0.7}}
{\color{BurntOrange}\put(9,1){\circle*{0.7}}}
\put(9,1){\circle{0.7}}
{\color{Mahogany}\put(11,1){\circle*{0.7}}}
\put(11,1){\circle{0.7}}
\thicklines
\put(1.35,1){\line(1,0){1.3}}
\put(3.35,1){\line(1,0){1.3}}
\put(5.35,1){\line(1,0){1.3}}
\put(7.35,1){\line(1,0){1.3}}
\put(9.35,1){\line(1,0){1.3}}
\put(7,1.35){\line(0,1){1.3}}
\put(0.8,1.5){\textcolor{CornflowerBlue}{\em 1}}
\put(2.8,1.5){\textcolor{WildStrawberry}{\em 2}}
\put(4.8,1.5){\textcolor{RoyalPurple}{\em 3}}
\put(6.3,1.5){\textcolor{ForestGreen}{\em 4}}
\put(6.1,2.8){\textcolor{Mulberry}{\em 5}}
\put(8.8,1.5){\textcolor{BurntOrange}{\em 5$'$}}
\put(10.8,1.5){\textcolor{Mahogany}{\em 6$'$}}
\end{picture}
\end{center} 
Let $J_{5} := \{1,2,3,4,5\}$ and $J_{6} := \{1,2,3,4,5',6'\}$, both to be thought of as subsets of the $\myE_{7}$-indexing set $I_{7} := \{1,2,3,4,5,5',6'\}$. 
We regard $\mathfrak{g}(\myA_{5})$ to be the Lie subalgebra of $\mathfrak{g}(\myE_{7})$ generated by $\{\myqx_{i},\myqy_{i},\myqh_{i}\}_{i \in J_{5}}$, while $\mathfrak{g}(\myA_{6}) \hookrightarrow \mathfrak{g}(\myE_{7})$ is generated by $\{\myqx_{i},\myqy_{i},\myqh_{i}\}_{i \in J_{6}}$. 
Let $I_{6} := \{1',2',3',4',5',6'\}$, and let $\psi: I_{6} \longrightarrow \{2,3,4,5,5',6'\}$ be the correspondence of colors induced by viewing our $\myE_{6}$ Dynkin diagram as a subgraph of our $\myE_{7}$ Dynkin diagram. 
We regard $\mathfrak{g}(\myE_{6}) \hookrightarrow \mathfrak{g}(\myE_{7})$ to be generated by $\{\myqx_{\psi(i)},\myqy_{\psi(i)},\myqh_{\psi(i)}\}_{i \in I_{6}}$. 
A bit more bookkeeping: We view $\{0,1,2,3,4,5',6',7\}$ as totally ordered (with respect to the obvious ordering) with `$i-1$' preceding and `$i+1$' succeeding $i$ when $i \in J_{6}$. 

For $\myE_{6}$, the fundamental weights $\omega_{1'}$ and $\omega_{6'}$ are minuscule; for $\myE_{7}$, only $\omega_{1}$ is minuscule. 
In the $\myE_{7}$ case, associated to the minuscule weight $\omega_{1}$ is a vertex-colored poset $P_{\mytinyE_{7}}(\omega_{1})$ (see \PosetFigureEseven), which we call a minuscule compression poset, such that the corresponding diamond-colored distributive lattice $L_{\mytinyE_{7}}(\omega_{1})$ is a representation diagram for the associated minuscule representation when we take all edge coefficients to be unity. 
This observation seems to be original to \cite{PrEur}. 
Analogously, in the $\myE_{6}$ case there are (vertex-colored) minuscule compression posets $P_{\mytinyE_{6}}(\omega_{1'})$ and 

\newpage 
\newcommand{\TypeEboxDot}[1]{
\setlength{\unitlength}{0.5cm}
\begin{picture}(1,1)
\thicklines
\put(0.5,0.5){\color{#1}\circle*{0.7}}
\put(0.5,0.5){\circle{0.7}}
\end{picture}
}
\newcommand{\TypeEboxOneDigit}[3]{
\setlength{\unitlength}{0.4cm}
\begin{picture}(1.6,1)
\thicklines
\put(0,0){\color{#3}\line(0,1){1}}
\put(1.6,0){\color{#3}\line(0,1){1}}
\put(0,0){\color{#3}\line(1,0){1.6}}
\put(0,1){\color{#3}\line(1,0){1.6}}
\put(0.25,0.4){\tiny \color{#3}$\mbox{\footnotesize $c$}_{#1,#2}$}
\end{picture}
}
\newcommand{\TypeEboxTwoDigit}[3]{
\setlength{\unitlength}{0.4cm}
\begin{picture}(1.6,1)
\thicklines
\put(0,0){\color{#3}\line(0,1){1}}
\put(1.6,0){\color{#3}\line(0,1){1}}
\put(0,0){\color{#3}\line(1,0){1.6}}
\put(0,1){\color{#3}\line(1,0){1.6}}
\put(0.1,0.4){\tiny \color{#3}$\mbox{\footnotesize $c$}_{#1,#2}$}
\end{picture}
}
\newcommand{\TypeEboxOneDigitPrime}[3]{
\setlength{\unitlength}{0.4cm}
\begin{picture}(1.6,1)
\thicklines
\put(0,0){\color{#3}\line(0,1){1}}
\put(1.6,0){\color{#3}\line(0,1){1}}
\put(0,0){\color{#3}\line(1,0){1.6}}
\put(0,1){\color{#3}\line(1,0){1.6}}
\put(0.15,0.4){\tiny \color{#3}$\mbox{\footnotesize $c$}_{#1,#2}$}
\end{picture}
}
\newcommand{\TypeEboxTwoDigitPrime}[3]{
\setlength{\unitlength}{0.4cm}
\begin{picture}(1.6,1)
\thicklines
\put(0,0){\color{#3}\line(0,1){1}}
\put(1.8,0){\color{#3}\line(0,1){1}}
\put(0,0){\color{#3}\line(1,0){1.8}}
\put(0,1){\color{#3}\line(1,0){1.8}}
\put(0.1,0.4){\tiny \color{#3}$\mbox{\footnotesize $c$}_{#1,#2}$}
\end{picture}
}

\newpage
\begin{figure}
\begin{center}
{{\bf \PosetFigureEseven}\ \   The vertex-colored minuscule compression poset $P_{\mytinyE_{7}}(\omega_{1})$.}\\
{\scriptsize See \ArrayFigureEseven\ for a concrete description of the partial order on $P_{\mytinyE_{7}}(\omega_{1})$.}

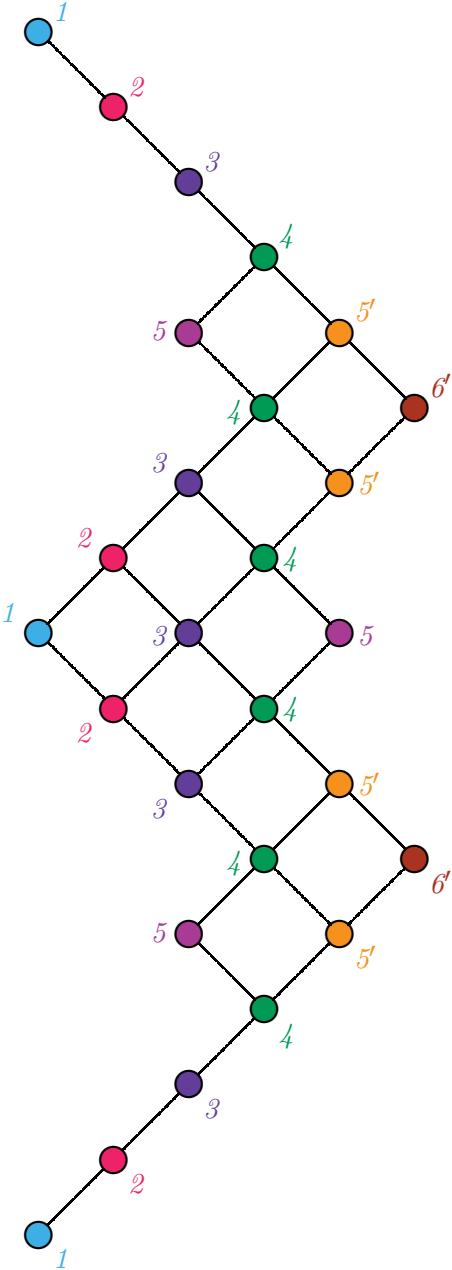
\setlength{\unitlength}{1cm}
\begin{picture}(6,17)
\put(0,16){\TypeEboxDot{CornflowerBlue}}
\put(0,8){\TypeEboxDot{CornflowerBlue}}
\put(0,0){\TypeEboxDot{CornflowerBlue}}
\put(0.7,-0.2){\textcolor{CornflowerBlue}{\em 1}}
\put(0,8.4){\textcolor{CornflowerBlue}{\em 1}}
\put(0.7,16.4){\textcolor{CornflowerBlue}{\em 1}}
\put(1,15){\TypeEboxDot{WildStrawberry}}
\put(1,9){\TypeEboxDot{WildStrawberry}}
\put(1,7){\TypeEboxDot{WildStrawberry}}
\put(1,1){\TypeEboxDot{WildStrawberry}}
\put(1.7,0.8){\textcolor{WildStrawberry}{\em 2}}
\put(1,6.8){\textcolor{WildStrawberry}{\em 2}}
\put(1,9.4){\textcolor{WildStrawberry}{\em 2}}
\put(1.7,15.4){\textcolor{WildStrawberry}{\em 2}}
\put(2,14){\TypeEboxDot{RoyalPurple}}
\put(2,10){\TypeEboxDot{RoyalPurple}}
\put(2,8){\TypeEboxDot{RoyalPurple}}
\put(2,6){\TypeEboxDot{RoyalPurple}}
\put(2,2){\TypeEboxDot{RoyalPurple}}
\put(2.7,1.8){\textcolor{RoyalPurple}{\em 3}}
\put(2,5.8){\textcolor{RoyalPurple}{\em 3}}
\put(2,8.1){\textcolor{RoyalPurple}{\em 3}}
\put(2,10.4){\textcolor{RoyalPurple}{\em 3}}
\put(2.7,14.4){\textcolor{RoyalPurple}{\em 3}}
\put(3,13){\TypeEboxDot{ForestGreen}}
\put(3,11){\TypeEboxDot{ForestGreen}}
\put(3,9){\TypeEboxDot{ForestGreen}}
\put(3,7){\TypeEboxDot{ForestGreen}}
\put(3,5){\TypeEboxDot{ForestGreen}}
\put(3,3){\TypeEboxDot{ForestGreen}}
\put(3.7,2.8){\textcolor{ForestGreen}{\em 4}}
\put(3,5.1){\textcolor{ForestGreen}{\em 4}}
\put(3.75,7.15){\textcolor{ForestGreen}{\em 4}}
\put(3.75,9.15){\textcolor{ForestGreen}{\em 4}}
\put(3,11.1){\textcolor{ForestGreen}{\em 4}}
\put(3.7,13.45){\textcolor{ForestGreen}{\em 4}}
\put(2,12){\TypeEboxDot{Mulberry}}
\put(4,8){\TypeEboxDot{Mulberry}}
\put(2,4){\TypeEboxDot{Mulberry}}
\put(2,4.15){\textcolor{Mulberry}{\em 5}}
\put(4.75,8.1){\textcolor{Mulberry}{\em 5}}
\put(2,12.15){\textcolor{Mulberry}{\em 5}}
\put(4,12){\TypeEboxDot{BurntOrange}}
\put(4,10){\TypeEboxDot{BurntOrange}}
\put(4,6){\TypeEboxDot{BurntOrange}}
\put(4,4){\TypeEboxDot{BurntOrange}}
\put(4.7,3.8){\textcolor{BurntOrange}{\em 5$'$}}
\put(4.75,6.1){\textcolor{BurntOrange}{\em 5$'$}}
\put(4.75,10.1){\textcolor{BurntOrange}{\em 5$'$}}
\put(4.7,12.4){\textcolor{BurntOrange}{\em 5$'$}}
\put(5,11){\TypeEboxDot{Mahogany}}
\put(5,5){\TypeEboxDot{Mahogany}}
\put(5.7,4.8){\textcolor{Mahogany}{\em 6$'$}}
\put(5.7,11.4){\textcolor{Mahogany}{\em 6$'$}}
\put(0.625,0.375){\qbezier(0,0)(0.375,0.375)(0.75,0.75)}
\put(1.625,1.375){\qbezier(0,0)(0.375,0.375)(0.75,0.75)}
\put(2.625,2.375){\qbezier(0,0)(0.375,0.375)(0.75,0.75)}
\put(3.625,3.375){\qbezier(0,0)(0.375,0.375)(0.75,0.75)}
\put(4.625,4.375){\qbezier(0,0)(0.375,0.375)(0.75,0.75)}
\put(2.625,4.375){\qbezier(0,0)(0.375,0.375)(0.75,0.75)}
\put(3.625,5.375){\qbezier(0,0)(0.375,0.375)(0.75,0.75)}
\put(2.625,6.375){\qbezier(0,0)(0.375,0.375)(0.75,0.75)}
\put(3.625,7.375){\qbezier(0,0)(0.375,0.375)(0.75,0.75)}
\put(1.625,7.375){\qbezier(0,0)(0.375,0.375)(0.75,0.75)}
\put(2.625,8.375){\qbezier(0,0)(0.375,0.375)(0.75,0.75)}
\put(3.625,9.375){\qbezier(0,0)(0.375,0.375)(0.75,0.75)}
\put(4.625,10.375){\qbezier(0,0)(0.375,0.375)(0.75,0.75)}
\put(0.625,8.375){\qbezier(0,0)(0.375,0.375)(0.75,0.75)}
\put(1.625,9.375){\qbezier(0,0)(0.375,0.375)(0.75,0.75)}
\put(2.625,10.375){\qbezier(0,0)(0.375,0.375)(0.75,0.75)}
\put(3.625,11.375){\qbezier(0,0)(0.375,0.375)(0.75,0.75)}
\put(2.625,12.375){\qbezier(0,0)(0.375,0.375)(0.75,0.75)}
\put(3.39,3.375){\qbezier(0,0)(-0.375,0.375)(-0.75,0.75)}
\put(4.39,4.375){\qbezier(0,0)(-0.375,0.375)(-0.75,0.75)}
\put(5.39,5.375){\qbezier(0,0)(-0.375,0.375)(-0.75,0.75)}
\put(3.39,5.375){\qbezier(0,0)(-0.375,0.375)(-0.75,0.75)}
\put(4.39,6.375){\qbezier(0,0)(-0.375,0.375)(-0.75,0.75)}
\put(2.39,6.375){\qbezier(0,0)(-0.375,0.375)(-0.75,0.75)}
\put(3.39,7.375){\qbezier(0,0)(-0.375,0.375)(-0.75,0.75)}
\put(4.39,8.375){\qbezier(0,0)(-0.375,0.375)(-0.75,0.75)}
\put(1.39,7.375){\qbezier(0,0)(-0.375,0.375)(-0.75,0.75)}
\put(2.39,8.375){\qbezier(0,0)(-0.375,0.375)(-0.75,0.75)}
\put(3.39,9.375){\qbezier(0,0)(-0.375,0.375)(-0.75,0.75)}
\put(4.39,10.375){\qbezier(0,0)(-0.375,0.375)(-0.75,0.75)}
\put(5.39,11.375){\qbezier(0,0)(-0.375,0.375)(-0.75,0.75)}
\put(3.39,11.375){\qbezier(0,0)(-0.375,0.375)(-0.75,0.75)}
\put(4.39,12.375){\qbezier(0,0)(-0.375,0.375)(-0.75,0.75)}
\put(3.39,13.375){\qbezier(0,0)(-0.375,0.375)(-0.75,0.75)}
\put(2.39,14.375){\qbezier(0,0)(-0.375,0.375)(-0.75,0.75)}
\put(1.39,15.375){\qbezier(0,0)(-0.375,0.375)(-0.75,0.75)}
\end{picture}
\end{center}
\end{figure}

\newpage
\begin{figure}
\begin{center}
{{\bf \ArrayFigureEseven}\ \   A naming of the positions of $P_{\mytinyE_{7}}(\omega_{1})$ for the purpose of building integer arrays.}

\vspace*{0.1in}
\parbox{6.25in}{\scriptsize Each vertex $v$ of $P_{\mytinyE_{7}}(\omega_{1})$ is uniquely identified by a pair $(p,q)$, where $p$ is the color of $v$ and $q$ is its rank. 
We can regard $P_{\mytinyE_{7}}(\omega_{1})$ to be a poset with partial order $(p,q) \leq (r,s)$ whenever $(p,q)$ and $(r,s)$ are the pairs identified with vertices of $P_{\mytinyE_{7}}(\omega_{1})$. 
We build integer arrays using the positions $P_{\mytinyE_{7}}(\omega_{1})$ as follows.  
The notation ``$c_{p,q}$'' identifies the array entry at position $(p,q)$ of $P_{\mytinyE_{7}}(\omega_{1})$.   
For $k \in \mathbb{Z}_{\geq 0}$, $L_{\mytinyE_{7}}(k\omega_{1})$ is the set integer arrays $\telt = \big(c_{p,q}(\telt)\big)$ such that $0 \leq c_{p,q}(\telt) \leq c_{r,s}(\telt) \leq k$ whenever $(p,q) \geq (r,s)$ in $P_{\mytinyE_{7}}(\omega_{1})$.}

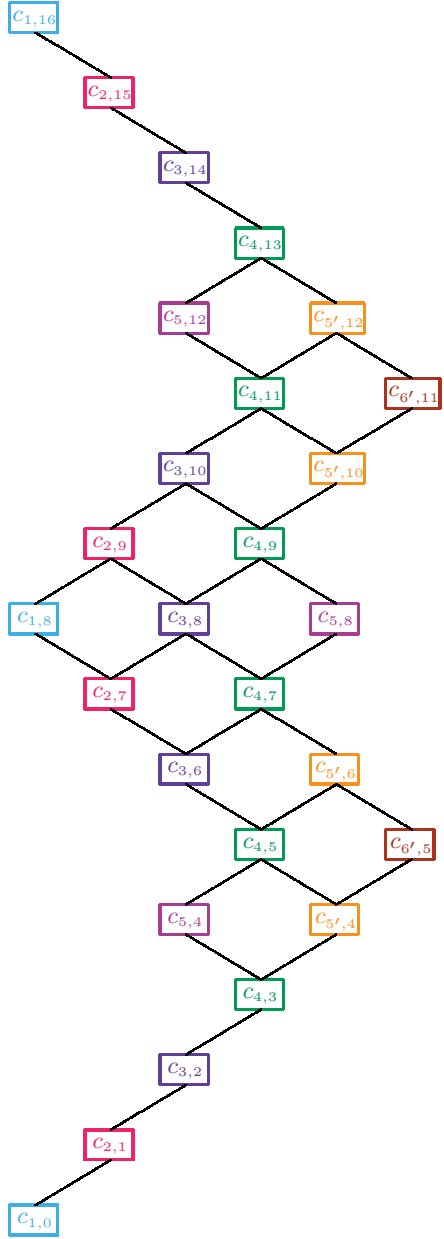
\setlength{\unitlength}{1cm}
\begin{picture}(6,17)
\put(0,16){\TypeEboxTwoDigit{1}{16}{CornflowerBlue}}
\put(0,8){\TypeEboxOneDigit{1}{8}{CornflowerBlue}}
\put(0,0){\TypeEboxOneDigit{1}{0}{CornflowerBlue}}
\put(1,15){\TypeEboxTwoDigit{2}{15}{WildStrawberry}}
\put(1,9){\TypeEboxOneDigit{2}{9}{WildStrawberry}}
\put(1,7){\TypeEboxOneDigit{2}{7}{WildStrawberry}}
\put(1,1){\TypeEboxOneDigit{2}{1}{WildStrawberry}}
\put(2,14){\TypeEboxTwoDigit{3}{14}{RoyalPurple}}
\put(2,10){\TypeEboxTwoDigit{3}{10}{RoyalPurple}}
\put(2,8){\TypeEboxOneDigit{3}{8}{RoyalPurple}}
\put(2,6){\TypeEboxOneDigit{3}{6}{RoyalPurple}}
\put(2,2){\TypeEboxOneDigit{3}{2}{RoyalPurple}}
\put(3,13){\TypeEboxTwoDigit{4}{13}{ForestGreen}}
\put(3,11){\TypeEboxTwoDigit{4}{11}{ForestGreen}}
\put(3,9){\TypeEboxOneDigit{4}{9}{ForestGreen}}
\put(3,7){\TypeEboxOneDigit{4}{7}{ForestGreen}}
\put(3,5){\TypeEboxOneDigit{4}{5}{ForestGreen}}
\put(3,3){\TypeEboxOneDigit{4}{3}{ForestGreen}}
\put(2,12){\TypeEboxTwoDigit{5}{12}{Mulberry}}
\put(4,8){\TypeEboxOneDigit{5}{8}{Mulberry}}
\put(2,4){\TypeEboxOneDigit{5}{4}{Mulberry}}
\put(4,12){\TypeEboxTwoDigitPrime{5'}{12}{BurntOrange}}
\put(4,10){\TypeEboxTwoDigitPrime{5'}{10}{BurntOrange}}
\put(4,6){\TypeEboxOneDigitPrime{5'}{6}{BurntOrange}}
\put(4,4){\TypeEboxOneDigitPrime{5'}{4}{BurntOrange}}
\put(5,11){\TypeEboxTwoDigitPrime{6'}{11}{Mahogany}}
\put(5,5){\TypeEboxOneDigitPrime{6'}{5}{Mahogany}}
\put(0.6,0.4){\qbezier(0,0)(0.5,0.3)(1,0.6)}
\put(1.6,1.4){\qbezier(0,0)(0.5,0.3)(1,0.6)}
\put(2.6,2.4){\qbezier(0,0)(0.5,0.3)(1,0.6)}
\put(3.6,3.4){\qbezier(0,0)(0.5,0.3)(1,0.6)}
\put(4.6,4.4){\qbezier(0,0)(0.5,0.3)(1,0.6)}
\put(2.6,4.4){\qbezier(0,0)(0.5,0.3)(1,0.6)}
\put(3.6,5.4){\qbezier(0,0)(0.5,0.3)(1,0.6)}
\put(2.6,6.4){\qbezier(0,0)(0.5,0.3)(1,0.6)}
\put(3.6,7.4){\qbezier(0,0)(0.5,0.3)(1,0.6)}
\put(1.6,7.4){\qbezier(0,0)(0.5,0.3)(1,0.6)}
\put(2.6,8.4){\qbezier(0,0)(0.5,0.3)(1,0.6)}
\put(3.6,9.4){\qbezier(0,0)(0.5,0.3)(1,0.6)}
\put(4.6,10.4){\qbezier(0,0)(0.5,0.3)(1,0.6)}
\put(0.6,8.4){\qbezier(0,0)(0.5,0.3)(1,0.6)}
\put(1.6,9.4){\qbezier(0,0)(0.5,0.3)(1,0.6)}
\put(2.6,10.4){\qbezier(0,0)(0.5,0.3)(1,0.6)}
\put(3.6,11.4){\qbezier(0,0)(0.5,0.3)(1,0.6)}
\put(2.6,12.4){\qbezier(0,0)(0.5,0.3)(1,0.6)}
\put(3.6,3.4){\qbezier(0,0)(-0.5,0.3)(-1,0.6)}
\put(4.6,4.4){\qbezier(0,0)(-0.5,0.3)(-1,0.6)}
\put(5.6,5.4){\qbezier(0,0)(-0.5,0.3)(-1,0.6)}
\put(3.6,5.4){\qbezier(0,0)(-0.5,0.3)(-1,0.6)}
\put(4.6,6.4){\qbezier(0,0)(-0.5,0.3)(-1,0.6)}
\put(2.6,6.4){\qbezier(0,0)(-0.5,0.3)(-1,0.6)}
\put(3.6,7.4){\qbezier(0,0)(-0.5,0.3)(-1,0.6)}
\put(4.6,8.4){\qbezier(0,0)(-0.5,0.3)(-1,0.6)}
\put(1.6,7.4){\qbezier(0,0)(-0.5,0.3)(-1,0.6)}
\put(2.6,8.4){\qbezier(0,0)(-0.5,0.3)(-1,0.6)}
\put(3.6,9.4){\qbezier(0,0)(-0.5,0.3)(-1,0.6)}
\put(4.6,10.4){\qbezier(0,0)(-0.5,0.3)(-1,0.6)}
\put(5.6,11.4){\qbezier(0,0)(-0.5,0.3)(-1,0.6)}
\put(3.6,11.4){\qbezier(0,0)(-0.5,0.3)(-1,0.6)}
\put(4.6,12.4){\qbezier(0,0)(-0.5,0.3)(-1,0.6)}
\put(3.6,13.4){\qbezier(0,0)(-0.5,0.3)(-1,0.6)}
\put(2.6,14.4){\qbezier(0,0)(-0.5,0.3)(-1,0.6)}
\put(1.6,15.4){\qbezier(0,0)(-0.5,0.3)(-1,0.6)}
\end{picture}
\end{center}
\end{figure}

\newpage
\begin{figure}
\begin{center}
{{\bf \PosetFiguresESix}\ \   The vertex-colored minuscule compression posets $P_{\mytinyE_{6}}(\omega_{1'})$ and $P_{\mytinyE_{6}}(\omega_{6'})$.}

\setlength{\unitlength}{1cm}
\begin{picture}(5,11)
\put(0,10){\TypeEboxDot{WildStrawberry}}
\put(0,4){\TypeEboxDot{WildStrawberry}}
\put(0,4.1){\textcolor{WildStrawberry}{\em 1$'$}}
\put(0.7,10.4){\textcolor{WildStrawberry}{\em 1$'$}}
\put(1,9){\TypeEboxDot{RoyalPurple}}
\put(1,5){\TypeEboxDot{RoyalPurple}}
\put(1,3){\TypeEboxDot{RoyalPurple}}
\put(1,3.1){\textcolor{RoyalPurple}{\em 2$'$}}
\put(1,5.1){\textcolor{RoyalPurple}{\em 2$'$}}
\put(1.7,9.4){\textcolor{RoyalPurple}{\em 2$'$}}
\put(2,8){\TypeEboxDot{ForestGreen}}
\put(2,6){\TypeEboxDot{ForestGreen}}
\put(2,4){\TypeEboxDot{ForestGreen}}
\put(2,2){\TypeEboxDot{ForestGreen}}
\put(2,2.1){\textcolor{ForestGreen}{\em 3$'$}}
\put(2.75,4.15){\textcolor{ForestGreen}{\em 3$'$}}
\put(2,6.1){\textcolor{ForestGreen}{\em 3$'$}}
\put(2.7,8.45){\textcolor{ForestGreen}{\em 3$'$}}
\put(1,7){\TypeEboxDot{Mulberry}}
\put(3,3){\TypeEboxDot{Mulberry}}
\put(3.75,3.1){\textcolor{Mulberry}{\em 4$'$}}
\put(1,7.15){\textcolor{Mulberry}{\em 4$'$}}
\put(3,7){\TypeEboxDot{BurntOrange}}
\put(3,5){\TypeEboxDot{BurntOrange}}
\put(3,1){\TypeEboxDot{BurntOrange}}
\put(3,1.1){\textcolor{BurntOrange}{\em 5$'$}}
\put(3.75,5.1){\textcolor{BurntOrange}{\em 5$'$}}
\put(3.7,7.4){\textcolor{BurntOrange}{\em 5$'$}}
\put(4,6){\TypeEboxDot{Mahogany}}
\put(4,0){\TypeEboxDot{Mahogany}}
\put(4,0.1){\textcolor{Mahogany}{\em 6$'$}}
\put(4.7,6.4){\textcolor{Mahogany}{\em 6$'$}}
\put(2.625,2.375){\qbezier(0,0)(0.375,0.375)(0.75,0.75)}
\put(1.625,3.375){\qbezier(0,0)(0.375,0.375)(0.75,0.75)}
\put(2.625,4.375){\qbezier(0,0)(0.375,0.375)(0.75,0.75)}
\put(3.625,5.375){\qbezier(0,0)(0.375,0.375)(0.75,0.75)}
\put(0.625,4.375){\qbezier(0,0)(0.375,0.375)(0.75,0.75)}
\put(1.625,5.375){\qbezier(0,0)(0.375,0.375)(0.75,0.75)}
\put(2.625,6.375){\qbezier(0,0)(0.375,0.375)(0.75,0.75)}
\put(1.625,7.375){\qbezier(0,0)(0.375,0.375)(0.75,0.75)}
\put(4.39,0.375){\qbezier(0,0)(-0.375,0.375)(-0.75,0.75)}
\put(3.39,1.375){\qbezier(0,0)(-0.375,0.375)(-0.75,0.75)}
\put(2.39,2.375){\qbezier(0,0)(-0.375,0.375)(-0.75,0.75)}
\put(3.39,3.375){\qbezier(0,0)(-0.375,0.375)(-0.75,0.75)}
\put(1.39,3.375){\qbezier(0,0)(-0.375,0.375)(-0.75,0.75)}
\put(2.39,4.375){\qbezier(0,0)(-0.375,0.375)(-0.75,0.75)}
\put(3.39,5.375){\qbezier(0,0)(-0.375,0.375)(-0.75,0.75)}
\put(4.39,6.375){\qbezier(0,0)(-0.375,0.375)(-0.75,0.75)}
\put(2.39,6.375){\qbezier(0,0)(-0.375,0.375)(-0.75,0.75)}
\put(3.39,7.375){\qbezier(0,0)(-0.375,0.375)(-0.75,0.75)}
\put(2.39,8.375){\qbezier(0,0)(-0.375,0.375)(-0.75,0.75)}
\put(1.39,9.375){\qbezier(0,0)(-0.375,0.375)(-0.75,0.75)}
\end{picture}
\hspace*{0.5in}
\begin{picture}(5,11)
\put(0,6){\TypeEboxDot{WildStrawberry}}
\put(0,0){\TypeEboxDot{WildStrawberry}}
\put(0.8,0.1){\textcolor{WildStrawberry}{\em 1$'$}}
\put(0,6.4){\textcolor{WildStrawberry}{\em 1$'$}}
\put(1,7){\TypeEboxDot{RoyalPurple}}
\put(1,5){\TypeEboxDot{RoyalPurple}}
\put(1,1){\TypeEboxDot{RoyalPurple}}
\put(1.8,1.1){\textcolor{RoyalPurple}{\em 2$'$}}
\put(1,5.1){\textcolor{RoyalPurple}{\em 2$'$}}
\put(1,7.4){\textcolor{RoyalPurple}{\em 2$'$}}
\put(2,8){\TypeEboxDot{ForestGreen}}
\put(2,6){\TypeEboxDot{ForestGreen}}
\put(2,4){\TypeEboxDot{ForestGreen}}
\put(2,2){\TypeEboxDot{ForestGreen}}
\put(2.8,2.1){\textcolor{ForestGreen}{\em 3$'$}}
\put(2,4.1){\textcolor{ForestGreen}{\em 3}}
\put(2.75,6.15){\textcolor{ForestGreen}{\em 3$'$}}
\put(2,8.4){\textcolor{ForestGreen}{\em 3$'$}}
\put(3,7){\TypeEboxDot{Mulberry}}
\put(1,3){\TypeEboxDot{Mulberry}}
\put(1,3.15){\textcolor{Mulberry}{\em 4$'$}}
\put(3.75,7.1){\textcolor{Mulberry}{\em 4$'$}}
\put(3,9){\TypeEboxDot{BurntOrange}}
\put(3,5){\TypeEboxDot{BurntOrange}}
\put(3,3){\TypeEboxDot{BurntOrange}}
\put(3.8,3.1){\textcolor{BurntOrange}{\em 5$'$}}
\put(3.75,5.1){\textcolor{BurntOrange}{\em 5$'$}}
\put(3,9.4){\textcolor{BurntOrange}{\em 5$'$}}
\put(4,10){\TypeEboxDot{Mahogany}}
\put(4,4){\TypeEboxDot{Mahogany}}
\put(4.8,4.1){\textcolor{Mahogany}{\em 6$'$}}
\put(4,10.4){\textcolor{Mahogany}{\em 6$'$}}
\put(0.625,0.375){\qbezier(0,0)(0.375,0.375)(0.75,0.75)}
\put(1.625,1.375){\qbezier(0,0)(0.375,0.375)(0.75,0.75)}
\put(2.625,2.375){\qbezier(0,0)(0.375,0.375)(0.75,0.75)}
\put(3.625,3.375){\qbezier(0,0)(0.375,0.375)(0.75,0.75)}
\put(1.625,3.375){\qbezier(0,0)(0.375,0.375)(0.75,0.75)}
\put(2.625,4.375){\qbezier(0,0)(0.375,0.375)(0.75,0.75)}
\put(1.625,5.375){\qbezier(0,0)(0.375,0.375)(0.75,0.75)}
\put(2.625,6.375){\qbezier(0,0)(0.375,0.375)(0.75,0.75)}
\put(0.625,6.375){\qbezier(0,0)(0.375,0.375)(0.75,0.75)}
\put(1.625,7.375){\qbezier(0,0)(0.375,0.375)(0.75,0.75)}
\put(2.625,8.375){\qbezier(0,0)(0.375,0.375)(0.75,0.75)}
\put(3.625,9.375){\qbezier(0,0)(0.375,0.375)(0.75,0.75)}
\put(2.39,2.375){\qbezier(0,0)(-0.375,0.375)(-0.75,0.75)}
\put(3.39,3.375){\qbezier(0,0)(-0.375,0.375)(-0.75,0.75)}
\put(4.39,4.375){\qbezier(0,0)(-0.375,0.375)(-0.75,0.75)}
\put(2.39,4.375){\qbezier(0,0)(-0.375,0.375)(-0.75,0.75)}
\put(3.39,5.375){\qbezier(0,0)(-0.375,0.375)(-0.75,0.75)}
\put(1.39,5.375){\qbezier(0,0)(-0.375,0.375)(-0.75,0.75)}
\put(2.39,6.375){\qbezier(0,0)(-0.375,0.375)(-0.75,0.75)}
\put(3.39,7.375){\qbezier(0,0)(-0.375,0.375)(-0.75,0.75)}
\end{picture}
\end{center}
\end{figure}

\newpage
\begin{figure}
\begin{center}
{{\bf \ArrayFiguresESix}\ \   A naming of the positions of $P_{\mytinyE_{6}}(\omega_{1'})$ and of $P_{\mytinyE_{6}}(\omega_{6'})$.}\\
{\scriptsize The notational conventions here follow \ArrayFigureEseven; the lattices $L_{\mytinyE_{6}}(k\omega_{1'})$ and $L_{\mytinyE_{6}}(k\omega_{6'})$ are built accordingly.}  

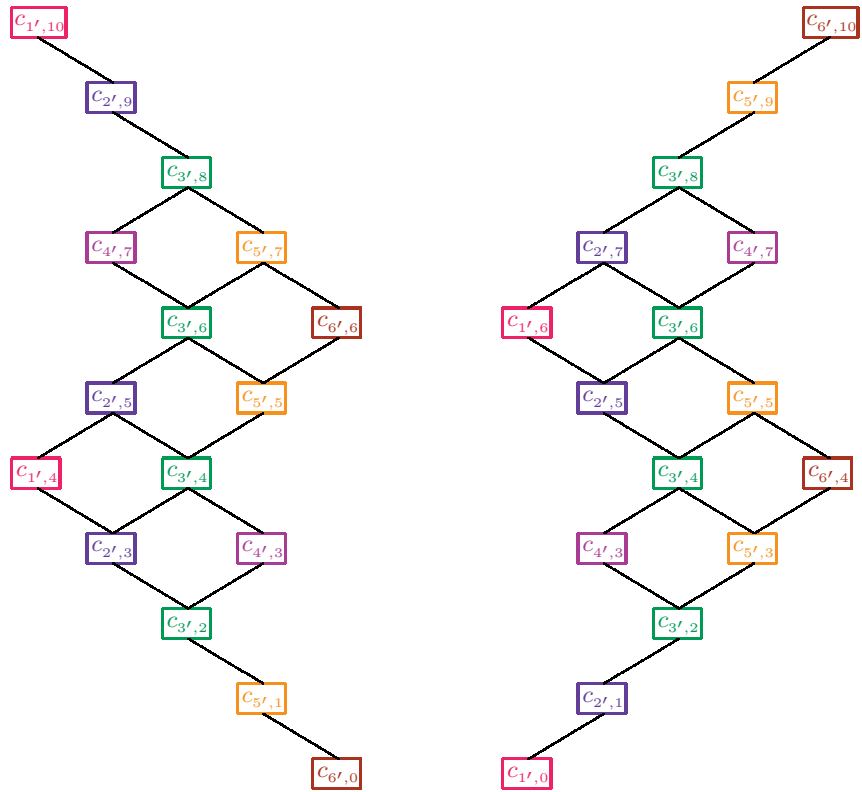
\setlength{\unitlength}{1cm}
\begin{picture}(5,11)
\put(0,10){\TypeEboxTwoDigitPrime{1'}{10}{WildStrawberry}}
\put(0,4){\TypeEboxOneDigitPrime{1'}{4}{WildStrawberry}}
\put(1,9){\TypeEboxOneDigitPrime{2'}{9}{RoyalPurple}}
\put(1,5){\TypeEboxOneDigitPrime{2'}{5}{RoyalPurple}}
\put(1,3){\TypeEboxOneDigitPrime{2'}{3}{RoyalPurple}}
\put(2,8){\TypeEboxOneDigitPrime{3'}{8}{ForestGreen}}
\put(2,6){\TypeEboxOneDigitPrime{3'}{6}{ForestGreen}}
\put(2,4){\TypeEboxOneDigitPrime{3'}{4}{ForestGreen}}
\put(2,2){\TypeEboxOneDigitPrime{3'}{2}{ForestGreen}}
\put(1,7){\TypeEboxOneDigitPrime{4'}{7}{Mulberry}}
\put(3,3){\TypeEboxOneDigitPrime{4'}{3}{Mulberry}}
\put(3,7){\TypeEboxOneDigitPrime{5'}{7}{BurntOrange}}
\put(3,5){\TypeEboxOneDigitPrime{5'}{5}{BurntOrange}}
\put(3,1){\TypeEboxOneDigitPrime{5'}{1}{BurntOrange}}
\put(4,6){\TypeEboxOneDigitPrime{6'}{6}{Mahogany}}
\put(4,0){\TypeEboxOneDigitPrime{6'}{0}{Mahogany}}
\put(2.6,2.4){\qbezier(0,0)(0.5,0.3)(1,0.6)}
\put(1.6,3.4){\qbezier(0,0)(0.5,0.3)(1,0.6)}
\put(2.6,4.4){\qbezier(0,0)(0.5,0.3)(1,0.6)}
\put(3.6,5.4){\qbezier(0,0)(0.5,0.3)(1,0.6)}
\put(0.6,4.4){\qbezier(0,0)(0.5,0.3)(1,0.6)}
\put(1.6,5.4){\qbezier(0,0)(0.5,0.3)(1,0.6)}
\put(2.6,6.4){\qbezier(0,0)(0.5,0.3)(1,0.6)}
\put(1.6,7.4){\qbezier(0,0)(0.5,0.3)(1,0.6)}
\put(4.6,0.4){\qbezier(0,0)(-0.5,0.3)(-1,0.6)}
\put(3.6,1.4){\qbezier(0,0)(-0.5,0.3)(-1,0.6)}
\put(2.6,2.4){\qbezier(0,0)(-0.5,0.3)(-1,0.6)}
\put(3.6,3.4){\qbezier(0,0)(-0.5,0.3)(-1,0.6)}
\put(1.6,3.4){\qbezier(0,0)(-0.5,0.3)(-1,0.6)}
\put(2.6,4.4){\qbezier(0,0)(-0.5,0.3)(-1,0.6)}
\put(3.6,5.4){\qbezier(0,0)(-0.5,0.3)(-1,0.6)}
\put(4.6,6.4){\qbezier(0,0)(-0.5,0.3)(-1,0.6)}
\put(2.6,6.4){\qbezier(0,0)(-0.5,0.3)(-1,0.6)}
\put(3.6,7.4){\qbezier(0,0)(-0.5,0.3)(-1,0.6)}
\put(2.6,8.4){\qbezier(0,0)(-0.5,0.3)(-1,0.6)}
\put(1.6,9.4){\qbezier(0,0)(-0.5,0.3)(-1,0.6)}
\end{picture}
\hspace*{0.5in}
\begin{picture}(5,11)
\put(0,6){\TypeEboxOneDigitPrime{1'}{6}{WildStrawberry}}
\put(0,0){\TypeEboxOneDigitPrime{1'}{0}{WildStrawberry}}
\put(1,7){\TypeEboxOneDigitPrime{2'}{7}{RoyalPurple}}
\put(1,5){\TypeEboxOneDigitPrime{2'}{5}{RoyalPurple}}
\put(1,1){\TypeEboxOneDigitPrime{2'}{1}{RoyalPurple}}
\put(2,8){\TypeEboxOneDigitPrime{3'}{8}{ForestGreen}}
\put(2,6){\TypeEboxOneDigitPrime{3'}{6}{ForestGreen}}
\put(2,4){\TypeEboxOneDigitPrime{3'}{4}{ForestGreen}}
\put(2,2){\TypeEboxOneDigitPrime{3'}{2}{ForestGreen}}
\put(3,7){\TypeEboxOneDigitPrime{4'}{7}{Mulberry}}
\put(1,3){\TypeEboxOneDigitPrime{4'}{3}{Mulberry}}
\put(3,9){\TypeEboxOneDigitPrime{5'}{9}{BurntOrange}}
\put(3,5){\TypeEboxOneDigitPrime{5'}{5}{BurntOrange}}
\put(3,3){\TypeEboxOneDigitPrime{5'}{3}{BurntOrange}}
\put(4,10){\TypeEboxTwoDigitPrime{6'}{10}{Mahogany}}
\put(4,4){\TypeEboxOneDigitPrime{6'}{4}{Mahogany}}
\put(0.6,0.4){\qbezier(0,0)(0.5,0.3)(1,0.6)}
\put(1.6,1.4){\qbezier(0,0)(0.5,0.3)(1,0.6)}
\put(2.6,2.4){\qbezier(0,0)(0.5,0.3)(1,0.6)}
\put(3.6,3.4){\qbezier(0,0)(0.5,0.3)(1,0.6)}
\put(1.6,3.4){\qbezier(0,0)(0.5,0.3)(1,0.6)}
\put(2.6,4.4){\qbezier(0,0)(0.5,0.3)(1,0.6)}
\put(1.6,5.4){\qbezier(0,0)(0.5,0.3)(1,0.6)}
\put(2.6,6.4){\qbezier(0,0)(0.5,0.3)(1,0.6)}
\put(0.6,6.4){\qbezier(0,0)(0.5,0.3)(1,0.6)}
\put(1.6,7.4){\qbezier(0,0)(0.5,0.3)(1,0.6)}
\put(2.6,8.4){\qbezier(0,0)(0.5,0.3)(1,0.6)}
\put(3.6,9.4){\qbezier(0,0)(0.5,0.3)(1,0.6)}
\put(2.6,2.4){\qbezier(0,0)(-0.5,0.3)(-1,0.6)}
\put(3.6,3.4){\qbezier(0,0)(-0.5,0.3)(-1,0.6)}
\put(4.6,4.4){\qbezier(0,0)(-0.5,0.3)(-1,0.6)}
\put(2.6,4.4){\qbezier(0,0)(-0.5,0.3)(-1,0.6)}
\put(3.6,5.4){\qbezier(0,0)(-0.5,0.3)(-1,0.6)}
\put(1.6,5.4){\qbezier(0,0)(-0.5,0.3)(-1,0.6)}
\put(2.6,6.4){\qbezier(0,0)(-0.5,0.3)(-1,0.6)}
\put(3.6,7.4){\qbezier(0,0)(-0.5,0.3)(-1,0.6)}
\end{picture}
\end{center}
\end{figure}

\newpage
\begin{figure}
\begin{center}
{{\bf \StackedFigureESix}\ \   A certain combination $\mathscr{P}$ of the vertex-colored posets $P_{\mytinyE_{6}}(\omega_{1'})$ and $P_{\mytinyE_{6}}(\omega_{6'})$.}

\setlength{\unitlength}{1cm}
\begin{picture}(5,15)
\put(0,14){\TypeEboxDot{WildStrawberry}}
\put(0,8){\TypeEboxDot{WildStrawberry}}
\put(0,6){\TypeEboxDot{WildStrawberry}}
\put(0,0){\TypeEboxDot{WildStrawberry}}
\put(0.7,-0.2){\textcolor{WildStrawberry}{\em 1$'$}}
\put(0,5.8){\textcolor{WildStrawberry}{\em 1$'$}}
\put(0,8.4){\textcolor{WildStrawberry}{\em 1$'$}}
\put(0.7,14.4){\textcolor{WildStrawberry}{\em 1$'$}}
\put(1,13){\TypeEboxDot{RoyalPurple}}
\put(1,9){\TypeEboxDot{RoyalPurple}}
\put(1,7){\TypeEboxDot{RoyalPurple}}
\put(1,5){\TypeEboxDot{RoyalPurple}}
\put(1,1){\TypeEboxDot{RoyalPurple}}
\put(1.7,0.8){\textcolor{RoyalPurple}{\em 2$'$}}
\put(1,4.8){\textcolor{RoyalPurple}{\em 2$'$}}
\put(1,7.1){\textcolor{RoyalPurple}{\em 2$'$}}
\put(1,9.4){\textcolor{RoyalPurple}{\em 2$'$}}
\put(1.7,13.4){\textcolor{RoyalPurple}{\em 2$'$}}
\put(2,12){\TypeEboxDot{ForestGreen}}
\put(2,10){\TypeEboxDot{ForestGreen}}
\put(2,8){\TypeEboxDot{ForestGreen}}
\put(2,6){\TypeEboxDot{ForestGreen}}
\put(2,4){\TypeEboxDot{ForestGreen}}
\put(2,2){\TypeEboxDot{ForestGreen}}
\put(2.7,1.8){\textcolor{ForestGreen}{\em 3$'$}}
\put(2,4.1){\textcolor{ForestGreen}{\em 3$'$}}
\put(2.75,6.15){\textcolor{ForestGreen}{\em 3$'$}}
\put(2.75,8.15){\textcolor{ForestGreen}{\em 3$'$}}
\put(2,10.1){\textcolor{ForestGreen}{\em 3$'$}}
\put(2.7,12.45){\textcolor{ForestGreen}{\em 3$'$}}
\put(1,11){\TypeEboxDot{Mulberry}}
\put(3,7){\TypeEboxDot{Mulberry}}
\put(1,3){\TypeEboxDot{Mulberry}}
\put(1,3.15){\textcolor{Mulberry}{\em 4$'$}}
\put(3.75,7.1){\textcolor{Mulberry}{\em 4$'$}}
\put(1,11.15){\textcolor{Mulberry}{\em 4$'$}}
\put(3,11){\TypeEboxDot{BurntOrange}}
\put(3,9){\TypeEboxDot{BurntOrange}}
\put(3,5){\TypeEboxDot{BurntOrange}}
\put(3,3){\TypeEboxDot{BurntOrange}}
\put(3.7,2.8){\textcolor{BurntOrange}{\em 5$'$}}
\put(3.75,5.1){\textcolor{BurntOrange}{\em 5$'$}}
\put(3.75,9.1){\textcolor{BurntOrange}{\em 5$'$}}
\put(3.7,11.4){\textcolor{BurntOrange}{\em 5$'$}}
\put(4,10){\TypeEboxDot{Mahogany}}
\put(4,4){\TypeEboxDot{Mahogany}}
\put(4.7,3.8){\textcolor{Mahogany}{\em 6$'$}}
\put(4.7,10.4){\textcolor{Mahogany}{\em 6$'$}}
\put(0.625,0.375){\qbezier(0,0)(0.375,0.375)(0.75,0.75)}
\put(1.625,1.375){\qbezier(0,0)(0.375,0.375)(0.75,0.75)}
\put(2.625,2.375){\qbezier(0,0)(0.375,0.375)(0.75,0.75)}
\put(3.625,3.375){\qbezier(0,0)(0.375,0.375)(0.75,0.75)}
\put(1.625,3.375){\qbezier(0,0)(0.375,0.375)(0.75,0.75)}
\put(2.625,4.375){\qbezier(0,0)(0.375,0.375)(0.75,0.75)}
\put(1.625,5.375){\qbezier(0,0)(0.375,0.375)(0.75,0.75)}
\put(2.625,6.375){\qbezier(0,0)(0.375,0.375)(0.75,0.75)}
\put(0.625,6.375){\qbezier(0,0)(0.375,0.375)(0.75,0.75)}
\put(1.625,7.375){\qbezier(0,0)(0.375,0.375)(0.75,0.75)}
\put(2.625,8.375){\qbezier(0,0)(0.375,0.375)(0.75,0.75)}
\put(3.625,9.375){\qbezier(0,0)(0.375,0.375)(0.75,0.75)}
\put(0.625,8.375){\qbezier(0,0)(0.375,0.375)(0.75,0.75)}
\put(1.625,9.375){\qbezier(0,0)(0.375,0.375)(0.75,0.75)}
\put(2.625,10.375){\qbezier(0,0)(0.375,0.375)(0.75,0.75)}
\put(1.625,11.375){\qbezier(0,0)(0.375,0.375)(0.75,0.75)}
\put(2.39,2.375){\qbezier(0,0)(-0.375,0.375)(-0.75,0.75)}
\put(3.39,3.375){\qbezier(0,0)(-0.375,0.375)(-0.75,0.75)}
\put(4.39,4.375){\qbezier(0,0)(-0.375,0.375)(-0.75,0.75)}
\put(2.39,4.375){\qbezier(0,0)(-0.375,0.375)(-0.75,0.75)}
\put(3.39,5.375){\qbezier(0,0)(-0.375,0.375)(-0.75,0.75)}
\put(1.39,5.375){\qbezier(0,0)(-0.375,0.375)(-0.75,0.75)}
\put(2.39,6.375){\qbezier(0,0)(-0.375,0.375)(-0.75,0.75)}
\put(3.39,7.375){\qbezier(0,0)(-0.375,0.375)(-0.75,0.75)}
\put(1.39,7.375){\qbezier(0,0)(-0.375,0.375)(-0.75,0.75)}
\put(2.39,8.375){\qbezier(0,0)(-0.375,0.375)(-0.75,0.75)}
\put(3.39,9.375){\qbezier(0,0)(-0.375,0.375)(-0.75,0.75)}
\put(4.39,10.375){\qbezier(0,0)(-0.375,0.375)(-0.75,0.75)}
\put(2.39,10.375){\qbezier(0,0)(-0.375,0.375)(-0.75,0.75)}
\put(3.39,11.375){\qbezier(0,0)(-0.375,0.375)(-0.75,0.75)}
\put(2.39,12.375){\qbezier(0,0)(-0.375,0.375)(-0.75,0.75)}
\put(1.39,13.375){\qbezier(0,0)(-0.375,0.375)(-0.75,0.75)}
\end{picture}
\end{center}
\end{figure}

\newpage
\begin{figure}
\begin{center}
{{\bf \StackedArrayESix}\ \   {\small A naming of the positions of $\mathscr{P}$ (cf.\ \StackedFigESix) for the purpose of building integer arrays.}}

\vspace*{0.1in}
\parbox{6.25in}{\scriptsize The notational conventions here largely follow \ArrayFigureEseven, but in this case arrays are constructed somewhat differently.  
Let $a,b \in \mathbb{Z}_{\geq 0}$.  
Define $L_{\mytinyE_{6}}(a\omega_{1'}+b\omega_{6'})$ to be the set integer arrays $\telt = \big(c_{p,q}(\telt)\big)$ such that $0 \leq c_{p,q}(\telt) \leq c_{r,s}(\telt) \leq a+b$ whenever $(p,q) \geq (r,s)$ in $\mathscr{P}$ and such that $\textcolor{WildStrawberry}{c_{1',8}}(\telt) \leq a$, $\textcolor{WildStrawberry}{c_{1',6}}(\telt) \geq a$.}

\setlength{\unitlength}{1cm}
\begin{picture}(5,15)
\put(0,14){\TypeEboxTwoDigitPrime{1'}{14}{WildStrawberry}}
\put(0,8){\TypeEboxOneDigitPrime{1'}{8}{WildStrawberry}}
\put(0,6){\TypeEboxOneDigitPrime{1'}{6}{WildStrawberry}}
\put(0,0){\TypeEboxOneDigitPrime{1'}{0}{WildStrawberry}}
\put(1,13){\TypeEboxTwoDigitPrime{2'}{13}{RoyalPurple}}
\put(1,9){\TypeEboxOneDigitPrime{2'}{9}{RoyalPurple}}
\put(1,7){\TypeEboxOneDigitPrime{2'}{7}{RoyalPurple}}
\put(1,5){\TypeEboxOneDigitPrime{2'}{5}{RoyalPurple}}
\put(1,1){\TypeEboxOneDigitPrime{2'}{1}{RoyalPurple}}
\put(2,12){\TypeEboxTwoDigitPrime{3'}{12}{ForestGreen}}
\put(2,10){\TypeEboxTwoDigitPrime{3'}{10}{ForestGreen}}
\put(2,8){\TypeEboxOneDigitPrime{3'}{8}{ForestGreen}}
\put(2,6){\TypeEboxOneDigitPrime{3'}{6}{ForestGreen}}
\put(2,4){\TypeEboxOneDigitPrime{3'}{4}{ForestGreen}}
\put(2,2){\TypeEboxOneDigitPrime{3'}{2}{ForestGreen}}
\put(1,11){\TypeEboxTwoDigitPrime{4'}{11}{Mulberry}}
\put(3,7){\TypeEboxOneDigitPrime{4'}{7}{Mulberry}}
\put(1,3){\TypeEboxOneDigitPrime{4'}{3}{Mulberry}}
\put(3,11){\TypeEboxTwoDigitPrime{5'}{11}{BurntOrange}}
\put(3,9){\TypeEboxOneDigitPrime{5'}{9}{BurntOrange}}
\put(3,5){\TypeEboxOneDigitPrime{5'}{5}{BurntOrange}}
\put(3,3){\TypeEboxOneDigitPrime{5'}{3}{BurntOrange}}
\put(4,10){\TypeEboxTwoDigitPrime{6'}{10}{Mahogany}}
\put(4,4){\TypeEboxOneDigitPrime{6'}{4}{Mahogany}}
\put(0.6,0.4){\qbezier(0,0)(0.5,0.3)(1,0.6)}
\put(1.6,1.4){\qbezier(0,0)(0.5,0.3)(1,0.6)}
\put(2.6,2.4){\qbezier(0,0)(0.5,0.3)(1,0.6)}
\put(3.6,3.4){\qbezier(0,0)(0.5,0.3)(1,0.6)}
\put(1.6,3.4){\qbezier(0,0)(0.5,0.3)(1,0.6)}
\put(2.6,4.4){\qbezier(0,0)(0.5,0.3)(1,0.6)}
\put(1.6,5.4){\qbezier(0,0)(0.5,0.3)(1,0.6)}
\put(2.6,6.4){\qbezier(0,0)(0.5,0.3)(1,0.6)}
\put(0.6,6.4){\qbezier(0,0)(0.5,0.3)(1,0.6)}
\put(1.6,7.4){\qbezier(0,0)(0.5,0.3)(1,0.6)}
\put(2.6,8.4){\qbezier(0,0)(0.5,0.3)(1,0.6)}
\put(3.6,9.4){\qbezier(0,0)(0.5,0.3)(1,0.6)}
\put(0.6,8.4){\qbezier(0,0)(0.5,0.3)(1,0.6)}
\put(1.6,9.4){\qbezier(0,0)(0.5,0.3)(1,0.6)}
\put(2.6,10.4){\qbezier(0,0)(0.5,0.3)(1,0.6)}
\put(1.6,11.4){\qbezier(0,0)(0.5,0.3)(1,0.6)}
\put(2.6,2.4){\qbezier(0,0)(-0.5,0.3)(-1,0.6)}
\put(3.6,3.4){\qbezier(0,0)(-0.5,0.3)(-1,0.6)}
\put(4.6,4.4){\qbezier(0,0)(-0.5,0.3)(-1,0.6)}
\put(2.6,4.4){\qbezier(0,0)(-0.5,0.3)(-1,0.6)}
\put(3.6,5.4){\qbezier(0,0)(-0.5,0.3)(-1,0.6)}
\put(1.6,5.4){\qbezier(0,0)(-0.5,0.3)(-1,0.6)}
\put(2.6,6.4){\qbezier(0,0)(-0.5,0.3)(-1,0.6)}
\put(3.6,7.4){\qbezier(0,0)(-0.5,0.3)(-1,0.6)}
\put(1.6,7.4){\qbezier(0,0)(-0.5,0.3)(-1,0.6)}
\put(2.6,8.4){\qbezier(0,0)(-0.5,0.3)(-1,0.6)}
\put(3.6,9.4){\qbezier(0,0)(-0.5,0.3)(-1,0.6)}
\put(4.6,10.4){\qbezier(0,0)(-0.5,0.3)(-1,0.6)}
\put(2.6,10.4){\qbezier(0,0)(-0.5,0.3)(-1,0.6)}
\put(3.6,11.4){\qbezier(0,0)(-0.5,0.3)(-1,0.6)}
\put(2.6,12.4){\qbezier(0,0)(-0.5,0.3)(-1,0.6)}
\put(1.6,13.4){\qbezier(0,0)(-0.5,0.3)(-1,0.6)}
\end{picture}
\end{center}
\end{figure}
\clearpage

\noindent 
$P_{\mytinyE_{6}}(\omega_{6'})$ -- see \PosetFiguresESix\ -- whose corresponding diamond-colored distributive lattices $L_{\mytinyE_{6}}(\omega_{1'})$ and $L_{\mytinyE_{6}}(\omega_{6'})$ are representation diagrams for the associated minuscule representations; again \cite{PrEur} is an original source. 

For the remainder of the paper, $k$ represents a fixed nonnegative integer. 
Our initial aim is to declare a diamond-colored distributive lattice -- the $\myE_{7}$-polyminuscule lattice $L_{\mytinyE_{7}}(k\omega_{1})$ -- that will serve as the combinatorial environment for our construction of the irreducible representation of $\mathfrak{g}(\myE_{7})$ with highest weight $k\omega_{1}$. 
As in \ArrayFigureEseven, we declare $L_{\mytinyE_{7}}(k\omega_{1})$ to be the set{\small 
\[\left\{\mbox{integer arrays } \telt = \big(c_{p,q}(\telt)\big)_{(p,q) \in P_{\mytinyE_{7}}(\omega_{1})}\, \rule[-3.5mm]{0.2mm}{8mm}\ 0 \leq c_{p,q}(\telt) \leq c_{r,s}(\telt) \leq k \mbox{ when } (r,s) \leq (p,q) \mbox{ in } P_{\mytinyE_{7}}(\omega_{1})\right\}.\]}Regard arrays in $L_{\mytinyE_{7}}(k\omega_{1})$ to be partially ordered by the rule $\selt \leq \telt$ if and only if $c_{p,q}(\selt) \leq c_{p,q}(\telt)$ for all positions $(p,q)$ in $P_{\mytinyE_{7}}(\omega_{1})$. 
If $\selt < \telt$ and there is no array $\relt$ such that $\selt < \relt < \telt$, then it is clear there is some position $(i,j)$ in $P_{\mytinyE_{7}}(\omega_{1})$ such that $c_{p,q}(\selt) = c_{p,q}(\telt)$ for all $(p,q) \not= (i,j)$ while $c_{i,j}(\selt) + 1 = c_{i,j}(\telt)$. 
In this case, $\selt$ is covered by $\telt$, and we assign color $i$ to this covering relation by writing $\selt \myarrow{i} \telt$. 
It is a straightforward exercise to verify that the edge-colored poset $L_{\mytinyE_{7}}(k\omega_{1})$ is a diamond-colored distributive lattice. 
Its maximal element $\melt$ has $c_{p,q}(\melt)=k$ at all positions $(p,q)$ in $P_{\mytinyE_{7}}(\omega_{1})$. 
Note that $\wt(\melt) = k\omega_{1}$. 

The $\myE_{6}$-polyminuscule lattices $L_{\mytinyE_{6}}(k\omega_{1'})$ and $L_{\mytinyE_{6}}(k\omega_{6'})$ are formed analogously, see \EsixMinusculeFigures. 
We can also realize these $\myE_{6}$ lattices as certain $\psi(I_{6})$-colored components of $L_{\mytinyE_{7}}(k\omega_{1})$; by Theorem 3.8 of \cite{DonDistributive}, these components are distributive sublattices of $L_{\mytinyE_{7}}(k\omega_{1})$. 
Now let $\melt'$ be the array from $L_{\mytinyE_{7}}(k\omega_{1})$ with \textcolor{CornflowerBlue}{$c_{1,16}$}$(\melt')=0$ but $c_{p,q}(\melt')=k$ for any position $(p,q)$ in $P_{\mytinyE_{7}}(\omega_{1})$ other than \textcolor{CornflowerBlue}{$(1,16)$}.  
It is evident that $\comp_{\psi(I_{6})}(\melt')$ is isomorphic to $L_{\mytinyE_{6}}(k\omega_{1'})^{\psi}$ and has $\melt'$ as its unique maximal element. 
Similarly, let $\melt'' \in L_{\mytinyE_{7}}(k\omega_{1})$ have $c_{p,q}(\melt'') = 0$ if $(p,q)$ is at or above position \textcolor{CornflowerBlue}{$(1,8)$} in $P_{\mytinyE_{7}}(\omega_{1})$ and $c_{p,q}(\melt'') = k$ if $(p,q)$ is at or below position \textcolor{Mahogany}{$(6',11)$}. 
One can readily see that $\comp_{\psi(I_{6})}(\melt'')$ is isomorphic to $L_{\mytinyE_{6}}(k\omega_{6'})^{\psi}$ and has $\melt''$ as its unique maximal element. 
 
There is one other family of $\myE_{6}$-polyminuscule lattices for us to consider, namely those of the form $L_{\mytinyE_{6}}(a\omega_{1'}+b\omega_{6'})$. 
These are obtained via \StackedFigures\ (understanding that the same conventions used above to build the $\myE_{7}$-polyminuscule lattices apply) and, like $L_{\mytinyE_{6}}(k\omega_{1'})$ and $L_{\mytinyE_{6}}(k\omega_{6'})$, can be realized as a $\psi(I_{6})$-component of some $\myE_{7}$-polyminuscule lattice. 
Let $a$ and $b$ be nonnegative integers with $a + b = k$, and let $\tilde{\melt}$ be the array from $L_{\mytinyE_{7}}(k\omega_{1})$ with \textcolor{CornflowerBlue}{$c_{1,16}$}$(\tilde{\melt})=0$, \textcolor{CornflowerBlue}{$c_{1,8}$}$(\tilde{\melt})=a$, and \textcolor{CornflowerBlue}{$c_{1,0}$}$(\tilde{\melt})=k$ and other entries determined as follows: 
Set $c_{p,q}(\tilde{\melt})=a$ at any position $(p,q)$ in the interval $[\textcolor{WildStrawberry}{(2,9)},\textcolor{WildStrawberry}{(2,15)}]$, which consists of $\textcolor{WildStrawberry}{(2,9)}$, $\textcolor{WildStrawberry}{(2,15)}$, and all positions in $P_{\mytinyE_{7}}(\omega_{1})$ between them; similarly, set $c_{p,q}(\tilde{\melt})=k$ for any position $(p,q)$ in the interval $[\textcolor{WildStrawberry}{(2,1)},\textcolor{Mahogany}{(6',11)}]$. 
One can see that $\comp_{\psi(I_{6})}(\tilde{\melt})$ is isomorphic to $L_{\mytinyE_{6}}(a\omega_{1'}+b\omega_{6'})$ and has $\tilde{\melt}$ as its unique maximal element. 

Of the polyminuscule lattices\footnote{In \cite{DD2} we use the more specific language `$\mysmallE_{7}$ prismatic minuscule lattice' to describe $L_{\mytinyE_{7}}(k\omega_{1})$, since its associated vertex-colored compression poset is naturally isomorphic to $P_{\mytinyE_{7}}(\omega_{1}) \times \mathbf{[k]}$, where $\mathbf{[k]}$ denotes a $k$-element chain. Similarly, each of $L_{\mytinyE_{6}}(k\omega_{1'})$, and $L_{\mytinyE_{6}}(k\omega_{6'})$ is an `$\mysmallE_{6}$ prismatic minuscule lattice'.} introduced above, only the $\myE_{6}$ lattices $L_{\mytinyE_{6}}(a\omega_{1'}+b\omega_{6'})$ are genuinely new. 
The lattices $L_{\mytinyE_{7}}(k\omega_{1})$, $L_{\mytinyE_{6}}(k\omega_{1'})$, and $L_{\mytinyE_{6}}(k\omega_{6'})$ seem to have first appeared in \cite{PrEur}. 
In that paper, Proctor, in collaboration with Stanley, connected these and related distributive lattices to Seshadri's `standard monomial theory' (\cite{Seshadri}, \cite{LMS}) in order to establish the now-canonical enumerative result that all minuscule compression posets are Gaussian. 
Moreover, identities {\sl (1)}--{\sl (3)} of the next result can be viewed as consequences of that combined work.\footnote{The arrays comprising these lattices -- called `$k$-multichains' in \cite{PrEur} -- and their distributive lattice ordering are evident from Proctor's paper.  
Less explicit is the natural edge-coloring and $\myE_{6} / \myE_{7}$-structure of these lattices.}
Identity {\sl (4)} is a consequence of results that will appear in \cite{DD2}. 

\noindent 
{\bf \EsixAndEsevenCharacter}\ \ {\sl Let $a$, $b$, and $k$ be nonnegative integers with $a+b = k$, and also let} 
$L_{k} := L_{\mytinyE_{7}}(k\omega_{1})$, $M_{k \cdot 1'} := L_{\mytinyE_{6}}(k\omega_{1'})$, $M_{k \cdot 6'} := L_{\mytinyE_{6}}(k\omega_{6'})$, {\sl and} $M_{a\cdot 1' + b\cdot 6'} := L_{\mytinyE_{6}}(a\omega_{1'}+b\omega_{6'})$.  
{\sl Then $L_{k}$ is} $\myE_{7}${\sl -structured, and its unique maximal element has weight $k\omega_{1}$. 
Each of $M_{k \cdot 1'}$, $M_{k \cdot 6'}$, and $M_{a\cdot 1' + b\cdot 6'}$ is} $\myE_{6}${\sl -structured with unique maximal element of respective weight $k\omega_{1'}$, $k\omega_{6'}$, $a\omega_{1'} + b\omega_{6'}$. 
Moreover:}
\begin{eqnarray*}
\parbox{3.35in}{{\sl (1) (Seshadri--Proctor--Stanley)}\hfill $\WGF(L_{k};\myvarZ)$} & = & \chi_{_{k\omega_{1}}}^{\mytinyE_{7}} = \mychar_{\mytinyE_{7}}(k\omega_{1})\\
\parbox{3.35in}{{\sl (2) (Seshadri--Proctor--Stanley)}\hfill $\WGF(M_{k \cdot 1'};\myvarZ)$} & = & \chi_{_{k\omega_{1'}}}^{\mytinyE_{6}} = \mychar_{\mytinyE_{6}}(k\omega_{1'})\\
\parbox{3.35in}{{\sl (3)  (Seshadri--Proctor--Stanley)}\hfill $\WGF(M_{k \cdot 6'};\myvarZ)$} & = & \chi_{_{k\omega_{6'}}}^{\mytinyE_{6}} = \mychar_{\mytinyE_{6}}(k\omega_{6'})\\
\parbox{3.35in}{{\sl (4) \cite{DD2}}\hfill $\WGF(M_{a\cdot 1' + b\cdot 6'};\myvarZ)$} & = & \chi_{_{a\omega_{1'}+b\omega_{6'}}}^{\mytinyE_{6}} = \mychar_{\mytinyE_{6}}(a\omega_{1'}+b\omega_{6'}).
\end{eqnarray*}

{\em Proof.} It follows from Theorem 13.2 of \cite{DonDistributive} that $L_{k}$ is $\myE_{7}$-structured. 
Each of $M_{k \cdot 1'}$, $M_{k \cdot 6'}$, and $M_{a\cdot 1' + b\cdot 6'}$ is a $\psi(I_{6})$-component of $L_{k}$ and so is $\myE_{6}$-structured. 
Alternatively, the fact that $L_{k}$ is $\myE_{7}$-structured (and hence that $M_{k \cdot 1'}$, $M_{k \cdot 6'}$, and $M_{a\cdot 1' + b\cdot 6'}$ are $\myE_{6}$-structured) follows from \JLemmas\ together with the observation that every $\{5,5'\}$-component of $L_{k}$ is the poset product of a $\{5\}$-component and a $\{5'\}$-component, and similarly for $\{5,6'\}$-components of $L_{k}$. 
The identities {\sl (1)}, {\sl (2)}, and {\sl (3)} are consequences of work from \cite{Seshadri} and \cite{PrEur}, as noted in the paragraph preceding the proposition statement.  
For different proofs of these identities, see Corollaries 9.4 and 9.7 of \cite{DonPosetModels}.  
Identity {\sl (4)} is to be established in the forthcoming paper \cite{DD2}.\hfill\QED

\vspace*{0.2in}
\noindent {\bf \S \ConstructionSection\ Our constructions of some irreducible} $\mathfrak{g}(\myE_{6})${\bf - and} $\mathfrak{g}(\myE_{7})${\bf -modules.}  
In the previous section we developed the precise combinatorial settings for our representation constructions, namely our $\myE_{6}$- and $\myE_{7}$-polyminuscule lattices.  
Our main objective now is to realize each of these lattices as a representation diagram for the appropriate representation of $\mathfrak{g}(\myE_{6})$ or $\mathfrak{g}(\myE_{7})$. 

To achieve this objective, we will attach certain pairs of coefficients to the edges of our polyminuscule lattices, define generator actions according to the equations of (1) in \S \Setup, and then confirm that \CombinatorialSerre\ applies.
Since each of our $\myE_{6}$-polyminuscule lattices arises as a component of some $\myE_{7}$-polyminuscule lattice, it suffices to carry out this process for $L_{\mytinyE_{7}}(k\omega_{1})$. 
Our choices of edge coefficients for $L_{\mytinyE_{7}}(k\omega_{1})$ are presented in \AllCoefficientFigures, although some prior discussion is needed to make those coefficient presentations easier to follow. 
The culmination of the work of this section is \ConstructionTheorem\ and \ConstructionCorollaries. 

Here is an outline of our approach. 
The key step in our constructions is to identify the $J_{5}$- and $J_{6}$-components of our $I_{7}$-edge-colored DCDL $L_{\mytinyE_{7}}(k\omega_{1})$ as, respectively, type $\myA_{5}$ and type $\myA_{6}$ skew-tabular lattices. 
So, within each $J_{5}$- and $J_{6}$-component of $L_{\mytinyE_{7}}(k\omega_{1})$, we will need to convert each array to a GT parallelogram; \AllCoefficientFigures\ will be central to this process. 
This will allow us to supply the edges of these components with (a version of) the edge coefficients from \cite{DD1}. 
Note that this will result in twice declaring coefficients on edges with colors from $J_{5} \cap J_{6} = \{1,2,3,4\}$ but only once declaring coefficients on edges with colors from $\{5,5',6'\}$. 
For well-definedness, we will need to check that our edge-coefficient definitions agree on edges with colors from $J_{5} \cap J_{6} = \{1,2,3,4\}$; this is done in \AgreeLemma. 
Well-definedness allows us to regard the $J_{5}$- and $J_{6}$-components of $L_{\mytinyE_{7}}(k\omega_{1})$ as (respectively) $\mathfrak{g}(\myA_{5})$ and $\mathfrak{g}(\myA_{6})$ representation diagrams, from which we immediately conclude that all crossing relations are satisfied and that diamond relations are satisfied for any diamonds with colors $\{i,j\} \subset I_{7}$ except for $\{5,5'\}$ and $\{5,6'\}$. 
Then, \CombinatorialSerre\ will apply to $L_{\mytinyE_{7}}(k\omega_{1})$ once we check that the diamond relations are satisfied on diamonds whose colors comprise the set $\{5,5'\}$ or $\{5,6'\}$. 
\ConstructionTheorem\ and \ConstructionCorollaries\ will then follow. 

Now for the particulars. 
To begin our analysis of the $J_{5}$-components of $L_{\mytinyE_{7}}(k\omega_{1})$, we pick some $\telt \in L_{\mytinyE_{7}}(k\omega_{1})$. 
As a set, $\comp_{J_{5}}(\telt)$ consists of those arrays from $L_{\mytinyE_{7}}(k\omega_{1})$ that have the same values as $\telt$ at all positions of colors \textcolor{BurntOrange}{$5'$} and \textcolor{Mahogany}{$6'$}. 
That is, within $\comp_{J_{5}}(\telt)$ we regard as fixed the array entries \textcolor{BurntOrange}{$c_{5',4}$}$:=$\textcolor{BurntOrange}{$c_{5',4}$}$(\telt)$, \textcolor{BurntOrange}{$c_{5',6}$}$:=$\textcolor{BurntOrange}{$c_{5',6}$}$(\telt)$, \textcolor{BurntOrange}{$c_{5',10}$}$:=$\textcolor{BurntOrange}{$c_{5',10}$}$(\telt)$, \textcolor{BurntOrange}{$c_{5',12}$}$:=$\textcolor{BurntOrange}{$c_{5',12}$}$(\telt)$, \textcolor{Mahogany}{$c_{6',5}$}$:=$\textcolor{Mahogany}{$c_{6',5}$}$(\telt)$, and \textcolor{Mahogany}{$c_{6',11}$}$:=$\textcolor{Mahogany}{$c_{6',11}$}$(\telt)$. 
Define integer partitions $\mysmallP$ and $\mysmallQ$ as the nine-tuples $\mysmallP := (\mbox{\textcolor{BurntOrange}{$c_{5',10}$}},\mbox{\textcolor{BurntOrange}{$c_{5',10}$}},\mbox{\textcolor{BurntOrange}{$c_{5',6}$}},\mbox{\textcolor{BurntOrange}{$c_{5',6}$}},k,k,k,k,k)$ and $\mysmallQ := (0,\mbox{\textcolor{BurntOrange}{$c_{5',12}$}},\mbox{\textcolor{BurntOrange}{$c_{5',12}$}},\mbox{\textcolor{BurntOrange}{$c_{5',12}$}},\mbox{\textcolor{BurntOrange}{$c_{5',12}$}},\mbox{\textcolor{BurntOrange}{$c_{5',4}$}},\mbox{\textcolor{BurntOrange}{$c_{5',4}$}},\mbox{\textcolor{BurntOrange}{$c_{5',4}$}},\mbox{\textcolor{BurntOrange}{$c_{5',4}$}})$.  
In \EsevenParallelogramFiguresAfive, we relate the arrays in the $J_{5}$-component $\comp_{J_{5}}(\telt)$ to the GT $5$-parallelograms framed by $\mysmallP/\mysmallQ$. 
The following lemma amounts to an observation. 

\noindent
{\bf \JfiveLemma}\ \ {\sl In the notation of the preceding paragraph and \EsevenParallelogramFiguresAfive, our $J_{5}$-component} $\comp_{J_{5}}(\telt)$ {\sl of} $L_{\mytinyE_{7}}(k\omega_{1})$ {\sl is isomorphic to the skew-tabular lattice} $L_{\mytinyA_{5}}^{\mbox{\tiny skew}}(\mysmallP/\mysmallQ)$ {\sl of \cite{DD1}, and hence, in the notation of \EsevenParallelogramFigureAfivePartOne,} $\displaystyle \mym_{i}(\telt) = \sum_{q=0}^{8}\left(\rule[-2.5mm]{0mm}{6.5mm}\left[\rule[-1.25mm]{0mm}{4.5mm}2g_{i,i-q}(\telt)-g_{i+1,i+1-q}(\telt)-g_{i-1,i-1-q}(\telt)\rule[-1.25mm]{0mm}{4.5mm}\right]\right)$ {\sl for each $i \in J_{5}$.}  

We similarly analyze the $J_{6}$-components of $L_{\mytinyE_{7}}(k\omega_{1})$. 
Again pick any $\telt' \in L_{\mytinyE_{7}}(k\omega_{1})$. 
The associated $J_{6}$-component of $L_{\mytinyE_{7}}(k\omega_{1})$ is the subset of $L_{\mytinyE_{7}}(k\omega_{1})$ consisting of those arrays that have the same values as $\telt'$ at all positions of color \textcolor{Mulberry}{$5$}. 
So, within $\comp_{J_{6}}(\telt')$ we regard as fixed the array entries \textcolor{Mulberry}{$c_{5,4}$}$:=$\textcolor{Mulberry}{$c_{5,4}$}$(\telt')$, \textcolor{Mulberry}{$c_{5,8}$}$:=$\textcolor{Mulberry}{$c_{5,8}$}$(\telt')$, and \textcolor{Mulberry}{$c_{5,12}$}$:=$\textcolor{Mulberry}{$c_{5,12}$}$(\telt')$. 
Define integer partitions $\mysmallP'$ and $\mysmallQ'$ as the nine-tuples $\mysmallP' = (\mbox{\textcolor{Mulberry}{$c_{5,8}$}},\mbox{\textcolor{Mulberry}{$c_{5,8}$}},\mbox{\textcolor{Mulberry}{$c_{5,8}$}},k,k,k,k,k,k)$ and $\mysmallQ' = (0,\mbox{\textcolor{Mulberry}{$c_{5,12}$}},\mbox{\textcolor{Mulberry}{$c_{5,12}$}},\mbox{\textcolor{Mulberry}{$c_{5,12}$}},\mbox{\textcolor{Mulberry}{$c_{5,12}$}},\mbox{\textcolor{Mulberry}{$c_{5,4}$}},\mbox{\textcolor{Mulberry}{$c_{5,4}$}},\mbox{\textcolor{Mulberry}{$c_{5,4}$}},\mbox{\textcolor{Mulberry}{$c_{5,4}$}})$.  
In \EsevenParallelogramFiguresAsix, we relate the arrays in the $J_{6}$-component $\comp_{J_{6}}(\telt')$ to the GT $6$-parallelograms framed by $\mysmallP'/\mysmallQ'$. 
We now observe that: 

\noindent
{\bf \JsixLemma}\ \ {\sl In the notation of the preceding paragraph and \EsevenParallelogramFiguresAsix, our $J_{6}$-component} $\comp_{J_{6}}(\telt')$ {\sl of} $L_{\mytinyE_{7}}(k\omega_{1})$ {\sl is isomorphic to} $L_{\mytinyA_{6}}^{\mbox{\tiny skew}}(\mysmallP'/\mysmallQ')$ {\sl of \cite{DD1}, and hence, in the notation of \EsevenParallelogramFigureAsixPartOne,} $\displaystyle \mym_{i}(\telt') = \sum_{q=0}^{8}\left(\rule[-2.5mm]{0mm}{6.5mm}\left[\rule[-1.25mm]{0mm}{4.5mm}2g_{i,i-q}(\telt')-g_{i+1,i+1-q}(\telt')-g_{i-1,i-1-q}(\telt')\rule[-1.25mm]{0mm}{4.5mm}\right]\right)$ {\sl for each $i \in J_{6}$.} 

In view of \JLemmas, we will, as needed, identify any $P_{\mytinyE_{7}}(\omega_{1})$-framed array $\telt$ from $L_{\mytinyE_{7}}(k\omega_{1})$\, using\, the\, coordinates\, of\, \EsevenParallelogramFiguresAfive\ or\, of\, \EsevenParallelogramFiguresAsix.\hfill 
In \CoefficientFiguresPartOne

\newpage
\begin{figure}[t] 
\begin{center}
{\small {\bf \EsevenParallelogramFigureAfivePartOne}\ \  A labelling of positions for arrays from our $J_{5}$-component $\comp_{J_{5}}(\telt)$ of \JfiveLemma\ when} 

\vspace*{-0.05in}
{\small viewed as GT 5-parallelograms framed by $\mysmallP/\mysmallQ$.  \EsevenParallelogramFigureAfivePartTwo\ gives the values in each of these positions.}
 
{\footnotesize The nonnegative integer entries below must weakly increase along NW-to-SE and NE-to-SW diagonals.}  

\setlength{\unitlength}{1in}
\begin{picture}(0,0)
\put(-3.3,-8){\line(0,1){8.9}}
\put(-3.3,-8){\line(1,0){6.6}}
\put(3.3,-8){\line(0,1){8.9}}
\put(-3.3,0.9){\line(1,0){6.6}}
\end{picture}

\vspace*{0.15in}
\begin{tabular}{ccccccc}
$g_{0,-8}$ & & & & & & \\
 & $g_{1,-7}$ & & & & & \\
 $g_{0,-7}$ & & $g_{2,-6}$ & & & & \\
 & $g_{1,-6}$ & & $g_{3,-5}$ & & & \\
 $g_{0,-6}$ & & $g_{2,-5}$ & & $g_{4,-4}$ & & \\
 & $g_{1,-5}$ & & $g_{3,-4}$ & & $g_{5,-3}$ & \\
 $g_{0,-5}$ & & $g_{2,-4}$ & & $g_{4,-3}$ & & $g_{6,-2}$ \\
 & $g_{1,-4}$ & & $g_{3,-3}$ & & $g_{5,-2}$ & \\
 $g_{0,-4}$ & & $g_{2,-3}$ & & $g_{4,-2}$ & & $g_{6,-1}$ \\
 & $g_{1,-3}$ & & $g_{3,-2}$ & & $g_{5,-1}$ & \\
 $g_{0,-3}$ & & $g_{2,-2}$ & & $g_{4,-1}$ & & $g_{6,0}$ \\
 & $g_{1,-2}$ & & $g_{3,-1}$ & & $g_{5,0}$ & \\
 $g_{0,-2}$ & & $g_{2,-1}$ & & $g_{4,0}$ & & $g_{6,1}$ \\
 & $g_{1,-1}$ & & $g_{3,0}$ & & $g_{5,1}$ & \\
 $g_{0,-1}$ & & $g_{2,0}$ & & $g_{4,1}$ & & $g_{6,2}$ \\
 & $g_{1,0}$ & & $g_{3,1}$ & & $g_{5,2}$ & \\
$g_{0,0}$ & & $g_{2,1}$ & & $g_{4,2}$ & & $g_{6,3}$ \\
 & $g_{1,1}$ & & $g_{3,2}$ & & $g_{5,3}$ & \\
 & & $g_{2,2}$ & & $g_{4,3}$ & & $g_{6,4}$ \\
 & & & $g_{3,3}$ & & $g_{5,4}$ & \\
 & & & & $g_{4,4}$ & & $g_{6,5}$ \\ 
 & & & & & $g_{5,5}$ & \\
 & & & & & & $g_{6,6}$  
 \end{tabular}
\end{center}
Suppose $\relt \myarrow{i} \selt$ for GT 5-parallelograms $\relt$ and $\selt$ within the skew-tabular lattice isomorphic to our $J_{5}$-component, with $g_{i,j}(\relt)+1=g_{i,j}(\selt)$.  
So, $i \in \{1,2,3,4,5\}$. 
Following \cite{DD1}, we set 
\[\myqP_{\relt,\selt}^{(i)} := -\ \frac{\mbox{$\displaystyle \prod_{p\in C_{i+1}}$} \left(\rule[-2.25mm]{0mm}{5.75mm}g_{i,j}-g_{i+1,p}+j-p\right)\mbox{$\displaystyle \prod_{p\in C_{i-1}}$} \left(\rule[-2.25mm]{0mm}{5.75mm}g_{i,j}-g_{i-1,p}+j-p-1\right)}{\mbox{$\displaystyle \prod_{p\in C_{i} \setminus \{j\}}$} \left(\rule[-2.25mm]{0mm}{5.75mm}g_{i,j}-g_{i,p}+j-p-1\right)\left(\rule[-2.25mm]{0mm}{5.75mm}g_{i,j}-g_{i,p}+j-p\right)},\]
where we regard each $g_{p,q}$ to be the array entry $g_{p,q}(\selt)$ for the GT 5-parallelogram $\selt$ and where $C_{p}$ denotes the set of indices that are valid in the $p^{\mbox{\tiny th}}$ column of the array.  (So, for example, $C_{2} = \{2,1,0,-1,-2,-3,-4,-5,-6\}$.)  
It is not hard to see that $\myqP_{\relt,\selt}^{(i)}$ is a positive rational number. 
\end{figure}

\newpage
\begin{figure}[t] 
\begin{center}
{\small {\bf \EsevenParallelogramFigureAfivePartTwo}\ \  Array values at each of the GT parallelogram positions from \EsevenParallelogramFigureAfivePartOne.} 

{\footnotesize For this $J_{5}$-component, we regard \textcolor{BurntOrange}{$c_{5',4}$}, \textcolor{BurntOrange}{$c_{5',6}$}, \textcolor{BurntOrange}{$c_{5',10}$}, \textcolor{BurntOrange}{$c_{5',12}$}, \textcolor{Mahogany}{$c_{6',5}$}, and \textcolor{Mahogany}{$c_{6',11}$} to be fixed values.}

\setlength{\unitlength}{1in}
\begin{picture}(0,0)
\put(-3.3,-6.7){\line(0,1){7.4}}
\put(-3.3,-6.7){\line(1,0){6.6}}
\put(3.3,-6.7){\line(0,1){7.4}}
\put(-3.3,0.7){\line(1,0){6.6}}
\end{picture}

\vspace*{0.15in}
\begin{tabular}{ccccccc}
$0$ & & & & & & \\
 & \textcolor{CornflowerBlue}{$c_{1,16}$} & & & & & \\
 \textcolor{BurntOrange}{$c_{5',12}$} & & \textcolor{WildStrawberry}{$c_{2,15}$} & & & & \\
 & \textcolor{BurntOrange}{$c_{5',12}$} & & \textcolor{RoyalPurple}{$c_{3,14}$} & & & \\
 \textcolor{BurntOrange}{$c_{5',12}$} & & \textcolor{BurntOrange}{$c_{5',12}$} & & \textcolor{ForestGreen}{$c_{4,13}$} & & \\
 & \textcolor{BurntOrange}{$c_{5',12}$} & & \textcolor{BurntOrange}{$c_{5',12}$} & & \textcolor{Mulberry}{$c_{5,12}$} & \\
 \textcolor{BurntOrange}{$c_{5',12}$} & & \textcolor{BurntOrange}{$c_{5',12}$} & & \textcolor{ForestGreen}{$c_{4,11}$} & & \textcolor{BurntOrange}{$c_{5',10}$} \\
 & \textcolor{BurntOrange}{$c_{5',12}$} & & \textcolor{RoyalPurple}{$c_{3,10}$} & & \textcolor{BurntOrange}{$c_{5',10}$} & \\
 \textcolor{BurntOrange}{$c_{5',12}$} & &\textcolor{WildStrawberry}{$c_{2,9}$} & & \textcolor{ForestGreen}{$c_{4,9}$} & & \textcolor{BurntOrange}{$c_{5',10}$} \\
 & \textcolor{CornflowerBlue}{$c_{1,8}$} & & \textcolor{RoyalPurple}{$c_{3,8}$} & & \textcolor{Mulberry}{$c_{5,8}$} & \\
\textcolor{BurntOrange}{$c_{5',4}$} & & \textcolor{WildStrawberry}{$c_{2,7}$} & & \textcolor{ForestGreen}{$c_{4,7}$} & & \textcolor{BurntOrange}{$c_{5',6}$} \\
 & \textcolor{BurntOrange}{$c_{5',4}$} & & \textcolor{RoyalPurple}{$c_{3,6}$} & & \textcolor{BurntOrange}{$c_{5',6}$} & \\
 \textcolor{BurntOrange}{$c_{5',4}$} & & \textcolor{BurntOrange}{$c_{5',4}$}  & & \textcolor{ForestGreen}{$c_{4,5}$} & & \textcolor{BurntOrange}{$c_{5',6}$} \\
 & \textcolor{BurntOrange}{$c_{5',4}$} & & \textcolor{BurntOrange}{$c_{5',4}$} & & \textcolor{Mulberry}{$c_{5,4}$} & \\
 \textcolor{BurntOrange}{$c_{5',4}$} & & \textcolor{BurntOrange}{$c_{5',4}$} & & \textcolor{ForestGreen}{$c_{4,3}$} & & $k$ \\
 & \textcolor{BurntOrange}{$c_{5',4}$} & & \textcolor{RoyalPurple}{$c_{3,2}$} & & $k$ & \\
\textcolor{BurntOrange}{$c_{5',4}$} & & \textcolor{WildStrawberry}{$c_{2,1}$} & & $k$ & & $k$ \\
 & \textcolor{CornflowerBlue}{$c_{1,0}$} & & $k$ & & $k$ & \\
 & & $k$ & & $k$ & & $k$ \\
 & & & $k$ & & $k$ & \\
 & & & & $k$ & & $k$ \\ 
 & & & & & $k$ & \\
 & & & & & & $k$  
 \end{tabular}
\end{center}
If $\relt \myarrow{i} \selt$ for arrays in our $J_{5}$-component of $L_{\mytinyE_{7}}(k\omega_{1})$, then we must have $g_{i,j}(\relt) + 1 = g_{i,j}(\selt)$ where $(i,j)$ is not one of the above-identified fixed-value positions.  Then we set
\[\myqX_{\selt,\relt}^{(i)} := \sqrt{\myqP_{\relt,\selt}^{(i)}\, } =: \myqY_{\relt,\selt}^{(i)}.\]
\end{figure}

\newpage 
\begin{figure}[t] 
\begin{center}
{\small {\bf \EsevenParallelogramFigureAsixPartOne}\ \  A labelling of positions for arrays from our $J_{6}$-component $\comp_{J_{6}}(\telt')$ of \JsixLemma\ when} 

\vspace*{-0.05in}
{\small viewed as GT 6-parallelograms framed by $\mysmallP'/\mysmallQ'$.  \EsevenParallelogramFigureAsixPartTwo\ gives the values in each of these positions.}

{\footnotesize The nonnegative integer entries below must weakly increase along NW-to-SE and NE-to-SW diagonals.}

\setlength{\unitlength}{1in}
\begin{picture}(0,0)
\put(-3.3,-8.25){\line(0,1){9.15}}
\put(-3.3,-8.25){\line(1,0){6.6}}
\put(3.3,-8.25){\line(0,1){9.15}}
\put(-3.3,0.9){\line(1,0){6.6}}
\end{picture}

\vspace*{0.15in}
\begin{tabular}{cccccccc}
$g_{0,-8}$ & & & & & & & \\
 & $g_{1,-7}$ & & & & & & \\
 $g_{0,-7}$ & & $g_{2,-6}$ & & & & & \\
 & $g_{1,-6}$ & & $g_{3,-5}$ & & & & \\
 $g_{0,-6}$ & & $g_{2,-5}$ & & $g_{4,-4}$ & & & \\
 & $g_{1,-5}$ & & $g_{3,-4}$ & & $g_{5',-3}$ & & \\
 $g_{0,-5}$ & & $g_{2,-4}$ & & $g_{4,-3}$ & & $g_{6',-2}$ & \\
 & $g_{1,-4}$ & & $g_{3,-3}$ & & $g_{5',-2}$ & & $g_{7,-1}$ \\
 $g_{0,-4}$ & & $g_{2,-3}$ & & $g_{4,-2}$ & & $g_{6',-1}$ & \\
 & $g_{1,-3}$ & & $g_{3,-2}$ & & $g_{5',-1}$ & & $g_{7,0}$ \\
 $g_{0,-3}$ & & $g_{2,-2}$ & & $g_{4,-1}$ & & $g_{6',0}$ & \\
 & $g_{1,-2}$ & & $g_{3,-1}$ & & $g_{5',0}$ & & $g_{7,1}$ \\
 $g_{0,-2}$ & & $g_{2,-1}$ & & $g_{4,0}$ & & $g_{6',1}$ & \\
 & $g_{1,-1}$ & & $g_{3,0}$ & & $g_{5',1}$ & & $g_{7,2}$ \\
 $g_{0,-1}$ & & $g_{2,0}$ & & $g_{4,1}$ & & $g_{6',2}$ & \\
 & $g_{1,0}$ & & $g_{3,1}$ & & $g_{5',2}$ & & $g_{7,3}$ \\
 $g_{0,0}$ & & $g_{2,1}$ & & $g_{4,2}$ & & $g_{6',3}$ & \\
 & $g_{1,1}$ & & $g_{3,2}$ & & $g_{5',3}$ & & $g_{7,4}$ \\
 & & $g_{2,2}$ & & $g_{4,3}$ & & $g_{6',4}$ & \\
 & & & $g_{3,3}$ & & $g_{5',4}$ & & $g_{7,5}$ \\
 & & & & $g_{4,4}$ & & $g_{6',5}$ & \\
 & & & & & $g_{5',5}$ & & $g_{7,6}$ \\
 & & & & & & $g_{6',6}$ & \\
 & & & & & & & $g_{7,7}$
\end{tabular}
\end{center}
Suppose $\relt' \myarrow{i} \selt'$ for GT 6-parallelograms $\relt'$ and $\selt'$ within the skew-tabular lattice isomorphic to our $J_{6}$-component, with $g_{i,j}(\relt')+1=g_{i,j}(\selt')$.  
So, $i \in \{1,2,3,4,5',6'\}$. 
Following \cite{DD1}, we set 
\[\myqQ_{\relt',\selt'}^{(i)} := -\ \frac{\mbox{$\displaystyle \prod_{p\in C_{i+1}}$} \left(\rule[-2.25mm]{0mm}{5.75mm}g_{i,j}-g_{i+1,p}+j-p\right)\mbox{$\displaystyle \prod_{p\in C_{i-1}}$} \left(\rule[-2.25mm]{0mm}{5.75mm}g_{i,j}-g_{i-1,p}+j-p-1\right)}{\mbox{$\displaystyle \prod_{p\in C_{i} \setminus \{j\}}$} \left(\rule[-2.25mm]{0mm}{5.75mm}g_{i,j}-g_{i,p}+j-p-1\right)\left(\rule[-2.25mm]{0mm}{5.75mm}g_{i,j}-g_{i,p}+j-p\right)},\]
where we regard each $g_{p,q}$ to be the array entry $g_{p,q}(\selt')$ for the GT 6-parallelogram $\selt'$ and where $C_{p}$ denotes the set of indices that are valid in the $p^{\mbox{\tiny th}}$ column of the array.  (So, for example, $C_{5'} = \{5,4,3,2,1,0,-1,-2,-3\}$.)  
It is not hard to see that $\myqQ_{\relt',\selt'}^{(i)}$ is a positive rational number. 
\end{figure}

\newpage 
\begin{figure}[t] 
\begin{center}
{\small {\bf \EsevenParallelogramFigureAsixPartTwo}\ \  Array values at each of the GT parallelogram positions from \EsevenParallelogramFigureAsixPartOne.}

{\footnotesize For this $J_{6}$-component, we regard \textcolor{Mulberry}{$c_{5,4}$}, \textcolor{Mulberry}{$c_{5,8}$}, and \textcolor{Mulberry}{$c_{5,12}$} to be fixed values.}

\setlength{\unitlength}{1in}
\begin{picture}(0,0)
\put(-3.3,-7){\line(0,1){7.7}}
\put(-3.3,-7){\line(1,0){6.6}}
\put(3.3,-7){\line(0,1){7.7}}
\put(-3.3,0.7){\line(1,0){6.6}}
\end{picture}

\vspace*{0.15in}
\begin{tabular}{cccccccc}
$0$ & & & & & & & \\
 & \textcolor{CornflowerBlue}{$c_{1,16}$} & & & & & & \\
 \textcolor{Mulberry}{$c_{5,12}$} & & \textcolor{WildStrawberry}{$c_{2,15}$} & & & & & \\
 & \textcolor{Mulberry}{$c_{5,12}$} & & \textcolor{RoyalPurple}{$c_{3,14}$} & & & & \\
 \textcolor{Mulberry}{$c_{5,12}$} & & \textcolor{Mulberry}{$c_{5,12}$} & & \textcolor{ForestGreen}{$c_{4,13}$} & & & \\
 & \textcolor{Mulberry}{$c_{5,12}$} & & \textcolor{Mulberry}{$c_{5,12}$} & & \textcolor{BurntOrange}{$c_{5',12}$} & & \\
 \textcolor{Mulberry}{$c_{5,12}$} & & \textcolor{Mulberry}{$c_{5,12}$} & & \textcolor{ForestGreen}{$c_{4,11}$} & & \textcolor{Mahogany}{$c_{6',11}$} & \\
 & \textcolor{Mulberry}{$c_{5,12}$} & & \textcolor{RoyalPurple}{$c_{3,10}$} & & \textcolor{BurntOrange}{$c_{5',10}$} & & \textcolor{Mulberry}{$c_{5,8}$} \\
 \textcolor{Mulberry}{$c_{5,12}$} & & \textcolor{WildStrawberry}{$c_{2,9}$} & & \textcolor{ForestGreen}{$c_{4,9}$} & & \textcolor{Mulberry}{$c_{5,8}$} & \\
 & \textcolor{CornflowerBlue}{$c_{1,8}$} & & \textcolor{RoyalPurple}{$c_{3,8}$} & & \textcolor{Mulberry}{$c_{5,8}$} & & \textcolor{Mulberry}{$c_{5,8}$} \\
 \textcolor{Mulberry}{$c_{5,4}$}& & \textcolor{WildStrawberry}{$c_{2,7}$} & & \textcolor{ForestGreen}{$c_{4,7}$} & & \textcolor{Mulberry}{$c_{5,8}$} & \\
 & \textcolor{Mulberry}{$c_{5,4}$}& & \textcolor{RoyalPurple}{$c_{3,6}$} & & \textcolor{BurntOrange}{$c_{5',6}$} & & \textcolor{Mulberry}{$c_{5,8}$} \\
 \textcolor{Mulberry}{$c_{5,4}$}& & \textcolor{Mulberry}{$c_{5,4}$}& & \textcolor{ForestGreen}{$c_{4,5}$} & & \textcolor{Mahogany}{$c_{6',5}$} & \\
 & \textcolor{Mulberry}{$c_{5,4}$}& & \textcolor{Mulberry}{$c_{5,4}$} & & \textcolor{BurntOrange}{$c_{5',4}$} & & $k$ \\
 \textcolor{Mulberry}{$c_{5,4}$}& & \textcolor{Mulberry}{$c_{5,4}$}& & \textcolor{ForestGreen}{$c_{4,3}$} & & $k$ & \\
 & \textcolor{Mulberry}{$c_{5,4}$}& & \textcolor{RoyalPurple}{$c_{3,2}$} & & $k$ & & $k$ \\
 \textcolor{Mulberry}{$c_{5,4}$}& & \textcolor{WildStrawberry}{$c_{2,1}$} & & $k$ & & $k$ & \\
 & \textcolor{CornflowerBlue}{$c_{1,0}$} & & $k$ & & $k$ & & $k$ \\
 & & $k$ & & $k$ & & $k$ & \\
 & & & $k$ & & $k$ & & $k$ \\
 & & & & $k$ & & $k$ & \\
 & & & & & $k$ & & $k$ \\
 & & & & & & $k$ & \\
 & & & & & & & $k$
\end{tabular}
\end{center}
If $\relt' \myarrow{i} \selt'$ for arrays in our $J_{6}$-component of $L_{\mytinyE_{7}}(k\omega_{1})$, then we must have $g_{i,j}(\relt') + 1 = g_{i,j}(\selt')$ where $(i,j)$ is not one of the above-identified fixed-value positions.  Then we set
\[\myqX_{\selt',\relt'}^{(i)} := \sqrt{\myqQ_{\relt',\selt'}^{(i)}\, } =: \myqY_{\relt',\selt'}^{(i)}.\]
\end{figure}
\clearpage

\noindent \CoefficientFiguresPartTwo, we supply coefficients to the edges of the $J_{5}$- and $J_{6}$-components of $L_{\mytinyE_{7}}(k\omega_{1})$. 
But, as observed in the third paragraph of this section, this results in twice defining coefficients for edges of color $i \in \{1,2,3,4\}$. 
That the two definitions agree is the content of the next result.

\noindent 
{\bf \AgreeLemma}\ \ {\sl Let $\relt \myarrow{i} \selt$ be an edge in} $L_{\mytinyE_{7}}(k\omega_{1})$ {\sl with color $i \in \{1,2,3,4\}$.  Let $\myqP^{(i)}_{\relt,\selt}$ be the product associated to this edge in \EsevenParallelogramFigureAfivePartTwo, and let $\myqQ^{(i)}_{\relt,\selt}$ denote the product from \EsevenParallelogramFigureAsixPartTwo.  Then $\myqP^{(i)}_{\relt,\selt} = \myqQ^{(i)}_{\relt,\selt}$.}

{\em Proof.} Throughout the proof, `$c_{p,q}$' by itself is short for `$c_{p,q}(\selt)$'. 
To start, suppose $\relt \myarrow{4} \selt$, with $\mbox{\textcolor{ForestGreen}{$c_{4,j'}$}}(\relt)+1=\mbox{\textcolor{ForestGreen}{$c_{4,j'}$}}(\selt)$ and $g_{4,j}(\selt)=\mbox{\textcolor{ForestGreen}{$c_{4,j'}$}}(\selt)$. 
The defining factors for $\myqP^{(\mbox{\scriptsize \em 4})}_{\relt,\selt}$ and $\myqQ^{(\mbox{\scriptsize \em 4})}_{\relt,\selt}$ can only differ where the entries of \CoefficientFigures\ differ on either side of the $4^{\mbox{\scriptsize th}}$ column.  For $\myqP^{(\mbox{\scriptsize \em 4})}_{\relt,\selt}$, these factors are
\[(\mbox{\textcolor{ForestGreen}{$c_{4,j'}$}} - \mbox{\textcolor{Mulberry}{$c_{5,4}$}} + j - 1)(\mbox{\textcolor{ForestGreen}{$c_{4,j'}$}} - \mbox{\textcolor{Mulberry}{$c_{5,12}$}} + j + 3)(\mbox{\textcolor{ForestGreen}{$c_{4,j'}$}} - \mbox{\textcolor{BurntOrange}{$c_{5',4}$}} + j - 0 - 1)(\mbox{\textcolor{ForestGreen}{$c_{4,j'}$}} - \mbox{\textcolor{BurntOrange}{$c_{5',12}$}} + j + 4 - 1),\]
and for $\myqQ^{(\mbox{\scriptsize \em 4})}_{\relt,\selt}$ these factors are
\[(\mbox{\textcolor{ForestGreen}{$c_{4,j'}$}} - \mbox{\textcolor{BurntOrange}{$c_{5',4}$}} + j - 1)(\mbox{\textcolor{ForestGreen}{$c_{4,j'}$}} - \mbox{\textcolor{BurntOrange}{$c_{5',12}$}} + j + 3)(\mbox{\textcolor{ForestGreen}{$c_{4,j'}$}} - \mbox{\textcolor{Mulberry}{$c_{5,4}$}} + j - 0 - 1)(\mbox{\textcolor{ForestGreen}{$c_{4,j'}$}} - \mbox{\textcolor{Mulberry}{$c_{5,12}$}} + j + 4 - 1).\] 
So, $\myqP^{(\mbox{\scriptsize \em 4})}_{\relt,\selt}=\myqQ^{(\mbox{\scriptsize \em 4})}_{\relt,\selt}$. 
Now say $\relt \myarrow{3} \selt$, with $\mbox{\textcolor{RoyalPurple}{$c_{3,j'}$}}(\relt)+1=\mbox{\textcolor{RoyalPurple}{$c_{3,j'}$}}(\selt)$ and $g_{3,j}(\selt)=\mbox{\textcolor{RoyalPurple}{$c_{3,j'}$}}(\selt)$. 
The relevant factors for $\myqP^{(\mbox{\scriptsize \em 3})}_{\relt,\selt}$ are
\[\frac{(\mbox{\textcolor{RoyalPurple}{$c_{3,j'}$}} - \mbox{\textcolor{BurntOrange}{$c_{5',4}$}}+j+0-1)(\mbox{\textcolor{RoyalPurple}{$c_{3,j'}$}} - \mbox{\textcolor{BurntOrange}{$c_{5',4}$}}+j+1-1)(\mbox{\textcolor{RoyalPurple}{$c_{3,j'}$}} - \mbox{\textcolor{BurntOrange}{$c_{5',12}$}}+j+4-1)(\mbox{\textcolor{RoyalPurple}{$c_{3,j'}$}} - \mbox{\textcolor{BurntOrange}{$c_{5',12}$}}+j+5-1)}{(\mbox{\textcolor{RoyalPurple}{$c_{3,j'}$}} - \mbox{\textcolor{BurntOrange}{$c_{5',4}$}}+j+0-1)(\mbox{\textcolor{RoyalPurple}{$c_{3,j'}$}} - \mbox{\textcolor{BurntOrange}{$c_{5',4}$}}+j+0)(\mbox{\textcolor{RoyalPurple}{$c_{3,j'}$}} - \mbox{\textcolor{BurntOrange}{$c_{5',12}$}}+j+4)(\mbox{\textcolor{RoyalPurple}{$c_{3,j'}$}} - \mbox{\textcolor{BurntOrange}{$c_{5',12}$}}+j+4-1)}.\] 
and for $\myqQ^{(\mbox{\scriptsize \em 3})}_{\relt,\selt}$ the relevant factors are
\[\frac{(\mbox{\textcolor{RoyalPurple}{$c_{3,j'}$}} - \mbox{\textcolor{Mulberry}{$c_{5,4}$}}+j+0-1)(\mbox{\textcolor{RoyalPurple}{$c_{3,j'}$}} - \mbox{\textcolor{Mulberry}{$c_{5,4}$}}+j+1-1)(\mbox{\textcolor{RoyalPurple}{$c_{3,j'}$}} - \mbox{\textcolor{Mulberry}{$c_{5,12}$}}+j+4-1)(\mbox{\textcolor{RoyalPurple}{$c_{3,j'}$}} - \mbox{\textcolor{Mulberry}{$c_{5,12}$}}+j+5-1)}{(\mbox{\textcolor{RoyalPurple}{$c_{3,j'}$}} - \mbox{\textcolor{Mulberry}{$c_{5,4}$}}+j+0-1)(\mbox{\textcolor{RoyalPurple}{$c_{3,j'}$}} - \mbox{\textcolor{Mulberry}{$c_{5,4}$}}+j+0)(\mbox{\textcolor{RoyalPurple}{$c_{3,j'}$}} - \mbox{\textcolor{Mulberry}{$c_{5,12}$}}+j+4)(\mbox{\textcolor{RoyalPurple}{$c_{3,j'}$}} - \mbox{\textcolor{Mulberry}{$c_{5,12}$}}+j+4-1)}.\] 
Both expressions simplify to $1$, so $\myqP^{(\mbox{\scriptsize \em 3})}_{\relt,\selt}=\myqQ^{(\mbox{\scriptsize \em 3})}_{\relt,\selt}$. 
The argument that $\myqP^{(\mbox{\scriptsize \em i})}_{\relt,\selt}=\myqQ^{(\mbox{\scriptsize \em i})}_{\relt,\selt}$ when $i \in \{1,2\}$ is entirely similar to the $i = 3$ case.\hfill\QED 

The main result of the paper, together with \ConstructionCorollaries, is:

\noindent
{\bf \ConstructionTheorem}\ \ 
{\sl Let} $L := L_{\mytinyE_{7}}(k\omega_{1})$ {\sl with $I_{7} = \{1,2,3,4,5,5',6'\}$. 
For each $i \in I$, assign the scalar pairs $\{(\myqX_{\telt,\selt}^{(i)},\myqY_{\selt,\telt}^{(i)})\}_{\selt \myarrow{i} \telt \mbox{\scriptsize \ in } L}$ prescribed in \CoefficientFigures\ to the color $i$ edges of $L$. 
Then $L$ is} $\myE_{7}${\sl -structured and the scalars satisfy all diamond and crossing relations. 
Therefore the action of the generators of} $\mathfrak{g}(\myE_{7})$ {\sl on the vector space $V[L]$ as defined by the formulas} (1) {\sl in} \S \Setup\ {\sl is well-defined; with $\mym_{i}(\relt)$ as identified in \JLemmas, we have $\myqh_{i}.v_{\relt} = \mym_{i}(\relt)v_{\relt}$ for all $\relt \in L$ and $i \in I_{7}$, and in particular $\{v_{\relt}\}_{\relt \in L}$ is a weight basis for the} $\mathfrak{g}(\myE_{7})${\sl -module $V[L]$; and the lattice $L$ together with $\displaystyle \bigcup_{i \in I}\{(\myqX_{\telt,\selt}^{(i)},\myqY_{\selt,\telt}^{(i)})\}_{\selt \myarrow{i} \telt \mbox{\scriptsize \ in } L}$ is its representation diagram. 
Moreover, the} $\mathfrak{g}(\myE_{7})${\sl -module $V[L]$ is irreducible with highest weight $k\omega_{1}$, and} $\WGF(L;\myvarZ) = \mychar(V[L]) = \chi_{_{k\omega_{1}}}^{\mytinyE_{7}}$.  

{\em Proof.} The last sentence of the theorem statement follows directly from \EsixAndEsevenCharacter, as does the claim that $L$ is $\myE_{7}$-structured. 
The remaining claims will follow from \CombinatorialSerre\ once we confirm that $L$ together with the said assignment of edge coefficients is a representation diagram.  
As observed in the third paragraph of this section, this latter task only requires us to check that on any edge in $L$ with color $i \in \{1,2,3,4\}$, the coefficients supplied by \CoefficientFigures\ agree -- which has been established in \AgreeLemma\ -- and that diamond relations hold on diamonds of the form \parbox{1.4cm}{\begin{center}
\setlength{\unitlength}{0.2cm}
\begin{picture}(5,3)
\put(2.5,0){\circle*{0.5}} \put(0.5,2){\circle*{0.5}}
\put(2.5,4){\circle*{0.5}} \put(4.5,2){\circle*{0.5}}
\put(0.5,2){\line(1,1){2}} \put(2.5,0){\line(-1,1){2}}
\put(4.5,2){\line(-1,1){2}} \put(2.5,0){\line(1,1){2}}
\put(1.05,0.55){\em \small 5$\ \!\! '$} \put(3.05,0.55){\em \small 5}
\put(1.05,2.55){\em \small 5} \put(3.05,2.55){\em \small 5$\ \!\! '$}
\put(3,-0.75){\footnotesize $\relt$} \put(5.25,1.75){\footnotesize $\telt$}
\put(3,4){\footnotesize $\uelt$} \put(-1,1.75){\footnotesize $\selt$}
\end{picture} \end{center}} and \parbox{1.4cm}{\begin{center}
\setlength{\unitlength}{0.2cm}
\begin{picture}(5,3)
\put(2.5,0){\circle*{0.5}} \put(0.5,2){\circle*{0.5}}
\put(2.5,4){\circle*{0.5}} \put(4.5,2){\circle*{0.5}}
\put(0.5,2){\line(1,1){2}} \put(2.5,0){\line(-1,1){2}}
\put(4.5,2){\line(-1,1){2}} \put(2.5,0){\line(1,1){2}}
\put(1.05,0.55){\em \small 6$\ \!\! '$} \put(3.05,0.55){\em \small 5}
\put(1.05,2.55){\em \small 5} \put(3.05,2.55){\em \small 6$\ \!\! '$}
\put(3,-0.75){\footnotesize $\relt'$} \put(5.25,1.75){\footnotesize $\telt'$}
\put(3,4){\footnotesize $\uelt'$} \put(-1,1.75){\footnotesize $\selt'$}
\end{picture} \end{center}}. 
In fact, for the latter diamond we will show that $\myqP^{(\mbox{\scriptsize \em 5})}_{\relt',\telt'} = \myqP^{(\mbox{\scriptsize \em 5})}_{\selt',\uelt'}$ and $\myqQ^{(\mbox{\scriptsize \em 6$'$})}_{\relt',\selt'} = \myqQ^{(\mbox{\scriptsize \em 6$'$})}_{\telt',\uelt'}$, and similarly for the former diamond.

On the color-$5$-and-$6'$ diamonds, exactly one of the following must be true: ({\sl i}) $\mbox{\textcolor{Mahogany}{$c_{6',5}$}}(\relt')+1=\mbox{\textcolor{Mahogany}{$c_{6',5}$}}(\selt') = \mbox{\textcolor{Mahogany}{$c_{6',5}$}}(\telt')+1=\mbox{\textcolor{Mahogany}{$c_{6',5}$}}(\uelt')$ or ({\sl ii})  $\mbox{\textcolor{Mahogany}{$c_{6',11}$}}(\relt')+1=\mbox{\textcolor{Mahogany}{$c_{6',11}$}}(\selt') = \mbox{\textcolor{Mahogany}{$c_{6',11}$}}(\telt')+1=\mbox{\textcolor{Mahogany}{$c_{6',11}$}}(\uelt')$. 
In case ({\sl i}), if we have $\mbox{\textcolor{Mulberry}{$c_{5,4}$}}(\relt')+1=\mbox{\textcolor{Mulberry}{$c_{5,4}$}}(\telt')=\mbox{\textcolor{Mulberry}{$c_{5,4}$}}(\selt')+1=\mbox{\textcolor{Mulberry}{$c_{5,4}$}}(\uelt')$ or $\mbox{\textcolor{Mulberry}{$c_{5,12}$}}(\relt')+1=\mbox{\textcolor{Mulberry}{$c_{5,12}$}}(\telt')=\mbox{\textcolor{Mulberry}{$c_{5,12}$}}(\selt')+1=\mbox{\textcolor{Mulberry}{$c_{5,12}$}}(\uelt')$, then we get $\myqQ^{(\mbox{\scriptsize \em 6$'$})}_{\relt',\selt'} = \myqQ^{(\mbox{\scriptsize \em 6$'$})}_{\telt',\uelt'}$, as there is no interdependence of the quantities used to define the left-hand and right-hand quantities;  similarly, we get $\myqP^{(\mbox{\scriptsize \em 5})}_{\relt',\telt'} = \myqP^{(\mbox{\scriptsize \em 5})}_{\selt',\uelt'}$. 
So now suppose that $\mbox{\textcolor{Mulberry}{$c_{5,8}$}}(\relt')+1=\mbox{\textcolor{Mulberry}{$c_{5,8}$}}(\telt')=\mbox{\textcolor{Mulberry}{$c_{5,8}$}}(\selt')+1=\mbox{\textcolor{Mulberry}{$c_{5,8}$}}(\uelt')$. 
It is evident again that $\myqP^{(\mbox{\scriptsize \em 5})}_{\relt',\telt'} = \myqP^{(\mbox{\scriptsize \em 5})}_{\selt',\uelt'}$, as neither computation depends on any `$\mbox{\textcolor{Mahogany}{$c_{6',5}$}}$' value. 
For each of the quantities $\myqQ^{(\mbox{\scriptsize \em 6$'$})}_{\relt',\selt'}$ and  $\myqQ^{(\mbox{\scriptsize \em 6$'$})}_{\telt',\uelt'}$, we will only consider those factors that include one of $\mbox{\textcolor{Mulberry}{$c_{5,8}$}}(\selt')$ or $\mbox{\textcolor{Mulberry}{$c_{5,8}$}}(\uelt')$. 
For brevity, in what follows both of `$\mbox{\textcolor{Mulberry}{$c_{5,8}$}}$' and `$\mbox{\textcolor{Mahogany}{$c_{6',5}$}}$' refer to $\uelt'$.  
In $\myqQ^{(\mbox{\scriptsize \em 6$'$})}_{\telt',\uelt'}$, we have 
\[\frac{(\mbox{\textcolor{Mahogany}{$c_{6',5}$}} - \mbox{\textcolor{Mulberry}{$c_{5,8}$}} + 1 - 1)(\mbox{\textcolor{Mahogany}{$c_{6',5}$}} - \mbox{\textcolor{Mulberry}{$c_{5,8}$}} + 1 - 0)(\mbox{\textcolor{Mahogany}{$c_{6',5}$}} - \mbox{\textcolor{Mulberry}{$c_{5,8}$}} + 1 - (-1))(\mbox{\textcolor{Mahogany}{$c_{6',5}$}} - \mbox{\textcolor{Mulberry}{$c_{5,8}$}} + 1 - (-1) - 1)}{(\mbox{\textcolor{Mahogany}{$c_{6',5}$}} - \mbox{\textcolor{Mulberry}{$c_{5,8}$}} + 1 - 0 - 1)(\mbox{\textcolor{Mahogany}{$c_{6',5}$}} - \mbox{\textcolor{Mulberry}{$c_{5,8}$}} + 1 - (-1) - 1)(\mbox{\textcolor{Mahogany}{$c_{6',5}$}} - \mbox{\textcolor{Mulberry}{$c_{5,8}$}} + 1 - 0)(\mbox{\textcolor{Mahogany}{$c_{6',5}$}} - \mbox{\textcolor{Mulberry}{$c_{5,8}$}} + 1 - (-1))},\]
which simplifies to $1$.  
Now, $\mbox{\textcolor{Mahogany}{$c_{6',5}$}}(\selt') - \mbox{\textcolor{Mulberry}{$c_{5,8}$}}(\selt') = \mbox{\textcolor{Mahogany}{$c_{6',5}$}}(\uelt') - (\mbox{\textcolor{Mulberry}{$c_{5,8}$}}(\uelt')-1) = \mbox{\textcolor{Mahogany}{$c_{6',5}$}} - \mbox{\textcolor{Mulberry}{$c_{5,8}$}} + 1$, so in $\myqQ^{(\mbox{\scriptsize \em 6$'$})}_{\relt',\selt'}$, we get  
\[\frac{(\mbox{\textcolor{Mahogany}{$c_{6',5}$}} - \mbox{\textcolor{Mulberry}{$c_{5,8}$}}+1)(\mbox{\textcolor{Mahogany}{$c_{6',5}$}} - \mbox{\textcolor{Mulberry}{$c_{5,8}$}}+2)(\mbox{\textcolor{Mahogany}{$c_{6',5}$}} - \mbox{\textcolor{Mulberry}{$c_{5,8}$}}+3)(\mbox{\textcolor{Mahogany}{$c_{6',5}$}} - \mbox{\textcolor{Mulberry}{$c_{5,8}$}}+2)}{(\mbox{\textcolor{Mahogany}{$c_{6',5}$}} - \mbox{\textcolor{Mulberry}{$c_{5,8}$}}+1)(\mbox{\textcolor{Mahogany}{$c_{6',5}$}} - \mbox{\textcolor{Mulberry}{$c_{5,8}$}}+2)(\mbox{\textcolor{Mahogany}{$c_{6',5}$}} - \mbox{\textcolor{Mulberry}{$c_{5,8}$}}+2)(\mbox{\textcolor{Mahogany}{$c_{6',5}$}} - \mbox{\textcolor{Mulberry}{$c_{5,8}$}}+3)},\] 
which also simplifies to $1$. 
It follows that $\myqQ^{(\mbox{\scriptsize \em 6$'$})}_{\relt',\selt'} = \myqQ^{(\mbox{\scriptsize \em 6$'$})}_{\telt',\uelt'}$. 
This completes our analysis of case ({\sl i}). 

For case ({\sl ii}), as with case ({\sl i}), demonstration of the desired inequalities easily reduces to the hypothesis that $\mbox{\textcolor{Mulberry}{$c_{5,8}$}}(\relt')+1=\mbox{\textcolor{Mulberry}{$c_{5,8}$}}(\telt')=\mbox{\textcolor{Mulberry}{$c_{5,8}$}}(\selt')+1=\mbox{\textcolor{Mulberry}{$c_{5,8}$}}(\uelt')$. 
As before, it is apparent that $\myqP^{(\mbox{\scriptsize \em 5})}_{\relt',\telt'} = \myqP^{(\mbox{\scriptsize \em 5})}_{\selt',\uelt'}$. 
Within the products defining $\myqQ^{(\mbox{\scriptsize \em 6$'$})}_{\relt',\selt'}$ and  $\myqQ^{(\mbox{\scriptsize \em 6$'$})}_{\telt',\uelt'}$, we only consider those factors that include one of $\mbox{\textcolor{Mulberry}{$c_{5,8}$}}(\selt')$ or $\mbox{\textcolor{Mulberry}{$c_{5,8}$}}(\uelt')$. 
For the remainder of the paragraph, `$\mbox{\textcolor{Mulberry}{$c_{5,8}$}}$' and `$\mbox{\textcolor{Mahogany}{$c_{6',11}$}}$' refer to $\uelt'$. 
Then $\myqQ^{(\mbox{\scriptsize \em 6$'$})}_{\telt',\uelt'}$, we have 
\[\frac{(\mbox{\textcolor{Mahogany}{$c_{6',11}$}} - \mbox{\textcolor{Mulberry}{$c_{5,8}$}} - 2 - 1)(\mbox{\textcolor{Mahogany}{$c_{6',11}$}} - \mbox{\textcolor{Mulberry}{$c_{5,8}$}} - 2 - 0)(\mbox{\textcolor{Mahogany}{$c_{6',11}$}} - \mbox{\textcolor{Mulberry}{$c_{5,8}$}} - 2 - (-1))(\mbox{\textcolor{Mahogany}{$c_{6',11}$}} - \mbox{\textcolor{Mulberry}{$c_{5,8}$}} - 2 - (-1) - 1)}{(\mbox{\textcolor{Mahogany}{$c_{6',11}$}} - \mbox{\textcolor{Mulberry}{$c_{5,8}$}} - 2 - 0 - 1)(\mbox{\textcolor{Mahogany}{$c_{6',11}$}} - \mbox{\textcolor{Mulberry}{$c_{5,8}$}} - 2 - (-1) - 1)(\mbox{\textcolor{Mahogany}{$c_{6',11}$}} - \mbox{\textcolor{Mulberry}{$c_{5,8}$}} - 2 - 0)(\mbox{\textcolor{Mahogany}{$c_{6',11}$}} - \mbox{\textcolor{Mulberry}{$c_{5,8}$}} - 2 - (-1))},\]
which simplifies to $1$. 
Now, $\mbox{\textcolor{Mahogany}{$c_{6',11}$}}(\selt') - \mbox{\textcolor{Mulberry}{$c_{5,8}$}}(\selt') = \mbox{\textcolor{Mahogany}{$c_{6',11}$}}(\uelt') - (\mbox{\textcolor{Mulberry}{$c_{5,8}$}}(\uelt')-1) = \mbox{\textcolor{Mahogany}{$c_{6',11}$}} - \mbox{\textcolor{Mulberry}{$c_{5,8}$}} + 1$, so in $\myqQ^{(\mbox{\scriptsize \em 6$'$})}_{\relt',\selt'}$, we get  
\[\frac{(\mbox{\textcolor{Mahogany}{$c_{6',11}$}} - \mbox{\textcolor{Mulberry}{$c_{5,8}$}}-2)(\mbox{\textcolor{Mahogany}{$c_{6',11}$}} - \mbox{\textcolor{Mulberry}{$c_{5,8}$}}-1)(\mbox{\textcolor{Mahogany}{$c_{6',11}$}} - \mbox{\textcolor{Mulberry}{$c_{5,8}$}})(\mbox{\textcolor{Mahogany}{$c_{6',11}$}} - \mbox{\textcolor{Mulberry}{$c_{5,8}$}}-1)}{(\mbox{\textcolor{Mahogany}{$c_{6',11}$}} - \mbox{\textcolor{Mulberry}{$c_{5,8}$}}-2)(\mbox{\textcolor{Mahogany}{$c_{6',11}$}} - \mbox{\textcolor{Mulberry}{$c_{5,8}$}}-1)(\mbox{\textcolor{Mahogany}{$c_{6',11}$}} - \mbox{\textcolor{Mulberry}{$c_{5,8}$}}-1)(\mbox{\textcolor{Mahogany}{$c_{6',11}$}} - \mbox{\textcolor{Mulberry}{$c_{5,8}$}})},\] 
which also simplifies to $1$. 
So $\myqQ^{(\mbox{\scriptsize \em 6$'$})}_{\relt',\selt'} = \myqQ^{(\mbox{\scriptsize \em 6$'$})}_{\telt',\uelt'}$, completing our analysis of case ({\sl ii}).

On diamonds of the form \parbox{1.4cm}{\begin{center}
\setlength{\unitlength}{0.2cm}
\begin{picture}(5,3)
\put(2.5,0){\circle*{0.5}} \put(0.5,2){\circle*{0.5}}
\put(2.5,4){\circle*{0.5}} \put(4.5,2){\circle*{0.5}}
\put(0.5,2){\line(1,1){2}} \put(2.5,0){\line(-1,1){2}}
\put(4.5,2){\line(-1,1){2}} \put(2.5,0){\line(1,1){2}}
\put(1.05,0.55){\em \small 5$\ \!\! '$} \put(3.05,0.55){\em \small 5}
\put(1.05,2.55){\em \small 5} \put(3.05,2.55){\em \small 5$\ \!\! '$}
\put(3,-0.75){\footnotesize $\relt$} \put(5.25,1.75){\footnotesize $\telt$}
\put(3,4){\footnotesize $\uelt$} \put(-1,1.75){\footnotesize $\selt$}
\end{picture} \end{center}}, similar case analysis shows that $\myqP^{(\mbox{\scriptsize \em 5})}_{\relt,\telt} = \myqP^{(\mbox{\scriptsize \em 5})}_{\selt,\uelt}$ and 
$\myqQ^{(\mbox{\scriptsize \em 5$'$})}_{\relt,\selt} = \myqQ^{(\mbox{\scriptsize \em 5$'$})}_{\telt,\uelt}$. 
This accounts for all diamond relations and completes the proof.\hfill\QED 

Given how we have identified our $\myE_{6}$-polyminuscule lattices as components of $\myE_{7}$-polyminuscule lattices, the next result follows immediately from \ConstructionTheorem.

\noindent
{\bf \ConstructionCorollaries}\ \ 
{\sl Let $a+b=k$ for nonnegative integers $a$, $b$, and $k$, and let $M$ be any of the} $\myE_{6}$- {\sl polyminuscule lattices} $L_{\mytinyE_{6}}(k\omega_{1'})$, $L_{\mytinyE_{6}}(k\omega_{6'})$, {\sl or} $L_{\mytinyE_{6}}(a\omega_{1'}+b\omega_{6'})$ {\sl as colored by $I_{6} = \{1',2',3',4',5',6'\}$. 
Regard $M$ to be a $\psi(I_{6})$-component of} $L_{k} := L_{\mytinyE_{7}}(k\omega_{1}),$ {\sl and for each $i \in I_{6}$ assign to any edge $\selt \myarrow{i} \telt$ in $M$ the scalar pair $(\myqX_{\telt,\selt}^{(\psi(i))},\myqY_{\selt,\telt}^{(\psi(i))})$ of the corresponding edge in $L_{k}$. 
Then $M$ is} $\myE_{6}${\sl -structured and the scalars satisfy all diamond and crossing relations. 
Therefore the action of the generators of} $\mathfrak{g}(\myE_{6})$ {\sl on the vector space $V[M]$ as defined by the formulas} (1) {\sl in} \S \Setup\ {\sl is well-defined; $\{v_{\relt}\}_{\relt \in M}$ is a weight basis for the} $\mathfrak{g}(\myE_{6})${\sl -module $V[M]$; and the lattice $M$ together with the assigned edge coefficients is its representation diagram. 
Moreover, the} $\mathfrak{g}(\myE_{6})${\sl -module $V[M]$ is irreducible, and its highest weight and character are as identified in \EsixAndEsevenCharacter.}

\vspace*{0.1in} 
\noindent {\bf \S \CombinatorialDiscussion\ Some further considerations.} 
To close the paper, we mention certain combinatorial distinctions enjoyed by our $\myE_{6}$- and $\myE_{7}$-polyminuscule lattices. 
(These combinatorial distinctions are, in fact, shared by all  supporting graphs for irreducible representations of semisimple Lie algebras.) 
We also pose some `extremal' questions relating to our constructions. 
Then we briefly discuss some generalizations of ideas from this paper being pursued in \cite{DD2}. 

Let $q$ be an indeterminate, and for any positive integer $m$, let $[m]$ denote the $q$-integer $\frac{1-q^{m}}{1-q}$. 
A degree $\ell$ polynomial $\sum_{i=0}^{\ell}c_{i}q^{i}$ is {\em symmetric} if $c_{\ell - i}=c_{i}$ for all $i \in \{0,1,\ldots,\ell\}$ and {\em unimodal} if there is some $u \in \{0,1,\ldots,\ell\}$ such that $c_{0} \leq c_{1} \cdots \leq c_{u} \geq \cdots \geq c_{\ell -1} \geq c_{\ell}$. 
A ranked poset $R$ with rank function $\rho$ is {\em rank symmetric} (respectively, {\em rank unimodal}) if its rank generating function $\RGF(R;q)$ is symmetric (respectively, unimodal). 
A {\em rank} of $R$ is a subset $\rho^{-1}(i)$, where $i \in \{0,1,\ldots,\ell\}$. 
We say $R$ is {\em strongly Sperner} if, for all positive integers $m$, no union of any $m$ antichains in $R$ is larger than the union of the $m$ largest distinct ranks of $R$. 
What follows is a simple application of Proposition 4.1 of \cite{DD1}, a result whose symmetry and unimodality aspects are ultimately due to Dynkin \cite{Dynkin} (see also the 3$^{\mbox{\tiny rd}}$ page of the commentary \cite{KOV}), whose `quotient-of-product' rank generating function identities can be traced to Jacobson \cite{Jac} via Lepowsky \cite{Lepowsky} and Proctor \cite{PrEur}, and whose Sperner aspects are due to Proctor \cite{PrPeck} and Stanley \cite{StanHard}. 

\noindent 
{\bf \CombinatorialResults}\ \ {\sl Let $L$ be any of the} $\myE_{6}$- {\sl and} $\myE_{7}$- {\sl polyminuscule lattices studied in \S\S \ESection\ and \ConstructionSection. 
Then $L$ is rank symmetric, rank unimodal, and strongly Sperner.  
Moreover, the rank generating  functions of these polyminuscule lattices can be expressed as a quotient of products as follows:}
{\small \begin{eqnarray*}
\RGF(L_{\mytinyE_{7}}(k\omega_{1});q) & = & \mbox{\footnotesize $\displaystyle \frac{[k+17]}{[17]}\cdot\frac{[k+16]}{[16]}\cdot\frac{[k+15]}{[15]}\cdot\frac{[k+14]}{[14]}\cdot\left(\frac{[k+13]}{[13]}\right)^{2}\cdot\left(\frac{[k+12]}{[12]}\right)^{2}\cdot\left(\frac{[k+11]}{[11]}\right)^{2}$}\\
&  & \cdot\ \mbox{\footnotesize $\displaystyle \left(\frac{[k+10]}{[10]}\right)^{2}\cdot\left(\frac{[k+9]}{[9]}\right)^{3}\cdot\left(\frac{[k+8]}{[8]}\right)^{2}\cdot\left(\frac{[k+7]}{[7]}\right)^{2}\cdot\left(\frac{[k+6]}{[6]}\right)^{2}\cdot\left(\frac{[k+5]}{[5]}\right)^{2}$}\\
&  & \cdot\ \mbox{\footnotesize $\displaystyle \frac{[k+4]}{[4]}\cdot\frac{[k+3]}{[3]}\cdot\frac{[k+2]}{[2]}\cdot\frac{[k+1]}{[1]}$}
\end{eqnarray*}
\begin{eqnarray*}
\RGF(L_{\mytinyE_{6}}(a\omega_{1}+b\omega_{6});q) & = & \mbox{\footnotesize $\displaystyle \frac{[a+7]}{[7]}\cdot\frac{[a+6]}{[6]}\cdot\frac{[a+5]}{[5]}\cdot\left(\frac{[a+4]}{[4]}\right)^{2}\cdot\frac{[a+3]}{[3]}\cdot\frac{[a+2]}{[2]}\cdot\frac{[a+1]}{[1]}\cdot\frac{[a+b+11]}{[11]}$}\\
&  & \cdot\ \mbox{\footnotesize $\displaystyle \frac{[a+b+10]}{[10]}\cdot\frac{[a+b+9]}{[9]}\cdot\left(\frac{[a+b+8]}{[8]}\right)^{2}\cdot\frac{[a+b+7]}{[7]}\cdot\frac{[a+b+6]}{[6]}\cdot\frac{[a+b+5]}{[5]}$}\\
&  & \cdot\ \mbox{\footnotesize $\displaystyle \frac{[b+7]}{[7]}\cdot\frac{[b+6]}{[6]}\cdot\frac{[b+5]}{[5]}\cdot\left(\frac{[b+4]}{[4]}\right)^{2}\cdot\frac{[b+3]}{[3]}\cdot\frac{[b+2]}{[2]}\cdot\frac{[b+1]}{[1]}$}
\end{eqnarray*}}

{\em Proof.} The claims of the second sentence of the theorem statement follow from Proposition 4.1 of \cite{DD1}. 
The expressions given for the rank generating functions are obtained by taking a concrete realization of the $\myE_{6}$ and $\myE_{7}$ root systems (e.g.\ \cite{Hum} or \cite{BMP}) and applying them to the quotient-of-products formula given in Proposition 4.1 of \cite{DD1}.\hfill\QED

Many of the weight bases that can be realized by combinatorial methods similar to those employed here enjoy certain of the `extremal' properties first studied in \cite{DonSupp}. 
A supporting graph, or its attendant weight basis, is {\em edge-minimal} if it contains no other supporting graph for the same representation as a proper subgraph. 
It is {\em solitary} if the only other weight bases with the same supporting graph are those obtained simply by re-scaling each of the vectors of the given weight basis. 
Some weight bases that are both edge-minimal and solitary include: The weight bases for the fundamental representations of $\mathfrak{g}(\myB_{n})$ and $\mathfrak{g}(\myC_{n})$ that are supported by the diamond-colored distributive lattices of \cite{DonSupp} and \cite{Beck}, the GT bases for the irreducible representations of the special linear Lie algebras (see \cite{DonSupp}), the weight bases for the `one-rowed' representations of $\mathfrak{g}(\myB_{n})$ and $\mathfrak{g}(\myG_{2})$ from \cite{DLP2}, and the weight bases of the `spin-node' representations of $\mathfrak{g}(\myB_{n})$ and $\mathfrak{g}(\myD_{n})$ from \cite{DD1}. 
The edge-minimal and solitary properties often occur together, and it is possible that these properties are, with some level of generality, equivalent.  
An interesting open question, at least to us, is whether the $\myE_{6}$- and $\myE_{7}$-polyminuscule lattices presented here enjoy either/both of these extremal properties. 

The $\myE_{6}$- and $\myE_{7}$-polyminuscule lattices of this paper are very special cases of the more general polyminuscule lattices to be presented in \cite{DD2}. 
These more general lattices include, but are not limited to, the $\myE_{6}$ and $\myE_{7}$ lattices of this paper, all skew-tabular lattices, and the even and odd orthogonal lattices of \cite{DD1}. 
At this time, we know how to use our general polyminuscule lattices to define, for all irreducible root systems and in a completely uniform way, analogs of the skew Schur functions. 
We have also developed a general rule for decomposing, as a sum of Weyl bialternants, the product of any Weyl bialternant with any of our polyminuscule symmetric functions. 
Evidence from our examination of many cases suggests it might be possible generally to construct representations on all polyminuscule lattices in a manner similar to the $\mathfrak{g}(\myE_{6})$- and  $\mathfrak{g}(\myE_{7})$-representation constructions obtained here.

%
\renewcommand{\refname}{\normalsize \bf References}
\renewcommand{\baselinestretch}{1.1}

\end{document}